\documentclass[12pt]{article}
\usepackage{latexsym,amsfonts,amsmath,graphics,amssymb}
\usepackage{epsfig}
\usepackage{verbatim}
\usepackage{stix}

\usepackage{pgfplots}

\DeclareMathOperator{\spann}{span}
\DeclareMathOperator{\wav}{wav}
\DeclareMathOperator{\supp}{supp}
\DeclareMathOperator{\wor}{wor}
\DeclareMathOperator{\scr}{sc}
\DeclareMathOperator{\rms}{rms}
\DeclareMathOperator{\dom}{dom}

\newtheorem{theorem}{Theorem}
\newtheorem{lemma}{Lemma}
\newtheorem{corollary}{Corollary}
\newtheorem{proposition}{Proposition}

\newtheorem{definition}{Definition}

\newtheorem{remark}{Remark}
\newenvironment{proof}{\begin{trivlist}
    \item[\hskip\labelsep{\it Proof.}]}{$\hfill\Box$\end{trivlist}}

\newcommand{\bsa}{\boldsymbol{a}}

\newcommand{\bsk}{\boldsymbol{k}}
\newcommand{\bsi}{\boldsymbol{i}}
\newcommand{\bsj}{\boldsymbol{j}}

\newcommand{\bsx}{\boldsymbol{x}}

\newcommand{\bst}{\boldsymbol{t}}

\newcommand{\bsp}{\boldsymbol{p}}

\newcommand{\bsy}{\boldsymbol{y}}

\newcommand{\bsalpha}{\boldsymbol{\alpha}}

\newcommand{\bszero}{\boldsymbol{0}}
\newcommand{\rd}{\,\mathrm{d}}
\newcommand{\NN}{\mathbb{N}}
\newcommand{\ZZ}{\mathbb{Z}}

\newcommand{\RR}{\mathbb{R}}
\newcommand{\real}{\mathbb{R}}

\newcommand{\HH}{\mathcal{H}}

\newcommand{\esssup}{\mathop{\mathrm{ess\,sup}}}

\newcommand{\cx}{{\mathcal X}}
\newcommand{\cs}{{\mathcal S}}
\newcommand{\cp}{{\mathcal P}}
\newcommand{\ce}{{\mathcal E}}
\newcommand{\re}{\mathbb{R}}
\newcommand{\cf}{{\mathcal{F}}}
\newcommand{\cfi}{{\mathcal{F}}^{-1}}
\newcommand{\cl}{{\mathcal L}}
\newcommand{\sign}{{\rm sign} \, }
\newcommand{\N}{\mathbb{N}}

\newcommand{\be}{\begin{equation}}
\newcommand{\ee}{\end{equation}}
\newcommand{\beq}{\begin{eqnarray}}
\newcommand{\beqq}{\begin{eqnarray*}}
\newcommand{\eeq}{\end{eqnarray}}
\newcommand{\eeqq}{\end{eqnarray*}}

\newcommand{\zz}{\mathbb{Z}}

\newcommand{\C}{\mathbb{C}}
\allowdisplaybreaks

\newcommand{\mig}[1]{\textcolor{blue}{#1}}
\newcommand{\red}[1]{\textcolor{red}{#1}}

\begin{document}

\title{\scshape 
QMC integration based on arbitrary $(t,m,s)$-nets yields optimal convergence rates on several scales of function spaces}

\author{$\quad$Michael Gnewuch\\$\quad$\\
{\small Institut f\"ur Mathematik,
Universit\"at Osnabr\"uck}\\
{\small Albrechtstra{\ss}e 28 a} \\
{\small 49076 Osnabr\"uck, Germany} \\
{\small email: michael.gnewuch@uni-osnabrueck.de}\\
$\quad$\\
Josef Dick\\$\quad$\\
{\small School of Mathematics and Statistics,
University of New South Wales}\\
{\small Sydney, NSW, 2052, Australia} \\
{\small email: josef.dick@unsw.edu.au}\\
$\quad$\\
Lev  Markhasin\\$\quad$\\
{\small 
Stuttgart Technology Center, Sony Semiconductor Solutions Europe}\\
{\small Bertha-Thalheimer-Weg 6}\\
{\small 70372 Stuttgart, Germany} \\
{\small  email: lev.markhasin@sony.com}\\
$\quad$\\
Winfried  Sickel\\$\quad$\\
{\small 
Mathematisches Institut, Friedrich-Schiller-Universit\"at Jena}\\
{\small Ernst-Abbe-Platz 1-4}\\
{\small 07740 Jena, Germany} \\
{\small  email: winfried.sickel@uni-jena.de}\\
}

\date{\today}
\maketitle

\pagebreak

\begin{abstract}
We study the integration problem over the $s$-dimensional unit cube 
on four types of Banach spaces of integrands. First we consider Haar wavelet spaces, consisting of functions whose Haar wavelet coefficients exhibit a certain decay behavior measured by a parameter $\alpha >0$.  We study the worst case error of integration over the norm unit ball
and provide upper error  bounds for 
quasi-Monte Carlo (QMC)
cubature rules based on arbitrary $(t,m,s)$-nets as well as matching lower error bounds for arbitrary cubature rules. These results show that using arbitrary $(t,m,s)$-nets as sample points yields the best possible rate of convergence. Afterwards we study spaces of integrands of fractional smoothness $\alpha \in (0,1)$ and 
state a sharp Koksma-Hlawka-type inequality. More precisely, we 
show that on those spaces the worst case error of integration is equal to the corresponding fractional discrepancy. 
Those spaces can be continuously embedded into tensor product Bessel potential spaces, also known as Sobolev spaces of dominated mixed smoothness, with the same set of parameters.
The latter spaces can be embedded into suitable Besov spaces of dominating mixed smoothness $\alpha$, which in turn can be embedded into the Haar wavelet spaces with the same set of parameters. 
Therefore our upper error bounds on Haar wavelet spaces
for 
QMC cubatures based on $(t,m,s)$-nets transfer (with possibly different constants) to the corresponding spaces of integrands of fractional smoothness and to Sobolev and Besov spaces of dominating mixed smoothness. 
Moreover, known lower error bounds for periodic Sobolev and Besov spaces of dominating mixed smoothness show that QMC integration based on arbitrary $(t,m,s)$-nets yields the best possible convergence rate on periodic as well as on non-periodic Sobolev and Besov spaces of dominating smoothness.
\end{abstract}

\noindent
{\bf MSC classification numbers:}
 Primary: 65D30, 46E35; 
 Secondary: 11K38, 42B35, 46B28.
 
 \vspace{1ex}
 


\tableofcontents


\section{Introduction}


\subsection{Goal and Motivation}

We consider the integration problem of approximating the $s$-dimensional integral
\begin{equation}\label{eq:I_s}
I_s(f) := \int_{[0,1]^s} f(\bsx) \, {\rm}d \bsx
\end{equation}
with a \emph{quasi-Monte Carlo (QMC) cubature rule}, i.e., an equal weight cubature rule whose weights sum up to one. The classical Koksma-Hlawka inequality 
ensures that the integration error of a QMC cubature rule 
\begin{equation}\label{eq:qmc_cub}
Q_P(f) := \frac{1}{|P|} \sum_{\bsp\in P} f(\bsp),
\end{equation}
where $|P|$ denotes the cardinality of the set of integration points $P$, 
can be bounded by the product of the star discrepancy of $P$ and the Hardy-Krause variation of the integrand $f$. Unfortunately, many integrands that appear in important applications do not exhibit a finite Hardy-Krause variation, see, e.g., \cite{Owen05}.

To overcome this fundamental issue, 
Dick~\cite{D08} studied the  integration problem for a suitable range of parameters $p,q \in [1,\infty]$ on Banach spaces $\mathcal{H}_{\alpha, s, p,q}$ of functions $f: [0,1]^s\to \RR$ of fractional smoothness $\alpha \in (0,1)$, defined with the help of the fractional derivative in the sense of Riemann and Liouville. These spaces contain functions that do not have finite Hardy-Krause variation. 
He introduced a notion of fractional discrepancy $D^*_{\alpha,s,p,q}$ and derived a Koksma-Hlawka-type inequality of fractional order. It ensures that the integration error 
of $Q_P$ 
can be bounded by the product of a fractional discrepancy of $P$ and the norm of the integrand in
$\mathcal{H}_{\alpha, s, p,q}$:
$$ |I_s(f) - Q_P(f)| \le  D^*_{\alpha,s,p',q'}(P) \|f\|_{\alpha,s,p,q}
\hspace{3ex}\text{for all $f\in \mathcal{H}_{\alpha, s, p,q}$\,,}$$
where $p'$ and $q'$ denote the adjoint exponents of $p$ and $q$, respectively, i.e.,
$$
\frac{1}{p} + \frac{1}{p'} = 1 = \frac{1}{q} + \frac{1}{q'}\,.
$$

Nevertheless, an essential question still remained open: 
\vspace{1.5ex}

\emph{What kind of integration points actually exhibit a small fractional discrepancy and will therefore lead to a small integration error on the spaces
$\mathcal{H}_{\alpha, s, p,q}$ of fractional smoothness?}
\vspace{1.5ex}

Here the error criterion we are interested in is the worst case integration error over the norm unit ball of the function space considered.

This important question was the starting point of our work.
Our conjecture was that $(t,m,s)$-nets, which are QMC point sets with very small star discrepancy, will also exhibit a very small fractional discrepancy and a small integration error.
Since the fractional discrepancy is more complicated to analyze than the common ($L_{p}$-)star discrepancy, we choose 
the following approach to prove our conjecture. It essentially consists of three main steps, each one leading to results that are interesting in their own right: 
\begin{enumerate}
\item First we prove matching upper and lower error bounds that show that \emph{arbitrary} $(t,m,s)$-nets lead to optimal QMC cubature rules on Haar wavelet  spaces $\mathcal{H}_{\wav, \alpha, s, p,q}$; here the parameter $\alpha$ characterizes the decay of the Haar wavelet coefficients of functions in $\mathcal{H}_{\wav, \alpha, s, p,q}$. 
\item Then we establish that 
the spaces of fractional smoothness can be continuously embedded into suitable Haar wavelet spaces
by taking an ``embedding detour'' over Sobolev spaces $S^\alpha_p H([0,1]^s)$ and Besov spaces $S^\alpha_{p,q}B([0,1]^s)$ of dominating mixed smoothness~$\alpha$.
In this way we are able to show that our upper bounds for the worst case error of integration on  
the Haar wavelet spaces transfer  
(with possibly different constants) 
to the Besov and Sobolev spaces of dominated mixed smoothness and to our spaces of fractional smoothness.
\item By showing  that the Koksma-Hlawka-type inequality of fractional order from \cite{D08} is actually sharp, i.e., that the worst case error of a QMC cubature rule $Q_P$ over $\mathcal{H}_{\alpha, s, p,q}$ is exactly the fractional discrepancy $D^*_{\alpha,s,p',q'}(P)$ of the underlying set $P$ of integration points,
we finally can verify that arbitrary $(t,m,s)$-nets exhibit a very small fractional discrepancy.
\end{enumerate}

\subsection{Main Results}

Let us explain the main results and the organization of our paper in more detail:
In Section~\ref{Sec:Haar} we discuss Haar wavelet spaces $\mathcal{H}_{\wav, \alpha, s, p,q}$ for parameters $\alpha >0$ and $p,q\in [1,\infty]$, cf. Definition~\ref{Def:Wavelet_Space_Multi}. Similarly as in Section~\ref{Sec:Fractional} and~\ref{Sec:Besov}, where we discuss spaces of functions of fractional smoothness and Besov and Sobolev spaces, respectively, we first consider the case of univariate functions, see Section~\ref{HAAR-1},
and afterwards the case of multivariate functions, see Section~\ref{HAAR-s}.
Let us just highlight one of the new results stated in Subsections~\ref{HAAR-1} and~\ref{HAAR-s}: In Theorem~\ref{tensor17} 
we show that the Haar wavelet spaces of $s$-variate functions exhibit a tensor product structure, which for $p=1$ is induced by the projective tensor norm and for  
$1<p<\infty$ by the $p$-nuclear tensor norm, which is also known as left $p$-Chevet-Saphar norm.
For the convenience of the reader, we added an appendix where we summarized the most important definitions and facts about tensor products of Banach spaces and where all proofs of tensor product results can be found, see Section~\ref{Sec:Appendix}.

In Section~\ref{Subsec:QMC} we study the approximation of integrals $I_s(f)$ as in \eqref{eq:I_s} with integrands $f$ from $\mathcal{H}_{\wav, \alpha, s, p,q}$ by QMC cubature rules $Q_P(f)$ as in \eqref{eq:qmc_cub}, where $P$ is an arbitrary $(t,m,s)$-net in some base $b\ge 2$. We show that the \emph{worst-case error} of integration, defined by
$$e^{\wor}(Q_P, \mathcal{H}_{\wav, \alpha, s, p,q}) := \sup |I_s(f) -Q_P(f)|,$$
where the supremum is taken over all functions from the norm unit ball of $\mathcal{H}_{\wav, \alpha, s, p,q}$, satisfies
\begin{equation}\label{eq:upper_bound_Haar}
e^{\wor}(Q_P, \mathcal{H}_{\wav, \alpha, s, p,q})  = O \left( N^{-\alpha} \ln(N)^{\frac{s-1}{q'}} \right),
\end{equation}   
where $N := |P| = b^m$, $p,q\in [0,\infty]$ and $\alpha \ge 1/p$ if $q=1$ and $\alpha > 1/p$ if $q>1$,
see Theorem~\ref{Nets}. 
For fixed parameters $s, \alpha, p,q$ the implicit constant in the big-$O$-notation depends on $t$ and $b$, but not on $m$ and also not on the specific construction of the $(t,m,s)$-net $P$. Note furthermore, that for $q=1$ we have $q'= \infty$, and consequently
\begin{equation*}
e^{\wor}(Q_P, \mathcal{H}_{\wav, \alpha, s, p,1})  = O \left( N^{-\alpha} \right),
\end{equation*}   
i.e., in this case the convergence rate is independent of the dimension $s$. 

By proving a matching lower bound for the worst-case error of arbitrary cubature rules,
see Theorem~\ref{LowerBound}
and Corollary~\ref{Cor_LowerBound},
we show that QMC cubature rules based on $(t,m,s)$-nets exhibit in fact the optimal rate of convergence. It is well-known that in the worst-case setting of numerical integration within the class of all deterministic algorithms cubature rules are indeed optimal, see \cite{TWW88}; hence our result also implies that non-linear and adaptive algorithms cannot achieve a better convergence rate of the worst-case error  than QMC cubature rules based on $(t,m,s)$-nets.

 Let us mention that in the special case $p=\infty =q$ the upper bound \eqref{eq:upper_bound_Haar} was proved earlier by Entacher in \cite{Ent97}.
 Furthermore, in the Hilbert space setting $p=2=q$ Heinrich et al. proved  in \cite{HHY04} the weaker result that the average worst case error over all full randomized scramblings \cite{Owen95} 
 of a fixed $(t,m,s)$-net satisfies the upper bound \eqref{eq:upper_bound_Haar}.
They additionally showed for $p=2=q$ a matching lower bound for arbitrary cubature rules. Nevertheless, our proof approach for the general upper bound \eqref{eq:upper_bound_Haar} relies on an exactness result on finite dimensional subspaces of $\mathcal{H}_{\wav, \alpha, s, p,q}$, see Lemma~\ref{Exakt}, and differs significantly from the proof appoaches for the upper bounds in \cite{Ent97, HHY04}. 
 For more details, please see
 Remark~\ref{Rem:Ent_HHY}.
 
A quantitative comparison of the performance of $(t,m,s)$-nets and of Smolyak algorithms \cite{Smo63}, i.e., cubatures whose integration nodes form a sparse grid, on Haar wavelet spaces is provided in Remark~\ref{Rem:Smolyak}.

In Section~\ref{Sec:Fractional} we discuss the spaces of fractional smoothness 
$\mathcal{H}_{\alpha, s, p,q}$ introduced in \cite{D08}. We remark that for fixed parameters $\alpha$, $s$, and $p$ different values of $q$ lead always to the same vector spaces endowed with equivalent norms; that is why in the course of the paper we often surpress the parameter $q$. In Subsections~\ref{Subsec:Frac_Univ} and 
\ref{Subsec:Frac_Multi} we define the spaces and (re-)state their important properties.
The new results in Subsection~\ref{Subsec:Frac_Univ} are Lemma~\ref{Lemma0} (iii) and the closure of a gap in the proof of \cite[Lemma~1]{D08}, see Remark~\ref{Rem:Proof_Gap}. New in Subsection~\ref{Subsec:Frac_Multi} is Lemma~\ref{s-Derivative} (ii),
which relates the representations of 
$\mathcal{H}_{\alpha, s, p,q}$ with the help of integral kernels, (partial) weak derivatives in the sense of Riemann and Liouville, and the anchored decomposition 
of functions in $\mathcal{H}_{\alpha, s, p,q}$. 
The anchored decomposition is also known as cut-HDMR, where HDMR stands for ``high-dimensional model representation'', see \cite{RA99} or \cite{Gri06, KSWW10}. Another new result stated in Subsection~\ref{Subsec:Frac_Multi} is Theorem~\ref{tensor3}, 
which describes for $1<p<\infty$ the tensor product structure of $\mathcal{H}_{\alpha, s, p,q}$ with the help of the $p$-nuclear tensor norm.
In Subsection~\ref{Subsec:SKHI} we show that the fractional Koksma-Hlawka inequality derived in \cite{D08} is actually sharp, i.e., that
\begin{equation}\label{Int:sharp_KH_inequality}
e^{\wor}(Q_P, \mathcal{H}_{\alpha, s,p, q}) = D^*_{\alpha,s,p',q'}(P)
\end{equation}
holds,
see Theorem~\ref{KoksmaHlawka}. 
This results may be interpreted as follows: If one considers as quality criterion for sample points  the fractional discrepancy, then the spaces of fractional smoothness are indeed the appropriate function spaces to study QMC integration.

In Section~\ref{Sec:Besov} we turn to Besov and Bessel potential spaces. 
We study  in Subsection~\ref{Subsec:Besov_Univ} the univariate case $s=1$, define the Besov space $B^\alpha_{p, q}([0,1])$ and the Bessel potential spaces $H^\alpha_p([0,1])$ on $[0,1]$.
Furthermore, we prove two embedding results, see Theorem~\ref{Winfried} and~\ref{BesovHaar}, which can be combined with known results, cf. Lemma~\ref{einbettung2},   
to obtain for $\alpha\in (0, 1)$ and $\alpha^{-1} < p  < \infty$ the continuous embeddings
\begin{equation*}
\mathcal{H}_{\alpha, p} := \mathcal{H}_{\alpha, 1, p}  \hookrightarrow H^\alpha_{p}([0,1])
\hookrightarrow 
\begin{cases}
B^\alpha_{p, 2}([0,1]) \hookrightarrow \HH_{\wav,\alpha,1, p, 2}, & \text{if $p < 2$,}\\
 B^\alpha_{2, 2}([0,1]) \hookrightarrow \HH_{\wav,\alpha, 1, 2, 2}, & \text{if $p \ge 2$,}
\end{cases}
\end{equation*}
see Corollary~\ref{DirEmb}.

In Section~\ref{Subsec:Besov_Multi} 
we define Bessel potential spaces $S^\alpha_p H([0,1]^s)$ of $s$-variate functions, also called (fractional) Sobolev spaces of dominated mixed smoothness, and Besov spaces of mixed dominated smoothness $S^\alpha_{p,q} B([0,1]^s)$, and state known results about their tensor product structure,
see Proposition~\ref{tensor4}.
Since the proofs of the univariate embedding results  are already involved, we aim to achieve 
corresponding $s$-variate embedding results mainly with the help of the theory of tensor products of Banach spaces. 

We provide embedding results between the different types of spaces of multivariate functions
in Corollary~\ref{Winfried2} 
and~\ref{Cor:Neu_Equation_60+}, 
and Theorem~\ref{tensor18a},
which can be summarized as follows: For $\alpha\in (0, 1)$ and $\alpha^{-1} < p  < \infty$ we have the continuous embeddings
\begin{equation*}
\mathcal{H}_{\alpha, s, p} \hookrightarrow S^\alpha_{p}H([0,1]^s)
\hookrightarrow 
\begin{cases}
S^\alpha_{p, 2}B([0,1]^s) \hookrightarrow \HH_{\wav,\alpha, s, p, 2}, & \text{if $p < 2$,}\\
S^\alpha_{2, 2}(B[0,1]^s) \hookrightarrow \HH_{\wav,\alpha, s, 2, 2}, & \text{if $p \ge 2$,}
\end{cases}
\end{equation*}
cf. Corollary~\ref{DirEmb_Multi}.

In Section~\ref{Subsec:QMC_Frac_Disc} we consider again quasi-Monte Carlo integration: Our results for QMC integration on Haar wavelet spaces and the embedding results from Section~\ref{Subsec:Besov_Multi} immediately yield upper error bounds for QMC intgration on Besov and Sobolev spaces of dominated mixed smoothness and on spaces of fractional smoothness, respectively.
Let us state the achieved upper bounds in more detail:
Let $Q_P(f)$ again be a QMC cubature rule as in \eqref{eq:qmc_cub}, based on an arbitrary $(t,m,s)$-net $P$ in some base $b\ge 2$, and  let $N := |P| = b^m$.  Let $p,q\in [0,\infty]$ and 
$\alpha \in [1/p,1)$ if $q=1$ and $\alpha \in (1/p,1)$ if $q>1$. Then we obtain for the corresponding Besov spaces of dominated mixed smoothness that
\begin{equation*}
e^{\wor}(Q_P, S^{\alpha}_{p,q}B([0,1]^s)) = O \left( N^{-\alpha} \ln(N)^{\frac{s-1}{q'}} \right),
\end{equation*}   
where 
the implicit constant in the big-$O$-notation 
does not depend on the specific construction of the $(t,m,s)$-net $P$, see Corollary~\ref{Cor:Nets_Besov_Spaces}. In particular, we have for $q=1$ 
\begin{equation*}
e^{\wor}(Q_P, S^{\alpha}_{p,1}B([0,1]^s))  = O \left( N^{-\alpha} \right).
\end{equation*}   
As pointed out in Remark~\ref{Rem:Known_Int_Besov_Sobolev}, 
these upper bounds are best possible, showing that QMC cubature rules based on $(t,m,s)$-nets achieve also the optimal rate of convergence on the corresponding Besov spaces of dominated mixed smoothness, and the same holds true if we restrict ourselves to Besov spaces of periodic functions of dominating mixed smoothness.

Now to Sobolev spaces of dominated mixed smoothness and to spaces of fractional smoothness,
see Corollary~\ref{Cor:Pre_Hauptresultat}
and~\ref{Hauptresultat}:
For $1<p<\infty$ and $\alpha \in (\max\{1/2, 1/p\}, 1)$ and for $X= S^\alpha_p([0,1]^s)$ as well as for 
$X= \mathcal{H}_{\alpha, s,p}$ we obtain
\begin{equation*}
e^{\wor}(Q_P, X) 
= O \left( N^{-\alpha} \ln(N)^{\frac{s-1}{2}} \right);
\end{equation*}
recall that for all $q_1, q_2 \in [1,\infty]$ the spaces $\mathcal{H}_{\alpha, s,p, q_1}$ and $\mathcal{H}_{\alpha, s,p,q_2}$ are equal as vector spaces and are endowed with equivalent norms -- that is why we surpressed the $q$-index.
As pointed out in Remark~\ref{Rem:Known_Int_Besov_Sobolev}, 
this upper bound is best possible for Sobolev spaces of dominating mixed smoothness, establishing that QMC cubature rules based on $(t,m,s)$-nets achieve also the optimal rate of convergence on these spaces, and the same holds true if we restrict ourselves to Sobolev spaces of periodic functions of dominating mixed smoothness.

Furthermore, due to \eqref{Int:sharp_KH_inequality}, we get also the upper bound
$$
D^*_{\alpha,s,p',q'}(P) = O \left( N^{-\alpha} \ln(N)^{\frac{s-1}{2}} \right),
$$
$q\in [1,\infty]$ arbitrary, for the fractional discrepancy, see again Corollary~\ref{Hauptresultat}. 
Note that this bound holds for \emph{all} $(t,m,s)$-nets $P$; as explained in Remark~\ref{Rem:No_Ext_to_alpha_1}, this result does not hold true anymore in the case where $\alpha = 1$.

We also have results for Sobolev spaces of dominated mixed smoothness and spaces of fractional smoothness for a lower degree of smoothness $1/p < \alpha \le 1/2$, see Corollary~\ref{Cor:Pre_Hauptresultat}(ii) 
and~\ref{Hauptresultat}(ii) for details.

As already mentioned above, we provide with our final Section~\ref{Sec:Appendix}
an appendix that contains on the one hand the most important definitions and results on tensor products of Banach spaces and of linear operators, and on the other hand the proofs of the results on the tensor product structures of Haar wavelet spaces, Besov and Sobolev spaces of dominating mixed smoothness, and spaces of fractional smoothness that are stated in Section~\ref{HAAR-s}, 
\ref{Subsec:Frac_Multi}, 
and~\ref{Subsec:Besov_Multi}.


\subsection{Some Remarks on Notation}

Let us close the introduction with some remarks concerning notation:
Throughout the paper we assume that $1 \le p, q \le \infty$. 
Furthermore, we use the convention that $1/\infty = 0$ and, consequently, that $r^{1/\infty} = 1$ for all $r\in \RR$.

For $s\in \N$ we put $S:= \{1,\ldots, s\}$. For a set $A$ we denote by $|A|$ its cardinality if $A$ is finite and its Lebesgue measure if $A$ is an infinite measurable subset of $\RR^s$. For a function $f: \dom(f) \to \C$  the support of $f$ is defined by
\begin{equation*}
\supp(f) := \{ x\in \dom(f) \,|\, f(x) \neq 0\}
\end{equation*}
and its closure is denoted by $\overline{\supp}(f)$.  By $\langle \cdot, \cdot \rangle$ we denote the scalar product in $L_2([0,1]^s)$ and by $\|\cdot\|_{L_p}$ the norm in $L_p([0,1]^s)$. 


\section{Haar Wavelet Spaces}
\label{Sec:Haar}


\subsection{The Univariate Case}
\label{HAAR-1}

Let $\varphi: \RR\to \RR$ be the characteristic function of the half-open unit interval,
i.e., $\varphi := 1_{[0,1)}$. Let us fix an integer base $b\ge 2$ and put 
$$\Delta_{-1} := \{0\} \hspace{2ex}\text{and}\hspace{2ex}\Delta_j := \{0,1,\ldots, b^j-1\} \hspace{2ex}\text{for $j\in{\NN}_0$.} $$
We consider the dilated and translated functions
$$\varphi^j_k(x) := b^{j/2}\varphi(b^jx -k) \hspace{2ex}
\text{for all $j\in{\NN}_0$ and $k\in \Delta_j$,}$$
and define the spaces 
$$V^j := \spann\{ \varphi^j_k \,|\, k\in\Delta_j\}.$$
The elements of $V^j$ are all functions that are constant on 
each \emph{elementary interval} of length $b^{-j}$
\begin{equation*}
E^j_k := [k b^{-j}, (k+1)b^{-j}), \hspace{2ex} k\in\Delta_j.
\end{equation*}
The spaces $V^j$ have obviously dimension $b^j$ and satisfy
$$V^0 \subset V^1 \subset V^2 \subset \cdots \subset L_2([0,1]).$$
For $j\in \NN$ we define the space $W^j$ of dimension $(b-1)b^{j-1}$
as the orthogonal complement of $V^{j-1}$ in $V^{j}$, i.e.,
$$W^j := \{\psi\in V^j \,|\, \langle \psi, \eta\rangle = 0
\hspace{2ex}\text{for all $\eta\in V^{j-1}$} \}.$$
If we put $W^0 := V^0$, then it is well
known that $\oplus^\infty_{j=0} W^j = L_2([0,1])$; 
here and in the rest of the paper $\oplus$ denotes an orthogonal sum in a Hilbert space.

For $j \in \NN_0$ let $P_j: L_2([0,1]) \to V^j$ be the orthogonal
projection onto $V^j$. The orthogonal projection onto
$W^j$ is given by $P_j - P_{j-1}$ (where, in the case $j=0$, we use the
convention that $P_{-1}:=0$ is the orthogonal projection onto
$V^{-1}:=\{0\}$).

We can represent the orthogonal projections $P_j-P_{j-1}$ with the help
of wavelets. For this purpose we can use an orthonormal basis of the
space $W^1$ and translated and scaled versions of those basis functions
to get a system of orthonormal wavelets, as it was done, e.g., in
\cite{GLS07}. Here we follow the approach used in \cite{Owen97, HHY04}
and consider not an orthonormal wavelet system, but a tight
wavelet frame.
We start by defining the univariate basic wavelet functions
\begin{equation*}
\psi_i(x) := b^{1/2} 1_{[b^{-1}i, b^{-1}(i+1))}(x)
-b^{-1/2} 1_{[0,1)}(x),
\hspace{2ex}\text{$i\in\Delta_1$.}
\end{equation*}
Obviously it follows $\int \psi_i (x)dx =0$ for all $i\in\Delta_1$.
These functions span the space $W^1$, and for $j\in \N$ the $b^{j}$ functions
\begin{equation*}
\begin{split}
\psi^j_{i,k}(x) :=& b^{(j-1)/2}\psi_i(b^{j-1}x-k)\\
=& b^{j/2-1} \left( b 1_{[b^{-j}(bk+i), b^{-j}(bk+i+1))}(x)
- 1_{[b^{1-j}k, b^{1-j}(k+1))}(x) \right)
, \hspace{2ex}
\text{$i\in\Delta_1$, $k\in\Delta_{j-1}$,}
\end{split}
\end{equation*}
span the space $W^j$.
Let us furthermore put
\begin{equation*}
\psi^0_{0,0} := \varphi \hspace{2ex} \text{and}\hspace{2ex}
\nabla_j :=
\left\{\begin{array}{ll} \{0\}, & \hspace{3ex}\text{if $j=0$,}
\\ \{0,1,\ldots, b -1\}, & \hspace{3ex}\text{if $j \in \N$.}\end{array} \right.
\end{equation*}

Clearly, 
\begin{equation*}
\langle \psi^0_{0,0}, \psi^j_{i,k} \rangle = \int^1_0 \psi^j_{i,k}(x)\, dx = 0
\hspace{2ex}\text{for all $j\in \N$, $i\in \nabla_j$, $k\in \Delta_{j-1}$.}
\end{equation*}
Furthermore,
\begin{equation*}
\supp(\psi^j_{i,k}) = E^{j-1}_k \hspace{2ex}\text{with}\hspace{2ex} |E^{j-1}_k| = b^{1-j} 
\hspace{2ex}\text{for all $j\in \N$, $i\in \nabla_j$, $k\in \Delta_{j-1}$,}
\end{equation*}
and 
\begin{equation*}
\supp(\psi^0_{0,0}) = [0,1) =: E^{-1}_0  \hspace{2ex}\text{with}\hspace{2ex}  |E^{-1}_0| = 1 \le b^{1}.
\end{equation*}
All wavelets $\psi^j_{i,k}$ are constant on elementary intervals $E^j_{\kappa}$, $\kappa \in \Delta_{j}$, of lengths $b^{-j}$.

The wavelets defined above are not orthogonal and even not linearly independent, but it is easily verified that they satisfy the following
two identities
\begin{equation}
\label{HHY6}
\sum_{i\in\nabla_j} \psi^j_{i,k}(x) = 0
\hspace{2ex}
\text{for all $j \in\NN$, $k\in\Delta_{j-1}$, $x\in [0,1]$,}
\end{equation}
and
\begin{equation}
\label{HHY7}
\langle \psi^j_{i,k}, \psi^{j'}_{i',k'} \rangle
= \delta_{jj'}\delta_{kk'}(\delta_{ii'}-b^{-1})
\hspace{2ex}
\text{for all $j,j' \in\NN$, $i\in \nabla_{j}$, $i'\in \nabla_{j'}$, $k\in\Delta_{j-1}$, $k'\in \Delta_{j'-1}$.}
\end{equation}
These two identities imply  for each $f\in L_{1}([0,1])$
\begin{equation}\label{fourier_sum_zero}
\sum_{i\in\nabla_j} \langle f, \psi^j_{i,k} \rangle = 0
\hspace{2ex}
\text{for all $j \in\NN$, $k\in\Delta_{j-1}$},
\end{equation}
and, if additionally $f\in L_{2}([0,1])$, 
\begin{equation}
\label{express}
(P_j-P_{j-1})f = \sum_{i\in\nabla_j} \sum_{k\in\Delta_{j-1}}
\langle f, \psi^j_{i,k} \rangle \psi^j_{i,k}
\hspace{3ex}\text{for all $j\in\NN_0$}
\end{equation}
(cf. \cite[Lemma~3]{Owen97}, where a corresponding
multi-resolution on the whole real line $\RR$ is considered).

Let $p\in [1, \infty]$. Due to H\"older's inequality and the fact that the wavelets $\psi^j_{i,k}$ are $L_\infty$-functions, 
we can define $P_j-P_{j-1}$ via identity (\ref{express}) for all $f\in L_p([0,1])$ resulting in a continuous linear operator $P_j-P_{j-1}: L_p([0,1]) \to L_\infty([0,1])$. 

For $f\in L_1([0,1])$ we may formally define 
\begin{equation}\label{sf_1}
\mathcal{S}(f) := \sum^\infty_{j=0} \sum_{k\in \Delta_{j-1}}
\sum_{i\in \nabla_{j}} \langle f, \psi^{j}_{i,k} \rangle
\psi^{j}_{i,k}.
\end{equation}
We write 
$$
\mathcal{S}(f) \equiv f
$$ 
if $\mathcal{S}(f)$ converges pointwise and $f= \mathcal{S}(f)$ almost everywhere;
in this case we identify $f$ with $\mathcal{S}(f)$ and consequently point evaluation of $f$ in $x$ is well-defined for 
every $x\in [0,1]$. 

For a finite $M\subset \NN_0$ we additionally use the notation
\begin{equation}\label{partial_sum_S}
\mathcal{S}_M(f) := \sum_{j\in M} \sum_{k\in \Delta_{j-1}}
\sum_{i\in \nabla_{j}} \langle f, \psi^{j}_{i,k} \rangle
\psi^{j}_{i,k}.
\end{equation}
Due to \eqref{HHY6} and \eqref{HHY7},
\begin{equation}\label{eq:coefficients_finite_sum}
\langle \mathcal{S}_M(f), \psi^j_{i,k} \rangle =
\begin{cases}
\langle f, \psi^j_{i,k} \rangle, & \text{if $j\in M$,}\\
0, & \text{otherwise.}
\end{cases}
\end{equation}

\begin{definition}\label{Def:Wavelet_Space_Uni}
Let $p \in [1,\infty]$ and
$\alpha\in\RR$. For $q\in [1,\infty)$ we define the \emph{Haar wavelet space}
\begin{equation}\label{def:wavelet_spaces_1}
 \mathcal{H}_{\wav,\alpha,p,q} :=
\{ f\in L_1([0,1]) \,|\,  \|f\|_{\wav,\alpha,p,q} <\infty\},
\end{equation}
where 
\begin{equation*}
\|f\|^q_{\wav,\alpha,p,q} :=
\sum_{j=0}^\infty b^{q (\alpha-1/p + 1/2) j} \left(
\sum_{k\in\Delta_{j-1}} \sum_{i\in\nabla_j}
|\langle f, \psi^j_{i,k} \rangle|^p \right)^{q/p},
\end{equation*}
with the obvious modifications for $p =\infty$.
For $q= \infty$ we put 
\begin{equation}\label{def:wavelet_spaces_2}
 \mathcal{H}_{\wav,\alpha,p,\infty} :=
\{ f\in L_1([0,1]) \,|\,   \|f\|_{\wav,\alpha,p,\infty} <\infty \hspace{1ex}\text{and}\hspace{1ex} f \equiv \mathcal{S}(f) \},
\end{equation}
where 
\begin{equation*}
\|f\|_{\wav,\alpha,p,\infty} :=
\sup_{j\in \N_0} \left( b^{(\alpha-1/p + 1/2) j} \left(
\sum_{k\in\Delta_{j-1}} \sum_{i\in\nabla_j}
|\langle f, \psi^j_{i,k} \rangle|^p \right)^{1/p} \right) ,
\end{equation*}
with the obvious modifications for $p =\infty$.
\end{definition}

The following result ensures that $\| \cdot \|_{\wav,\alpha,s,p,q}$ is a norm on 
$ \HH_{\wav,\alpha,s,p,\infty}$ and not only a semi-norm: Let $f\in L_1([0,1])$. 
\begin{equation}\label{cond_ensures_norm}
\text{If}\hspace{1ex} \langle f, \psi^j_{i,k} \rangle = 0 \hspace{1ex} \text{for all $j\in \N_0$, $k\in\Delta_{j-1}$, $i\in\nabla_j$, then $f=0$ almost everywhere,}
\end{equation}
see, e.g., \cite[Theorem~2.1(i)]{Tri10} for the case $b=2$.


\begin{lemma}\label{Rem_meaningful_1}
Let $p,q \in [1,\infty]$, $\alpha \in \RR$,  and $f \in \mathcal{H}_{\wav, \alpha,p,q}$. 
\begin{itemize}
\item[(i)] Let $q <\infty$. Then the sum $\mathcal{S}(f)$ 
converges unconditionally in $\HH_{\wav,\alpha,p,q}$ to $f$.
\item[(ii)] Let $\alpha \ge 1/p$ if $q=1$ and $\alpha>1/p$ if $q\in (1,\infty]$. 
Then 
\begin{equation}\label{est_L2_L1_Hwav}
\|f\|_{L_2} \le C^{1/2}_{b,\alpha,p,q}  \|f\|_{L_1}^{1/2} \|f\|^{1/2}_{\wav,\alpha,p,q},
\end{equation}
where  $C_{b,\alpha,p,q} := b^{1/p'} \big( 1-b^{-q'(\alpha -1/p)} \big)^{-1/q'}$.
In particular, we have the continuous embedding 
$\HH_{\wav,\alpha,p,q} \hookrightarrow L_2([0,1])$.
Furthermore, for every $x\in [0,1]$ the sum
$\mathcal{S}(f)(x)$ converges absolutely, and $\mathcal{S}(f) = f$ almost everywhere.
Therefore we may identify $f$ with its wavelet expansion $\mathcal{S}(f)$.
With this identification, for every $x\in [0,1]$ the point evaluation functional $\delta_x: \HH_{\wav,\alpha,p,q} \to \RR$, $g\mapsto g(x)$ is continuous and its operator norm is at most $C_{b,\alpha,p,q}$. 
\end{itemize}
\end{lemma}

\begin{remark}\label{Rem:RKBS}
Let $p,q\in [1,\infty]$ and $\alpha \ge 1/p$ if $q=1$ and $\alpha>1/p$ if $q\in (1,\infty]$. 
Lemma \ref{Rem_meaningful_1} shows that the sets $\mathcal{H}_{\wav,\alpha,p,q}$ satisfy
\begin{equation}\label{alt_def_Haar_wavelet_space}
 \mathcal{H}_{\wav,\alpha,p,q} =
\{ f\in L_2([0,1]) \,|\,  \|f\|_{\wav,\alpha,p,q} <\infty  \hspace{1ex}\text{and}\hspace{1ex} f \equiv \mathcal{S}(f)\}.
\end{equation}
We decided to define the wavelet spaces via \eqref{def:wavelet_spaces_1} and \eqref{def:wavelet_spaces_2}, since definition \eqref{def:wavelet_spaces_2} is analogous to the one that was used in \cite{Ent97}. It is straightforward to show that the normed space $\mathcal{H}_{\wav,\alpha,p,q}$ is complete, i.e., a Banach space. 
\end{remark}

\begin{proof}[Proof of Lemma \ref{Rem_meaningful_1}]
(i) Let $M \subset \N_0$ be finite and $\mathcal{S}_M(f)$ be as in \eqref{partial_sum_S}.
Due to \eqref{eq:coefficients_finite_sum} we obtain for $p\in [1,\infty)$
\begin{equation*}
\| f -  \mathcal{S}_M(f)\|^q_{\wav, \alpha, p,q} = \sum_{j\in \NN\setminus M} b^{q(\alpha - 1/p +1/2)j} 
\left(  \sum_{k\in \Delta_{j-1}}
\sum_{i\in \nabla_{j}} | \langle f, \psi^{j}_{i,k} \rangle  |^p \right)^{q/p}
\end{equation*}
and, for $p=\infty$,
\begin{equation*}
\| f -  \mathcal{S}_M(f)\|^q_{\wav, \alpha, \infty, q} = \sum_{j\in \NN\setminus M} b^{q(\alpha +1/2)j} 
 \max_{k\in \Delta_{j-1}}
\max_{i\in \nabla_{j}} | \langle f, \psi^{j}_{i,k} \rangle  |^q.
\end{equation*}
Since the sum $\| f \|^q_{\wav, \alpha, p,q}$ converges absolutely (and thus, in particular, unconditionally), we see that $\mathcal{S}(f)$ converges unconditionally to $f$ in $\mathcal{H}_{\wav,\alpha,p,q}$.

(ii) Consider a fixed $j\in\NN_0$. For each $i\in \nabla_j$ the supports of all $\psi^j_{i,k}$, $k\in \Delta_{j-1}$, are disjoint. Hence we have for $p\in (1,\infty]$
\begin{equation*}
\begin{split}
\sum_{k\in \Delta_{j-1}} \sum_{i\in \nabla_{j}} | \langle f, \psi^{j}_{i,k} \rangle  |^{p'}
&\le \sum_{k\in \Delta_{j-1}} \sum_{i\in \nabla_{j}} b^{jp'/2}
\| f|_{\supp(\psi_{i,k}^j)} \|_{L_1}^{p'}\\
&\le b^{1 + jp'/2} \|f\|^{p'}_{L_1},
\end{split}
\end{equation*}
and  for $p=1$
\begin{equation*}
 \max_{k\in \Delta_{j-1}} \max_{i\in \nabla_{j}} | \langle f, \psi^{j}_{i,k} \rangle  | \le b^{j/2} \|f\|_{L_1}.
\end{equation*}
Furthermore, we get for any $x\in [0,1]$ in the case where $p\in (1,\infty]$
\begin{equation*}
\sum_{k\in \Delta_{j-1}} \sum_{i\in \nabla_{j}} | \psi^{j}_{i,k} (x) |^{p'} \le b^{1 + jp'/2}, 
\end{equation*}
and, in the case where $p=1$,
\begin{equation*}
 \max_{k\in \Delta_{j-1}} \max_{i\in \nabla_{j}} | \psi^{j}_{i,k}(x) | \le b^{j/2}.
\end{equation*}
Therefore, applying H\"older's inequality twice, we get for $p\in (1,\infty)$ and $q\in (1,\infty]$
\begin{equation*}
\begin{split}
&\sum^\infty_{j=0} \sum_{k\in \Delta_{j-1}} \sum_{i\in \nabla_{j}} | \langle f, \psi^{j}_{i,k} \rangle  |^{2}\\
\le &\sum^\infty_{j=0} \left( \sum_{k\in \Delta_{j-1}} \sum_{i\in \nabla_{j}} | \langle f, \psi^{j}_{i,k} \rangle  |^{p} \right)^{1/p}    \left( \sum_{k\in \Delta_{j-1}} \sum_{i\in \nabla_{j}} | \langle f, \psi^{j}_{i,k} \rangle  |^{p'} \right)^{1/p'} \\ 
\le &\|f\|_{\wav,\alpha,p,q} \left( \sum_{j=0}^\infty b^{-q'(\alpha -1/p)j} b^{q'/p'}  \|f\|_{L_1}^{q'} \right)^{1/q'}\\
\le &b^{1/p'} \left( 1- b^{-q'(\alpha -1/p)} \right)^{-1/q'} \|f\|_{\wav,\alpha,p,q} \|f\|_{L_1} <\infty,
\end{split}
\end{equation*}
since $\alpha > 1/p$. The same estimate holds for $p\in\{1,\infty\}$ and $q\in (1,\infty]$.
For $q=1$ we get 
\begin{equation*}
\begin{split}
\sum^\infty_{j=0} \sum_{k\in \Delta_{j-1}} \sum_{i\in \nabla_{j}} | \langle f, \psi^{j}_{i,k} \rangle  |^{2}
\le &\|f\|_{\wav,\alpha,p,1} \left( \sup_{j=0}^\infty b^{-(\alpha -1/p)j} b^{1/p'}  \|f\|_{L_1} \right)\\
\le &b^{1/p'} \|f\|_{\wav,\alpha,p,1} \|f\|_{L_1} <\infty,
\end{split}
\end{equation*}
since $\alpha \ge 1/p$. 

In particular, we obtain $f\in L_2([0,1])$ and thus
$$
\|f\|^2_{L_2} = \sum^\infty_{j=0} \sum_{k\in \Delta_{j-1}} \sum_{i\in \nabla_{j}} | \langle f, \psi^{j}_{i,k} \rangle  |^{2},
$$
which gives \eqref{est_L2_L1_Hwav}. Note that \eqref{est_L2_L1_Hwav} and 
$\|\cdot \|_{L_1} \le \| \cdot \|_{L_2}$ imply that $\mathcal{H}_{\wav,\alpha,p,q}$ is continuously embedded in $L_2([0,1])$, and consequently $\mathcal{S}(f)$ converges in $L_2([0,1])$ to $f$. 

Moreover, we obtain for any $x\in [0,1]$ by again employing H\"older's inequality
\begin{equation*}
\begin{split}
|\mathcal{S}(f)(x)| 
\le &  \sum^\infty_{j=0} \sum_{k\in \Delta_{j-1}} \sum_{i\in \nabla_{j}} | \langle f, \psi^{j}_{i,k} \rangle |
| \psi^{j}_{i,k} (x) |\\
\le &\sum^\infty_{j=0} \left( \sum_{k\in \Delta_{j-1}} \sum_{i\in \nabla_{j}} | \langle f, \psi^{j}_{i,k} \rangle  |^{p} \right)^{1/p}    \left( \sum_{k\in \Delta_{j-1}} \sum_{i\in \nabla_{j}} | \psi^{j}_{i,k} (x)|^{p'} \right)^{1/p'} \\ 
\le &b^{1/p'} \left( 1- b^{-q'(\alpha -1/p)} \right)^{-1/q'} \|f\|_{\wav,\alpha,p,q}  <\infty
\end{split}
\end{equation*}
with the obvious modifications in the case where $p\in \{1,\infty\}$.
This shows that $\mathcal{S}f(x)$ converges absolutely. Since $\mathcal{S}f = f$ in $L_2([0,1])$, we have $\mathcal{S}(f) = f$ almost everywhere. Moreover, by identifying $\mathcal{S}(f)$ and $f$, we see that the point evaluation functional $\delta_x$ is continuous with operator norm $\|\delta_x\| \le C_{b,\alpha,p,q}$. 
\end{proof}

\begin{remark}\label{Rem:Haar_Inclusion}
Let $\alpha \in\RR$ be fixed. For  $p_1, p_2, q_1,q_2\in [1,\infty]$ with $p_1\ge p_2$ and $q_1\le q_2$ we have 
\begin{equation}\label{comparison_Haar_norms}
\|f\|_{\HH_{\wav,\alpha, p_1,q_1}} \ge  \|f\|_{\HH_{\wav,\alpha,p_2 ,q_2}}.
\end{equation}
This leads for  $q_1,q_2\in [1,\infty)$ or $q_1=\infty =q_2$ to
\begin{equation*}
\HH_{\wav,\alpha,p_1,q_1} \subseteq \HH_{\wav,\alpha,p_2,q_2};
\end{equation*}
the same holds for $q_1\in [1,\infty)$, $q_2=\infty$, and $\alpha \ge 1/p_1$ if $q_1=1$ and  $\alpha >1/p_1$ if $q_1>1$, cf. Lemma~\ref{Rem_meaningful_1}.
\end{remark}

Note that for $p=2=q$ and $\alpha \ge 0$ the space $\mathcal{H}_{\wav,\alpha,2,2}$ is continuously embedded
in $L_2([0,1])$ and is itself a Hilbert space.

\begin{lemma}\label{lemmalp}
Let $1\le p,q \le \infty$ and $\alpha >0$.
Then $\mathcal{H}_{\wav,\alpha,p,q}$ is continuously embedded into $L_p ([0,1])$.
Moreover, for all $f\in \mathcal{H}_{\wav,\alpha,p,q}$ the sum $\mathcal{S}(f)$ converges in 
$L_p ([0,1])$ unconditionally to $f$.
\end{lemma}

\begin{proof}
\emph{Step~1}: We show that the sum $\mathcal{S}(f)$ converges unconditionally in $L_p([0,1])$ to some function $\tilde{f}$. 
Let $M\subset \N_0$ be finite.
Let $N\in \N_0$, and assume that $M\cap \{0,1,\ldots ,N-1\} = \emptyset$. We denote the maximum of $M$ by $N'$.
Let us first consider the case $1\le p <\infty$. Then
\beqq
\|\mathcal{S}_{M} (f)\|_{L_p} &\le & \sum^\infty_{j=N} \Big\| \sum_{k\in \Delta_{j-1}}
\sum_{i\in \nabla_{j}} \langle f, \psi^{j}_{i,k} \rangle
\psi^{j}_{i,k}\Big\|_{L_p}  \, .
\eeqq
We fix $j\ge 1$.
By employing the disjointness of the supports and \eqref{fourier_sum_zero} we find
\beqq
\int_0^1 \Big|\sum_{k\in \Delta_{j-1}} \sum_{i\in \nabla_{j}} \langle f, \psi^{j}_{i,k} \rangle
\psi^{j}_{i,k}(x)\Big|^p dx & = &
 \sum_{k=0}^{b^{j-1}-1} \int_{b^{1-j}k}^{b^{1-j}(k+1)} \Big|\sum_{i\in \nabla_{j}} 
 \langle f, \psi^{j}_{i,k} \rangle\psi^{j}_{i,k}(x)\Big|^p dx
\\
&= &
\sum_{k=0}^{b^{j-1}-1} \sum_{r=0}^{b-1}
\int_{b^{-j}(bk+r)}^{b^{-j}(bk+r+1)} \Big|\langle f, \psi^{j}_{r,k} \rangle\, b^{j/2}
\Big|^p dx
\\
&=& b^{j(\frac{p}{2} - 1)}\, \sum_{k=0}^{b^{j-1}-1} \sum_{r=0}^{b-1}
|\langle f, \psi^{j}_{r,k} \rangle|^p\,.
\eeqq
The term $j=0$ can be handled similarly. Hence we obtain
\[
 \|\mathcal{S}_{M}(f)\|_{L_p} \le   \sum^\infty_{j= N}\,  b^{j(\frac 12 -\frac 1p)} 
 \Big( 
\sum_{k\in \Delta_{j-1}} \sum_{i\in \nabla_j} | \langle f, \psi^{j}_{i,k} \rangle |^p\, 
\Big)^{1/p}\, .
\]
Furthermore, we observe that  $\alpha >0 $ implies
\begin{equation}\label{eq1}
\|\mathcal{S}_{M}(f)\|_{L_p} \le 
\Big(\sum_{j=N}^\infty b^{-\alpha j}\Big)   \|f\|_{\wav,\alpha,p,\infty}\le \frac{b^{-\alpha N}}{1 - b^{-\alpha}}\, \|f\|_{\wav,\alpha,p,q}  
\end{equation}
for any $1 \le  q \le \infty$.
Analogously, we derive bound \eqref{eq1} also for $p=\infty$.
This shows that $\mathcal{S}(f)$ converges for $1\le p \le \infty$ unconditionally in the Banach space $L_p ([0,1])$ to some function $\tilde{f}$.

\emph{Step~2}: We show that $\tilde{f}=f$
almost everywhere. For arbitrary $N\in \N_0$ we put $M(N) := \{0,1,\ldots, N\}$. 
Then we obtain 
for all $j\in \N_0$, $k\in\Delta_{j-1}$ and $i\in\nabla_j$, due to the fact that $\psi^{j}_{i,k} \in L_{p'}([0,1])$ and identity \eqref{eq:coefficients_finite_sum}, that 
\begin{equation*}
 \langle \tilde{f}, \psi^{j}_{i,k} \rangle 
 = \lim_{N\to \infty}  \langle \mathcal{S}_{M(N)} f, \psi^{j}_{i,k} \rangle
 =  \langle f, \psi^{j}_{i,k} \rangle.
\end{equation*}
Hence \eqref{cond_ensures_norm} ensures
$f = \tilde{f}$ almost everywhere, i.e., $\mathcal{S}(f)$ converges unconditionally in $L_p([0,1])$ to $f$.
\\
This convergence result and  \eqref{eq1} for $N=0$ yield
\begin{equation*}
\| f \|_{L_p} 
\le \frac{1}{1 - b^{-\alpha}} \, \|f\|_{\wav,\alpha,p,q}
\end{equation*}
for all $1 \le p \le \infty$,
establishing the desired embedding result.
\end{proof}


\subsection{The Multivariate Case}
\label{HAAR-s}


For $\bsj\in \NN_0^s$ we put $|\bsj| := j_1+j_2+\cdots +j_s$.
We define the approximation space of level $L$ as
\begin{equation}
\label{approx}
V^{s,L} := \sum_{|\bsj|=L} \bigotimes^s_{\ell =1} V^{j_\ell}
\hspace{2ex}\text{and put}\hspace{2ex}
W^{s,L} := \bigoplus_{|\bsj|=L}\bigotimes^s_{\ell=1} W^{j_\ell}.
\end{equation}
In particular, $W^{s,0} = V^{s,0}$.
It is easily verified that
\begin{equation*}
V^{s,L} = V^{s,L-1}\oplus W^{s,L}
\hspace{2ex}\text{for all $L\in\NN$.}
\end{equation*}
Denote by $P_{s,L}$ the orthogonal projection in $L_2([0,1]^s)$ onto
$V^{s,L}$. Define $P_{s,-1} =0$.
Then $P_{s,L}-P_{s,L-1}$ is the orthogonal projection
onto $W^{s,L}$. Due to (\ref{approx}) we have
\begin{equation}
\label{tensorproj}
P_{s,L}-P_{s,L-1} = \sum_{|\bsj|=L} \otimes^s_{\ell=1}
(P_{j_\ell}-P_{j_\ell-1}).
\end{equation}
For $\bsj\in\NN_0^s$ let us put
\begin{equation}\label{neu-ws-1}
\nabla_{\bsj} := \nabla_{j_1}\times\cdots\times \nabla_{j_s}
\hspace{2ex}\text{and}\hspace{2ex}
\Delta_{\bsj} := \Delta_{j_1}\times\cdots\times \Delta_{j_s},
\end{equation}
as well as
\begin{equation*}
J=J(\bsj) := | \{ \ell \in S \,|\, j_{\ell} > 0 \}|.
\end{equation*}
Furthermore, let
\begin{equation*}
\Psi^{\bsj}_{\bsi,\bsk} := \otimes^s_{\ell=1}
\psi^{j_\ell}_{i_\ell, k_\ell}
\end{equation*}
and ${\bf 1} := (1,1,\ldots,1) \in \ZZ^s$. 
Note that 
\begin{equation*}
\supp( \Psi^{\bsj}_{\bsi,\bsk} ) = \prod_{\ell = 1}^s E^{j_\ell -1}_{k_\ell} =: E^{\bsj- {\bf 1}}_{\bsk} 
\hspace{2ex}\text{with}\hspace{2ex} |E^{\bsj- {\bf 1}}_{\bsk}| 
= b^{J-|\bsj|} \le b^{s-|\bsj|}
\end{equation*}
for all $\bsj \in \N_0^s$, $\bsi \in \nabla_{\bsj}$, $\bsk \in \Delta_{\bsj - {\bf 1}}$, and the wavelets $\Psi^{\bsj}_{\bsi,\bsk}$
are constant on all \emph{elementary $s$-dimensional intervals} $E^{\bsj}_{{\bf \kappa}}$, ${\bf \kappa}\in \Delta_{\bsj}$,  of volume $b^{-|\bsj|}$.
From (\ref{HHY7}) follows 
\begin{equation}\label{eq:integral_wavelet}
\int_{[0,1)^s} \Psi^{{\bf 0}}_{{\bf 0}, {\bf 0}}(\bsx) \, {\rm d}\bsx = 1
\hspace{3ex}\text{and}\hspace{3ex}
\int_{[0,1)^s} \Psi^{\bsj}_{\bsi,\bsk}(\bsx) \, {\rm d}\bsx = 0
\hspace{3ex}\text{if}\hspace{3ex}
|\bsj|\ge 1.
\end{equation}
This trivially implies for all $f\in L_1([0,1)^s)$
\begin{equation}\label{eq:integral_f}
\int_{[0,1)^s} f(\bsx) \, {\rm d}\bsx = \langle f, \Psi^{{\bf 0}}_{{\bf 0}, {\bf 0}} \rangle = \sum^\infty_{|\bsj|=0} \sum_{\bsk\in \Delta_{{\bf j-1}}}
\sum_{\bsi\in \nabla_{\bsj}} \langle f, \Psi^{\bsj}_{\bsi,\bsk} \rangle
\int_{[0,1)^s} \Psi^{\bsj}_{\bsi,\bsk}(\bsx) \, {\rm d}\bsx.
\end{equation}

Due to (\ref{tensorproj}) and (\ref{express}) we get
\begin{equation}
\label{s-express}
(P_{s,L}-P_{s,L-1})f =
\sum_{|\bsj|=L} \sum_{\bsi\in\nabla_{\bsj}}
\sum_{\bsk\in\Delta_{\bsj-{\bf 1}}} \langle f, \Psi^{\bsj}_{\bsi,\bsk}
\rangle \Psi^{\bsj}_{\bsi,\bsk}
\end{equation}
for all $L\in\NN_0$.
Since the right hand side is, for all $f\in L_1([0,1]^s)$,
a well-defined function in $L_\infty([0,1]^s)$, we may in particular
define
$P_{s,L}-P_{s,L-1}$ in this way as a continuous linear operator
on $L_p([0,1]^s)$ for all $p\in [1,\infty]$.

Analogously to the univariate case, we define for $f\in L_1([0,1]^s)$ formally
\begin{equation}
\label{sf}
\mathcal{S}(f) := \sum^\infty_{|\bsj|=0} \sum_{\bsk\in \Delta_{{\bf j-1}}}
\sum_{\bsi\in \nabla_{\bsj}} \langle f, \Psi^{\bsj}_{\bsi,\bsk} \rangle
\Psi^{\bsj}_{\bsi,\bsk}
\end{equation} 
and use the notation $f\equiv \mathcal{S}(f)$ to indicate that $\mathcal{S}(f)$ converges pointwise, coincides almost everywhere with $f$ and that we identify $f$ with $\mathcal{S}(f)$.

\begin{definition}\label{Def:Wavelet_Space_Multi}
Let $p \in [1,\infty]$ and $\alpha\in\RR$.  For $q\in [1, \infty)$ we define the \emph{$s$-variate Haar wavelet space}
\begin{equation*}
\HH_{\wav,\alpha,s,p,q} := \{ f\in L_1([0,1]^s) \,|\,
\|f\|_{\wav, \alpha,s,p,q} < \infty \},
\end{equation*}
where
\begin{equation*}
\|f\|^q_{\wav,\alpha,s,p,q} :=
\sum^\infty_{L=0} b^{q(\alpha-1/p+1/2)L} \sum_{|\bsj|=L}
\left( \sum_{\bsk\in \Delta_{\bsj -{\bf 1}}} 
\sum_{\bsi \in \nabla_{\bsj}}
|\langle f, \Psi^{\bsj}_{\bsi,\bsk} \rangle |^p
\right)^{q/p},
\end{equation*}
with the obvious modifications for $p = \infty$.
For $q = \infty$ we set
\begin{equation*}
\HH_{\wav,\alpha,s,p,\infty} := \{ f\in L_1([0,1]^s) \,|\,
\|f\|_{\wav, \alpha,s,p,\infty} < \infty \hspace{1ex}\text{and}\hspace{1ex} f \equiv \mathcal{S}(f) \},
\end{equation*}
where
\begin{equation*}
\|f\|_{\wav,\alpha,s,p,\infty} :=
\sup_{L\in \N_0}  \max_{|\bsj|=L} \left( b^{(\alpha-1/p+1/2)L}
\left( \sum_{\bsk\in \Delta_{\bsj -{\bf 1}}} 
\sum_{\bsi \in \nabla_{\bsj}}
|\langle f, \Psi^{\bsj}_{\bsi,\bsk} \rangle |^p
\right)^{1/p} \right),
\end{equation*}
again with the obvious modifications for $p= \infty$.
\end{definition}

\begin{lemma}\label{Rem_meaningful}
Let $p,q \in [1,\infty]$ and $f \in \mathcal{H}_{\wav, \alpha,s, p,q}$. 
\begin{itemize}
\item[(i)] Let $q <\infty$. Then the sum $\mathcal{S}(f)$ 
converges unconditionally in $\HH_{\wav,\alpha,s, p,q}$ to $f$.
\item[(ii)] Let $\alpha \ge 1/p$ if $q=1$ and $\alpha>1/p$ if $q\in (1,\infty]$. 
Then 
\begin{equation*}
\|f\|_{L_2} \le C^{s/2}_{b,\alpha,p,q}  \|f\|_{L_1}^{1/2} \|f\|^{1/2}_{\wav, \alpha, s, p,q},
\end{equation*}
with  $C_{b,\alpha,p,q}$
as in Lemma~\ref{Rem_meaningful_1}. 
In particular, we have the 
continuous embedding 
$\HH_{\wav,\alpha,s, p,q} \hookrightarrow L_2([0,1]^s)$.
Furthermore, for every $x\in [0,1]^s$ the sum
$\mathcal{S}(f)(x)$ converges absolutely, and $\mathcal{S}(f) = f$ almost everywhere.
Therefore we may identify $f$ with its wavelet expansion $\mathcal{S}(f)$.
With this identification, for every $x\in [0,1]^s$ the  operator norm of the point evaluation functional $\delta_x: \HH_{\wav,\alpha,s, p,q} \to \RR$ is bounded by $C_{b,\alpha,p,q}^s$ and, consequently,
convergence in $\HH_{\wav,\alpha, s, p,q}$ implies pointwise
convergence. 
\end{itemize}
\end{lemma}

The proof of Lemma~\ref{Rem_meaningful} follows exactly the lines of the proof of Lemma~\ref{Rem_meaningful_1}, which is why we omit it here. 

Furthermore, the statement of Remark~\ref{Rem:RKBS}, adapted in the obvious way, holds also true for $s\ge 2$.

Let $\alpha \in\RR$ be fixed. For  $p_1, p_2, q_1,q_2\in [1,\infty]$ with $p_1\ge p_2$ and $q_1\le q_2$ we have 
\begin{equation*}
\|f\|_{\HH_{\wav,\alpha, s, p_1,q_1}} \ge  \|f\|_{\HH_{\wav,\alpha,s, p_2 ,q_2}}.
\end{equation*}
This leads for  $q_1,q_2\in [1,\infty)$ or $q_1=\infty =q_2$ to
\begin{equation*}
\HH_{\wav,\alpha, s, p_1,q_1} \subseteq \HH_{\wav,\alpha, s, p_2 ,q_2};
\end{equation*}
the same holds for $q_1\in [1,\infty)$, $q_2=\infty$, and $\alpha \ge 1/ p_1$ if $q_1=1$ and  $\alpha >1/ p_1$ if $q_1>1$, cf. Lemma \ref{Rem_meaningful}.


\begin{lemma}\label{lemmalp_multivariate}
Let $s\in \N$, $1\le p,q \le \infty$, and $\alpha >0$.
Then $\mathcal{H}_{\wav,\alpha,s,p,q}$ is continuously embedded into $L_p ([0,1]^s)$.
Moreover, for all $f\in \mathcal{H}_{\wav,\alpha,s,p,q}$ the sum $\mathcal{S}(f)$ converges in 
$L_p ([0,1]^s)$ unconditionally to $f$.
\end{lemma}

\begin{proof}
The proof follows exactly the lines of the proof of Lemma~\ref{lemmalp}.
\end{proof}

For $p=2=q$ and $\alpha \ge 0$ the space $\HH_{\wav,\alpha,s,2,2}$
endowed with the scalar product
\begin{equation*}
\langle f,g \rangle_{\wav,\alpha, s,2,2}
= \sum^\infty_{L=0} b^{2\alpha L} \sum_{|\bsj|=L}
\sum_{\bsk\in \Delta_{\bsj -{\bf 1}}} \sum_{\bsi \in \nabla_{\bsj}}
\langle f, \Psi^{\bsj}_{\bsi,\bsk} \rangle
\langle \Psi^{\bsj}_{\bsi,\bsk}, g \rangle
\end{equation*}
is the complete $s$-fold tensor product Hilbert space of
$\HH_{\wav,\alpha,2,2}$ and is continuously embedded in $L^2([0,1]^s)$.
(For the definition of and elementary results about complete tensor
products of Hilbert spaces we refer the reader to \cite{Wei}.)
It is not hard to verify that
\begin{equation*}
\langle f,g \rangle_{\wav,\alpha,s,2,2} =
\sum^\infty_{L=0} b^{2\alpha L}
\langle f, (P_{s,L} - P_{s,L-1})g \rangle.\\
\end{equation*}
We now extend this tensor product result for the case $p=2=q$ to the case $p=q$ and arbitrary $p \in [1,\infty)$. 
To this purpose, the Hilbert space tensor norm has to be replaced in the case $p=1$ by  the
\emph{projective tensor norm} $\delta_1$, and in the case $1<p<\infty$ by the
\emph{$p$-nuclear tensor norm} $\delta_p$, which are both defined in Section~\ref{tpbs}
of the appendix on tensor products of Banach spaces.

\begin{theorem}\label{tensor17}
 Let $s\ge 2$, $\alpha >0$ and $1\le p < \infty$. Then
 \[
  \HH_{\wav,\alpha,s-1,p,p} \otimes_{\delta_p}  \HH_{\wav,\alpha,p,p} =
   \HH_{\wav,\alpha,s,p,p} 
 \]
 holds with equality of norms. 
 \end{theorem}

The proof of Theorem~\ref{tensor17} can be found in Section~\ref{APP:TPSS} of the appendix.


\subsection{$(t,m,s)$-Nets and Quasi-Monte Carlo Integration}
\label{Subsec:QMC}


Let $s\in \N$. We now investigate the approximation of integrals 
\begin{equation*}
I_s(f) := \int_{[0,1)^s} f(\bsx) \, {\rm d}\bsx,
\hspace{2ex}f\in \HH_{\wav, \alpha, s, p,q},
\end{equation*}
with the help of quasi-Monte Carlo (QMC) cubature rules.

Let $t,m \in \N_0$. A finite set $P\subset [0,1)^s$ is called a \emph{$(t,m,s)$-net in base $b$} if each elementary $s$-dimensional interval $E^{\bsj}_{\bsk}$ of volume $b^{t-m}$ contains exactly $b^t$
points. In particular, $|P| = b^m$. 

In the following 
we always consider Haar wavelet spaces $\HH_{\wav, \alpha, s, p,q}$ and $(t,m,s)$-nets with respect to the \emph{same} (fixed) integer base $b \ge 2$.

\begin{remark}
The abstract definition of $(t,m,s)$-nets 
is due to Niederreiter \cite{Nie87}. Previously, Sobol' studied $(t,m,s)$-nets in base $2$ \cite{Sob67}, and Faure provided constructions for arbitrary prime bases $b$ \cite{Fau82}. For further constructions and for background information about nets we refer the reader to \cite{Nie92, DP10}. 
An easy observation is that the $t$-parameter can essentially be seen as a quality parameter of a $(t,m,s)$-net: the smaller $t$, the more uniformly distributed the net. On the one-hand it is well-known that a necessary condition for the existence of a $(t,m,s)$-net in base $b$ is given by
\[
t \ge \frac{s-1}{b} - \log_b \frac{(b-1)(s-1) + b+ 1}{2},
\]
i.e., for fixed base $b$ the parameter $t$ has to increase with the dimension $s$.
On the other hand one can, e.g., obtain  for any prime number $b\ge s-1$ and any $m\in \N$ a $(0,m,s)$-net in base $b$ with the help of Faure's constructions from \cite{Fau82}. 
This shows that it is important to study not only the case where $b=2$, but to investigate arbitrary large prime bases $b$.
\end{remark}

We define the \emph{QMC-cubature } $Q_{P}$ corresponding to a $(t,m,s)$-net $P$
by
\begin{equation}
\label{qmc-alg}
Q_{P}(f) = b^{-m} \sum_{\bsp\in P} f(\bsp).
\end{equation}

We define the integration error of $Q_P$ applied to $f$ by
\begin{equation}\label{ie}
e(Q_P,f) := | I_s(f) - Q_P(f)|. 
\end{equation}
The \emph{worst-case error of $Q_P$ on $\HH_{\wav,\alpha,s,p,q}$} is defined by
\begin{equation}\label{wce}
e^{\wor}(Q_P, \HH_{\wav,\alpha,s,p,q}) :=
\sup_{f \in \HH_{\wav,\alpha,s,p,q}\,;\,\|f\|_{\wav,\alpha, s,p,q}\le 1} e(Q_P,f).
\end{equation}
Due to the definition of $\HH_{\wav,\alpha,s,p,q}$ and to Lemma~\ref{Rem_meaningful} these 
definitions are meaningful as long as 
$\alpha \ge 1/p$ if $q=1$ and  $\alpha > 1/p$ if $q>1$; if this condition is satisfied we can perform point evaluations by identifying each $f\in \HH_{\wav,\alpha,s,p,q}$ with the pointwise convergent series $\mathcal{S}(f)$.

\begin{lemma}\label{Exakt}
Let $P$ be a $(t,m,s)$-net in base $b$, and put $L:= m-t$.
Then the QMC-cubature  $Q_{P}$ defined in (\ref{qmc-alg})
is \emph{exact} on the aproximation space $V^{s,L}$, i.e.,
\begin{equation*}
Q_{P}(f) = I_s(f) 
\hspace{2ex}\text{for all $f\in V^{s,L}$.}
\end{equation*}
\end{lemma}

\begin{proof}
Since the integral $I_s$ and the QMC-algorithm $Q_{P}$ are linear,
we only need to show for arbitrary multi-indices $\bsj$,
$|\bsj| \le L$, $\bsi\in \nabla_{\bsj}$, and
$\bsk \in \Delta_{\bsj-{\bf 1}}$, that 
$I_s(\Psi^{\bsj}_{\bsi,\bsk}) = Q_{s,L}(\Psi^{\bsj}_{\bsi,\bsk})$.
Recall $\Psi^{\bsj}_{\bsi, \bsk}$ is constant on all
elementary $s$-dimensional intervals
$E^{\bsj}_{\bf a} =
\prod^s_{\ell=1} E^{j_\ell}_{a_\ell}$, which have volume
$b^{-|\bsj|} \ge b^{-L} = b^{t-m}$
and therefore contain $b^{m-|\bsj|}$ integration points of
$P$. Since $Q_{P}$ is properly normalized, we thus have
$I_s(\Psi^{\bsj}_{\bsi, \bsk}) = Q_{P}(\Psi^{\bsj}_{\bsi, \bsk})$.
\end{proof}

\begin{theorem}\label{Nets}
Let $p,q\in [1,\infty]$. Let $\alpha \ge 1/p$ if $q=1$ and $\alpha > 1/p$ if $q>1$.
Let $m\ge 1$, and let $P$ be a $(t,m,s)$-net in base $b$. 
Furthermore, put $N:=|P|=b^{m}$.
Then there exists a constant $C$, independent of $N$, such that
\begin{equation}
 \label{nets}
e^{\wor}(Q_P, \HH_{\wav,\alpha,s,p,q})
\le C N^{-\alpha}
\ln(N)^{\frac{s-1}{q'}}.
\end{equation}
In particular, we obtain for $q=1$ that the convergence rate of the worst-case error does not depend on the dimension $s$:
\begin{equation}
 \label{nets_q=1}
e^{\wor}(Q_P, \HH_{\wav,\alpha,s,p,1})
\le C N^{-\alpha}.
\end{equation}
\end{theorem}

\begin{proof}
Let $f\in \HH_{\wav,\alpha,s,p,q}$. Due to identity \eqref{eq:integral_f}, Lemma~\ref{Rem_meaningful}, and H\"older's inequality for sums we have
\begin{equation*}
\begin{split}
&|I_s(f) - Q_{P}(f)| =
\left| \sum^\infty_{|\bsj|=0} \sum_{\bsk\in\Delta_{\bsj-{\bf 1}}}
\sum_{\bsi \in \nabla_{\bsj}}\langle f, \Psi^{\bsj}_{\bsi,\bsk}
\rangle
\left( I_s(\Psi^{\bsj}_{\bsi,\bsk}) - Q_{P}(\Psi^{\bsj}_{\bsi,\bsk})
\right) \right|\\
\le &\|f\|_{\wav,\alpha,s,p,q} \left( \sum^\infty_{|\bsj|=0}
b^{-q'(\alpha-1/p+1/2) |\bsj|}
\left( \sum_{\bsk\in\Delta_{\bsj-{\bf 1}}} \sum_{\bsi\in\nabla_{\bsj}}
\left| I_s(\Psi^{\bsj}_{\bsi,\bsk}) - Q_{P}(\Psi^{\bsj}_{\bsi,\bsk})
\right|^{p'} \right)^{q'/p'} \right)^{1/q'},
\end{split}
\end{equation*}
with the obvious modifications in the cases where $p=1$ or $q=1$.
Put $L:=m-t$.
The identities \eqref{eq:integral_wavelet} and the exactness result Lemma~\ref{Exakt} give us
\begin{equation*}
\begin{split}
e^{\wor}(Q_P, \HH_{\wav,\alpha,s,p,q})^{q'}
\le& \sum^\infty_{|\bsj|=L+1} b^{-q'(\alpha-1/p+1/2) |\bsj|}
\left( \sum_{\bsk\in\Delta_{\bsj-{\bf 1}}} \sum_{\bsi\in\nabla_{\bsj}}
 \left|Q_{P} \left(\Psi^{\bsj}_{\bsi,\bsk}
\right) \right|^{p'} \right)^{q'/p'}\\
=&\, b^{-q'(L+t)} \sum^\infty_{|\bsj|=L+1} b^{-q'(\alpha-1/p+1/2)|\bsj|}
\left( \sum_{\bsi\in\nabla_{\bsj}}\sum_{\bsk\in\Delta_{\bsj-{\bf 1}}}
\left| \sum_{\bsp\in P}
\Psi^{\bsj}_{\bsi,\bsk}(\bsp) \right|^{p'} \right)^{q'/p'},
\end{split}
\end{equation*}
again with the obvious modifications in the cases where $p=1$ or $q=1$.

For $|\bsj| \ge L+1$ and $\bsi \in \nabla_{\bsj}$ the supports of the functions
$\Psi^{\bsj}_{\bsi, \bsk}$, $\bsk\in\Delta_{\bsj-{\bf 1}}$, are pairwise
disjoint, have volume at most $b^{J-L-1}$ and therefore contain at most $b^{t+J-1}$ points of $P$.
Recall that  $|P|=b^{L+t}$ and notice that $\|\Psi^{\bsj}_{\bsi,\bsk} \|_{\infty} = b^{|\bsj|/2}(1-b^{-1})^J$. 
Hence we get for $p\in (1, \infty]$
\begin{equation}\label{p>1}
\begin{split}
\left( \sum_{\bsi\in\nabla_{\bsj}}\sum_{\bsk\in\Delta_{\bsj-{\bf 1}}}
\left| \sum_{\bsp\in P} \Psi^{\bsj}_{\bsi,\bsk}(\bsp) \right|^{p'} \right)^{1/p'}
&\le \left( \sum_{\bsi \in \nabla_{\bsj}} b^{L-J+1} \left( b^{t+J-1}b^{|\bsj|/2}(1-b^{-1})^J \right)^{p'} \right)^{1/p'}\\
&\le b^{\frac{L+1}{p'} +t -1} (b-1)^{J} b^{|\bsj|/2},
\end{split}
\end{equation}
and for $p = 1$
\begin{equation}\label{p=1}
\max_{\bsi\in\nabla_{\bsj}}\max_{\bsk\in\Delta_{\bsj-{\bf 1}}}
\left| \sum_{\bsp\in P} \Psi^{\bsj}_{\bsi,\bsk}(\bsp) \right|
\le b^{t -1} (b-1)^{J} b^{|\bsj|/2}.
\end{equation}
These estimates yield for $q\in (1,\infty]$
\begin{equation*}
\begin{split}
e^{\wor}(Q_{P}, \HH_{\wav,\alpha,s,p,q})^{q'} 
\le& \,b^{-q'(L+1)/p} (b-1)^{sq'}
\sum^\infty_{|\bsj|=L+1} b^{-q'(\alpha-1/p)|\bsj|}\\
=& \,b^{-q'(L+1)/p} (b-1)^{sq'}
 \sum^\infty_{\nu = L+1} b^{ -q'\left(\alpha-1/p \right)\nu}
{\nu+s-1 \choose s-1}\\
=&
\left[ (b-1)^{sq'}\sum^\infty_{\nu=0} b^{-q'(\alpha-1/p)\nu}
\left( {\nu + L+s \choose s-1} m^{1-s} \right) \right]
 b^{-q'\alpha(L+1)} m^{s-1}.
 \end{split}
\end{equation*}
Hence
\begin{equation*}
e^{\wor}(Q_{P}, \HH_{\wav,\alpha,s,p,q}) \le  C(b,t,s,\alpha,p,q) N^{-\alpha} \log_b(N)^{\frac{s-1}{q'}},
\end{equation*}
where
\begin{equation*}
C(b,t,s,\alpha,p,q) :=  (b-1)^{s}b^{\alpha(t-1)} \left( \sum^\infty_{\nu=0} b^{-q'(\alpha-1/p)\nu}
\frac{1}{(s-1)!} \left( 1 + s+\nu \right)^{s-1} \right)^{1/q'},
\end{equation*}
and the sum in the last expression converges since $\alpha > 1/p$. 

For $q=1$ we get from (\ref{p>1}), (\ref{p=1}), and $\alpha \ge 1/p$ 
\begin{equation*}
\begin{split}
e^{\wor}(Q_{P}, \HH_{\wav,\alpha,s,p,1}) \le&
b^{-(L+t)} \sup_{|\bsj| \ge L+1} b^{-(\alpha -1/p +1/2)|\bsj|} \left( b^{\frac{L+1}{p'} +t -1} (b-1)^{J} b^{|\bsj|/2} \right)\\
\le& (b-1)^{s} b^{-\alpha(L+1)} = \left[ (b-1)^{s} b^{\alpha(t-1)} \right] N^{-\alpha}.
\end{split}
\end{equation*}
This concludes the proof of the theorem.
\end{proof}

\begin{remark}\label{Rem:Ent_HHY}
In \cite{Ent97} Entacher proves a result for the convergence of 
cubature rules based on nets on classes $\widetilde{E}^\alpha_s(C)$ of generalized Haar 
functions. Although his analysis
relies on  a different Haar wavelet expansion, it is not hard to see that his result is
equivalent to the special case of Theorem~\ref{Nets} where $p=\infty =q$. Notice that compared to \cite{Ent97}, our approach is conceptionally simpler and leads to a shorter proof.
Entacher shows that for $p=\infty= q$ the convergence rate in (\ref{nets}) is a lower bound in dimension $s=2$ for at least one cubature  rule based on a $(0,m,s)$-net in base $2$. 
But he does
not provide a lower bound for all cubature rules based on nets or, even more general, for arbitrary cubature rules.

In \cite{HHY04} Heinrich, Hickernell, and Yue consider also Haar wavelet
spaces, but confine themselves to the Hilbert space situation, i.e., to the spaces $\HH_{\wav,\alpha, s, p,q}$ for $p=2=q$.
They consider fully scrambled $(t,m,s)$-nets as introduced in \cite{Owen95} and prove that for a
$(t,m,s)$-net $P$ in base $b$ with $N=b^m$ points the average of the worst case error taken over the class $\mathcal{Q}_{\scr}(P)$ of all scrambled quasi-Monte Carlo cubature rules based on $P$, satisfies
\begin{equation}
\label{rms}
\rms_{Q\in\mathcal{Q}_{\scr}(P)} e^{\wor}(Q,\HH_{\wav,\alpha,s,2,2} )
= O \left( N^{-\alpha}\ln(N)^{\frac{s-1}{2}} \right)
\hspace{3ex}\text{for $\alpha > 1/2$.}
\end{equation}
(Here $\rms$ stands for ``root mean square''.)
Since Owen's scrambling maps $(t,m,s)$-nets again to $(t,m,s)$-nets,
our result Theorem~\ref{Nets} provides in the special case $p=2=q$ a stronger 
statement than 
(\ref{rms}), namely that not only the average
over $\mathcal{Q}_{\scr}(P)$, but already \emph{every single} 
$Q\in\mathcal{Q}_{\scr}(P)$ \emph{itself} satisfies the bound (\ref{rms}).

In \cite{HHY04} also upper bounds for the random-case and the average-case error
of cubature rules based on scrambled net are provided. The proofs of all upper bounds in \cite{HHY04} rely crucially on the Hilbert space structure of  $\HH_{\wav,\alpha, s, 2,2}$.
By establishing lower bounds for all three
error criteria for arbitrary cubature  rules, Heinrich et al. show that their upper
bounds exhibit optimal convergence rates. In particular, their lower bound for the worst-case error shows that our result Theorem \ref{Nets} is optimal in the case
where $p=2=q$. To demonstrate that the upper bound in Theorem \ref{Nets} is best possible in \emph{all cases}, we extend their lower bound in Corollary~\ref{Cor_LowerBound} to the case of arbitrary
$1\le p,q \le \infty$.
\end{remark}

\begin{theorem}
 \label{LowerBound}
Let $p,q \in [1, \infty]$, and let  $\alpha >0$.
Then there exists a constant $c$,
depending only on $\alpha, s,b$, and $q$, such that every cubature  rule
 $Q$ that uses $N\ge 3$ sample points satisfies
\begin{equation*}
 e^{\wor}(Q, \widetilde{\HH}_{\wav, \alpha, s,p,q}) \ge c N^{-\alpha}
\ln(N)^{\frac{s-1}{q'}},
\end{equation*}
where $\widetilde{\HH}_{\wav, \alpha, s,p,q}$ denotes the subspace $\cup_{L\in \N_0} V^{s,L}$ of $\HH_{\wav, \alpha, s,p,q}$, again endowed with the norm ${\|\cdot\|_{\wav, \alpha, s,p,q}}$.
\end{theorem}

Note that function evaluation is well-defined on $\widetilde{\HH}_{\wav, \alpha, s,p,q}$ for all $\alpha \in \RR$. Assume now that $\alpha \ge 1/p$ if $q=1$ and $\alpha > 1/p$ if $q>1$, so that function evaluation is also well-defined on $\HH_{\wav, \alpha, s,p,q}$.
Since $\widetilde{\HH}_{\wav, \alpha, s,p,q} \subseteq \HH_{\wav, \alpha, s,p,q}$, Theorem \ref{LowerBound} immediately implies the following corollary.

\begin{corollary}
 \label{Cor_LowerBound}
Let $p,q \in [1, \infty]$. Let  $\alpha \ge 1/p$ if $q=1$ and $\alpha > 1/p$ if $q>1$.
Then there exists a constant $c$,
depending only on $\alpha, s,b$, and $q$, such that every cubature  rule
 $Q$ that uses $N\ge 3$ sample points satisfies
\begin{equation*}
 e^{\wor}(Q, \HH_{\wav, \alpha, s,p,q}) \ge c N^{-\alpha}
\ln(N)^{\frac{s-1}{q'}}.
\end{equation*}
\end{corollary}

\begin{proof}[Proof of Theorem \ref{LowerBound}]
 We adapt the proof approach from \cite[Thm.~9]{HHY04}. Let $N$ be given, and let $m\in\NN$
be the uniquely determined number satisfying $b^{m-1} < 2N \le b^m$.
Let $P$ be an $N$-point set in $[0,1)^s$. For each multi-index $\bsj\in \N^s_0$ with $|\bsj|=m$ we
define functions $f_{\bsj}$ by
\begin{equation*}
 f_{\bsj}(x) = \left\{\begin{array}{ll} 1, & \hspace{3ex}\text{if $x\in E^{\bsj}_{\bsa}$
 with $P\cap E^{\bsj}_{\bsa}=\emptyset$, $\bsa \in \Delta_{\bsj}$,} \\ 0, & \hspace{3ex}\text{otherwise.}\end{array} \right.
\end{equation*}
Due to $b^m \ge 2N$, we obtain
\begin{equation*}
 \int_{[0,1]^s} f_{\bsj}(x)\,{\rm d}x \ge b^{-m}(b^m-N) \ge 1/2.
\end{equation*}
If we put
\begin{equation*}
 f_0 := \sum_{|\bsj|=m} f_{\bsj},
\end{equation*}
then $f_0\in V^{s,m}$ and $f_0(\bsp) = 0$ for all $\bsp\in P$. If $Q$ is any cubature  rule whose sample points are in $P$, then
\begin{equation*}
e(Q, f_0) = \int_{[0,1]^s} f_{0}(x)\,{\rm d}x = \sum_{|\bsj|=m} \int_{[0,1]^s} f_{\bsj}(x)\,{\rm d}x
\ge \frac{1}{2} { m+s-1 \choose s-1} \ge \frac{m^{s-1}}{2(s-1)!}.
\end{equation*}
Now we estimate $\|f_0\|_{\wav,\alpha,s,p,q}$. To this end we use the
estimate $|\supp(\psi^j_{i,k})| \|\psi^j_{i,k}\|_\infty  \le b^{-j/2}(b-1)$, which holds for all $j,k,i$.
H\"older's inequality implies
\begin{equation*}
 |\langle f_{\bsj'}, \Psi^{\bsj}_{\bsi, \bsk} \rangle |
\le |\supp(\Psi^{\bsj}_{\bsi,\bsk})| \|\Psi^{\bsj}_{\bsi,\bsk}\|_\infty 
\le b^{-\frac{|\bsj|}{2}}(b-1)^s,
\hspace{2ex}\text{$\bsj, \bsj'\in \N^s_0$, $\bsi \in \nabla_{\bsj}$, $\bsk\in \Delta_{\bsj -{\bf 1}}$.}
\end{equation*}
Furthermore, $j_\nu' < j_\nu$ for some $\nu\in S$ implies
$\langle f_{\bsj'}, \Psi^{\bsj}_{\bsi, \bsk} \rangle = 0$.
Put 
\begin{equation*}
A(\bsj,m) := \{\bsj'\in \NN^s_0 \,|\, |\bsj'|=m \, \wedge \, \forall \nu\in S: 
j_\nu' \ge j_\nu \}.
\end{equation*}
We get
\begin{equation*}
\|f_0\|^q_{\wav,\alpha,s,p,q}
= \sum_{L=0}^m b^{q(\alpha - \frac{1}{p} + \frac{1}{2})L} 
\sum_{|\bsj| =L}
\left( \sum_{k\in\Delta_{\bsj-{\bf 1}}} \sum_{\bsi\in\nabla_{\bsj}}
\left| \sum_{\bsj'\in A(\bsj,m)} 
\langle f_{\bsj'}, \Psi^{\bsj}_{\bsi, \bsk} \rangle
\right|^p \right)^{q/p},
\end{equation*}
with the obvious modifications if $p=\infty$ or $q=\infty$.
Let $|\bsj| = L \le m$.
Then
\begin{equation*}
|A(\bsj,m)| = { m-L +s-1 \choose s-1}.
\end{equation*} 
Let $p\in [1,\infty)$. Since $|\Delta_{\bsj-{\bf 1}}| |\nabla_{\bsj}|= b^{|\bsj|} = b^L$, we obtain
\begin{equation*}
\left( \sum_{k\in\Delta_{\bsj-{\bf 1}}} \sum_{\bsi\in\nabla_{\bsj}}
\left| \sum_{\bsj'\in A(\bsj,m)} 
\langle f_{\bsj'}, \Psi^{\bsj}_{\bsi, \bsk} \rangle
\right|^p \right)^{1/p} \le (b-1)^s b^{L(1/p -1/2)} { m-L +s-1 \choose s-1}.
\end{equation*}
Let $p = \infty$. Then we get the corresponding estimate
\begin{equation*}
\max_{k\in\Delta_{\bsj-{\bf 1}}} \max_{\bsi\in\nabla_{\bsj}}
\left| \sum_{\bsj'\in A(\bsj,m)} 
\langle f_{\bsj'}, \Psi^{\bsj}_{\bsi, \bsk} \rangle
\right|  \le (b-1)^s b^{-L/2} { m-L +s-1 \choose s-1}.
\end{equation*}
Hence we obtain for $p\in [1,\infty]$ and $q\in [1,\infty)$
\begin{equation*}
\begin{split}
\|f_0\|^q_{\wav,\alpha,s,p,q} 
\le& (b-1)^{qs} \sum^m_{L=0} b^{qL\alpha} \sum_{|\bsj|=L} 
{m-L+s-1 \choose s-1}^q \\
\le& (b-1)^{qs} \left[\sum_{L=0}^m b^{q\alpha(L-m)} 
{m-L+s-1 \choose s-1}^q \right] b^{q\alpha m} 
{m+s-1 \choose s-1}\\
\le& (b-1)^{qs} \left[\sum_{\nu=0}^\infty b^{-q\alpha\nu} 
{\nu+s-1 \choose s-1}^q \right] b^{q\alpha m} 
{m+s-1 \choose s-1}\\
\le& C^q  b^{\alpha qm}m^{s-1},
\end{split}
\end{equation*}
where $C$
is a constant only depending on $s,b,\alpha$, and $q$.

Furthermore, we obtain for $p\in [1,\infty]$ and $q = \infty$
\begin{equation*}
\begin{split}
\|f_0\|_{\wav,\alpha,s,p,\infty} 
\le (b-1)^{s} \max_{0\le L \le m} b^{L\alpha} 
{m-L+s-1 \choose s-1} \le C  b^{\alpha m},
\end{split}
\end{equation*}
where
$$C=C(b,s,\alpha) := (b-1)^s \sup_{0\le \nu < \infty} b^{-\nu \alpha} {\nu + s-1 \choose s-1}.$$

If we put $f_* = f_0 \|f_0\|_{\wav,\alpha,s,p,q}^{-1}$, then we get for all $p,q\in [1,\infty]$
\begin{equation*}
e^{\wor}(Q, \widetilde{\HH}_{\wav, \alpha, s,p,q}) \ge  e(Q, f_*) \ge \frac{m^{s-1}}{2(s-1)!\, C\, b^{\alpha m} \, m^{\frac{s-1}{q}}}
\ge C'\, b^{- \alpha m}\, m^{\frac{s-1}{q'}},
\end{equation*}
where $C'$ is a suitable constant only depending on $s,\alpha, b$, and $q$.
\end{proof}

\begin{remark}\label{Rem:Smolyak}
For integration on Haar wavelet spaces also Smolyak algorithms \cite{Smo63}, i.e., cubature rules whose set of integration nodes forms a sparse grid, have been rigorously analyzed, but their performance is clearly not as good as the one of QMC cubatures based on $(t,m,s)$-nets. Let us explain this in more detail: In the Hilbert space case where $p=2=q$, upper error bounds for deterministic Smolyak algorithms were provided in \cite{GLS07} and matching upper and lower error bounds for randomized Smolyak algorithms were derived \cite{WG23}.
In \cite{GLS07} actually so-called multiwavelet spaces $\mathcal{H}^s_{\alpha,n}$, with $s,n \in \N$, $\alpha > 1/2$,  based on the first $n$ Legendre polynomials, were studied. For base $b=2$ we have that $\HH_{\wav, \alpha, s,2,2} = \mathcal{H}^s_{\alpha,1}$. 
If $Q^{\rm Sm}_N$ denotes the Smolyak algorithm that relies, e.g., on a univariate iterated trapezoidal rule and uses altogether $N$ integration nodes, then it was shown in \cite[Corollary~4.7]{GLS07} that 
\begin{equation}\label{smo_det}
e^{\wor}(Q^{\rm Sm}_N, \HH_{\wav,\alpha,s,2,2}) 
= O \left( N^{-\alpha} \ln(N)^{(s-1)(\alpha + 1/2)} \right).
\end{equation}
The power of the logarithm in \eqref{smo_det} is obviously worse than the one in 
\eqref{nets} for $q=2=q'$. But since \cite{GLS07} provides solely a sharp general lower bound for arbitrary  cubature rules and not an additional matching lower bound specifically for Smolyak cubatures, it still could be that the upper bound in \eqref{smo_det} is simply too pessimistic. The known results for the randomized setting indicate that 
the latter is not the case. We know for the randomized error, i.e., the worst case root mean square error, that the following holds: In \cite{HHY04} it was shown 
that cubatures $Q^{\rm sc}_P$ based on randomly fully scrambled $(t,m,s)$-nets $P$ consisting of $N$ points 
exhibit the asymptotical behavior
\begin{equation}\label{nets_ran}
e^{\rm ran}(Q^{\rm sc}_P, \HH_{\wav,\alpha,s,2,2}) 
\asymp N^{-(\alpha + 1/2)},
\end{equation}
while it was shown in \cite{WG23} that randomized Smolyak algorithms $Q^{\rm raSm}_N$ with $N$ integration nodes, based on univariate stratified sampling (which is essentially the ``straightforward randomization'' of iterated trapezoidal rules), behave asymptotically like
\begin{equation}\label{smo_ran}
e^{\rm ran}(Q^{\rm raSm}_N, \HH_{\wav,\alpha,s,2,2}) 
\asymp N^{-(\alpha + 1/2)} \ln(N)^{(s-1)(\alpha + 1)}.
\end{equation}
For more details see \cite{GLS07, WG23} and the literature mentioned therein.
\end{remark}



\section{Spaces of Fractional  Smoothness}
\label{Sec:Fractional}


\subsection{The Univariate Case}
\label{Subsec:Frac_Univ}


The \emph{fractional derivative of order $\alpha > 0$ in the sense of Riemann and Liouville}
is given by
\begin{equation}\label{eq_fracderivative}
\frac{\rd^\alpha f(x)}{\rd x^\alpha} :=
\frac{1}{\Gamma(\beta-\alpha)} \frac{\rd^\beta}{\rd x^\beta}\int_0^x
f(t) (x-t)^{-\alpha+\beta-1}\,\rd t,
\end{equation}
where $\beta \in \NN$ satisfies $\beta-1 \le \alpha < \beta$. This means, in order to obtain a
derivative of fractional order, we first integrate the function with
fractional order $0 < \beta - \alpha \le 1$ and then differentiate
the resulting function $\beta$ times, where $\beta$ is an integer
(hence we
use the usual definition of differentiation for this part).

The differential operator  we use in this paper for
$\alpha > 0$ is given by (cf. \cite{D08})
\begin{equation*}
D^\alpha f(t) := \frac{\rd^\alpha(f(t) - f(0))}{\rd t^\alpha}.
\end{equation*}
Note that for $\alpha$ a positive integer $D^\alpha$ is just the
usual differential operator of order~$\alpha$. We use $f(t)-f(0)$
rather than just $f(t)$ in the definition above, since the fractional
derivative of a constant function is in general not $0$ as for
classical differentiation (see also the proof of
Lemma \ref{Lemma0}(iii) below).

Let $\alpha \in (0,1]$, $ \alpha^{-1} < p \le \infty$ and set
\begin{equation}\label{def_frac_H_alpha_p}
\HH_{\alpha,p} := \left\{f:[0,1]\rightarrow \real
\, \bigg|\,  \exists \tilde{f} \in
L_p([0,1])\, \exists c \in \real: f = c + \frac{1}{\Gamma(\alpha)}
\int_0^1 \tilde{f}(t) (\cdot-t)_+^{\alpha-1} \,\rd t \right\}.
\end{equation}
Here $(\cdot-t)_+^{\alpha-1}$ denotes the function that is
$(x-t)^{\alpha-1}$ for $x>t$ and $0$ otherwise.
Note that $(x-t)_+^{\alpha-1} \in L_q([0,1])$ for all $x$ as long as
$q(1-\alpha) < 1$.
Let $q = p'$.
Then the condition $q(1-\alpha)<1$ is equivalent to $p> \alpha^{-1}$.
Hence, $p>\alpha^{-1}$ and $\tilde{f} \in L_p([0,1])$ imply, due to H\"older's inequality,
$(x-\cdot)_+^{\alpha-1} \tilde{f}(\cdot) \in L_1([0,1])$.
Note that $\alpha^{-1} < p_1 \le p_2 \le \infty$ implies $\HH_{\alpha,p_2} \subseteq \HH_{\alpha,p_1}$.

In the case $\alpha = 1$, the function $(\cdot-t)_+^{\alpha-1}$ is 
nothing but the characteristic function $1_{(t,1]}$ of the half-open
interval $(t,1]$. 
Note that in this case  $\mathcal{H}_{\alpha,p}$ is also well-defined for $p=1=\alpha$.
Due to the fundamental theorem of calculus the space $\HH_{1,p}$ is the space of
\emph{absolutely continuous} functions on $[0,1]$ whose classical first derivative exists almost everywhere and belongs to $L_p([0,1])$.
A similar characterization of the spaces $\HH_{\alpha,p}$ can be 
given with the help of fractional differentiation in the sense of 
Riemann and Liouville.

\begin{lemma}\label{Lemma0}
Let $\alpha \in (0,1]$ and $\alpha^{-1} < p \le \infty$.
\begin{itemize}
\item[{\rm(}i\rm{)}] Each $f\in \HH_{\alpha,p}$ is
H\"older continuous with H\"older exponent $\alpha -1/p$.
\item[{\rm(}ii\rm{)}]
If $\tilde{f}\in L_p([0,1])$, $c\in\RR$, and
\begin{equation}
\label{rep}
f := c + \frac{1}{\Gamma(\alpha)} \int^1_0 \tilde{f}(t)
(\cdot - t)_+^{\alpha -1}\,dt,
\end{equation}
then for almost all $x\in [0,1]$ the fractional
derivative $D^{\alpha}f(x)$ exists
and is equal to $\tilde{f}(x)$. In particular, we have for all
$f\in\HH_{\alpha,p}$
\begin{equation}\label{eq_taylor}
f(x) = f(0) + \frac{1}{\Gamma(\alpha)} \int_0^1 D^\alpha f(t)
(x-t)_+^{\alpha-1} \,\rd t.
\end{equation}
\item[{\rm(}iii\rm{)}] There exist functions $f\in\HH_{\alpha,p}$ whose
fractional derivative $\frac{\rd^\alpha}{\rd x^\alpha} f$ is not
in $L_p([0,1])$.
\end{itemize}
\end{lemma}

\begin{proof}
Statement (i) of Lemma \ref{Lemma0} was verified in
\cite[Prop.~1]{D08}.

Statement (ii) follows from calculations in \cite[\S 2.2.1]{D08},
which we restate here briefly:
Observe
that $f(0) = c$.
We have
\begin{equation*}
\begin{split}
\int^x_0 (f(t)-f(0))(x-t)^{-\alpha}\,{\rm d}t
&= \frac{1}{\Gamma(\alpha)} \int^x_0\int^1_0 \tilde{f}(s)
(t-s)^{\alpha -1}_+\,{\rm d}s\, (x-t)^{-\alpha}\,{\rm d}t\\
&= \frac{1}{\Gamma(\alpha)} \int^x_0 \tilde{f}(s)\left( \int^x_s (t-s)^{\alpha-1}
(x-t)^{-\alpha}\,{\rm d}t \right) \,{\rm d}s\\
&= \Gamma(1-\alpha) \int^x_0 \tilde{f}(s)\,{\rm d}s,
\end{split}
\end{equation*}
where we used the well-known identity for the beta function (or Euler integral of the first kind)
\begin{equation}\label{beta-integral}
\int^x_s (t-s)^{\alpha-1} (x-t)^{-\alpha}\,\rd t
=\int^1_0 z^{\alpha-1}(1-z)^{-\alpha}\, \rd z
= \Gamma(\alpha)\Gamma(1-\alpha).
\end{equation}
Thus $D^\alpha f(t) = \frac{d}{dx} \int^x_0 \tilde{f}(s)\,\rd s$.
Since $\tilde{f}$ is in particular Lebesgue integrable on $[0,1]$,
the derivative $\frac{d}{dx} \int^x_0 \tilde{f}(s)\,\rd s$ exists
and is equal to $\tilde{f}(x)$ for almost all $x$, due to the fundamental
theorem of calculus and integration.

Let us now verify statement (iii): Let $f\equiv c$ for some $c\neq 0$. Then obviously $f\in \HH_{\alpha,p}$ and
\begin{equation*}
\frac{\rd^{\alpha}}{\rd x^\alpha} f(x) =
\frac{c}{\Gamma(1-\alpha)} \frac{\rd}{\rd x} \int^x_0 (x-t)^{-\alpha}
\,\rd t = \frac{c}{\Gamma(1-\alpha)} x^{-\alpha}.
\end{equation*}
Since $p>\alpha^{-1}$, necessarily
$\frac{\rd^{\alpha}}{\rd x^\alpha}f\notin L_p([0,1])$.
\end{proof}

Let us restate \cite[Lemma 1]{D08}.

\begin{lemma}\label{Note}
Let $\alpha \in (0, 1]$ and $\alpha^{-1} < p \le \infty$. Let
$g \in L_p([0,1])$ be
such that $\langle g, (x-\cdot)_+^{\alpha-1} \rangle = 0$
for all $x \in [0,1]$. Then  $g(x) =0$ for almost all
$x \in [0,1]$.
\end{lemma}

\begin{remark}\label{Rem:Proof_Gap}
The proof of Lemma \ref{Note} given in \cite{D08} actually contains
a gap. Let $\lambda$ be the Lebesgue measure on $[0,1]$.
For general measurable subsets $M_+$, $M_-$ of $[0,1]$ with
$M_+\cap M_-=\emptyset$ the condition
\begin{equation*}
\lambda(M_+ \cap [0,x)) >0
\hspace{2ex}\text{if and only if}\hspace{2ex}
\lambda(M_- \cap [0,x)) >0
\end{equation*}
for all $x\in [0,1]$ does not necessarily imply the condition
\begin{equation*}
\lambda(M_+ \cap B) >0
\hspace{2ex}\text{if and only if}\hspace{2ex}
\lambda(M_- \cap B) >0
\end{equation*}
for all Borel sets $B$ as stated in \cite{D08}.
A counterexample is provided by
\begin{equation*}
M_+ = \bigcup_{n\in\NN} \left[ \frac{1}{2n+1}, \frac{1}{2n} \right)
\hspace{2ex}\text{and}\hspace{2ex}
M_- = \bigcup_{n\in\NN} \left[ \frac{1}{2(n+1)}, \frac{1}{2n+1} \right).
\end{equation*}
Instead of fixing the proof approach from \cite{D08}, we 
may deduce Lemma~\ref{Note} from Lemma~\ref{Lemma0}: The statement is clear for $\alpha =1$ and $p\ge 1$. Let now $\alpha \in (0,1)$ and $p> \alpha^{-1}$. Define $f: [0,1] \to \RR$ by 
$$
f(x) := \frac{1}{\Gamma(\alpha)} \langle g, (x-\cdot)^{\alpha - 1}_+ \rangle.
$$
Hence $f$ is the zero function and due to Lemma~\ref{Lemma0}(ii) we have $g(x)= D^{\alpha} f(x) = 0 $ for almost all $x\in [0,1]$.
\end{remark}

For $p=2$, Lemma \ref{Note} and Parzen's integral representation theorem
\cite{P63} (see also \cite{K71}) imply that $\HH_{\alpha,2}$ is a
reproducing kernel Hilbert space with reproducing kernel
\begin{equation}\label{eq:rk_1}
K_\alpha(x,y) := 1 +  \int_0^1
(x-t)_+^{\alpha-1} (y-t)_+^{\alpha-1}\,\rd t
\end{equation}
and inner product
\begin{equation*}
\langle f, g\rangle_\alpha := f(0)g(0) + \frac{1}{\Gamma(\alpha)^2} \int_0^1 D^\alpha f(t)
D^\alpha g(t) \,\rd t.
\end{equation*}
Further note that for $\alpha = 1$ we obtain via this construction
the classical \emph{Sobolev space of order one}~$\HH_{1,2}$ on $[0,1]$, 
endowed with the \emph{anchored norm},
see \cite[Appendix~A.2.2]{NW08}.

\begin{definition}\label{Def:Space_Frac_Smooth_Uni}
For $p,q\in [1,\infty]$ and $\alpha\in (1/p,1]$ 
let the \emph{space of fractional smoothness} $\HH_{\alpha, p,q}$ 
be the vector space $\HH_{\alpha, p}$ 
defined in \eqref{def_frac_H_alpha_p}, 
equipped with the norm
\begin{equation*}
\| f \|_{\alpha,p,q} := \big( |f(0)|^{q} + \Gamma(\alpha)^{-q}
\|D^{\alpha}f\|^{q}_{L_p} \big)^{1/q}.
\end{equation*}
\end{definition}

Endowed with the norm above, $\HH_{\alpha,p,q}$ is obviously a complete normed space.
For $\alpha = 1$ and the choice $q=p$ we obtain 
the \emph{Sobolev space of order one}~$\HH_{1,p}$ on $[0,1]$,
endowed with the \emph{anchored norm},
cf., e.g., \cite{HS16}.

\begin{remark}\label{q_egal}
Since for fixed $\alpha$ and $p$ the norms $\|\cdot\|_{\alpha, p,q}$ are equivalent for all 
$1 \le q \le \infty$, we usually suppress the reference to $q$ and write simply $\HH_{\alpha,p}$  instead of $\HH_{\alpha,p,q}$. 
\end{remark}



\subsection{The Multivariate Case}
\label{Subsec:Frac_Multi}


Let $s\in \NN$, $\alpha>0$, and $p>\alpha^{-1}$.
We define the function space $\HH_{\alpha,s,p}$ as the vector space
\begin{equation}\label{def_frac_H_alpha_s_p}
\HH_{\alpha,s,p} := \{ f:[0,1]^s \to \RR \,|\,
\forall u \subseteq S \,\exists \tilde{f}_u \in
L_p([0,1]^{u}) : f =
\Phi((\tilde{f}_u)_{u\subseteq S}) \},
\end{equation}
where
\begin{equation}
\label{Phi}
\Phi((\tilde{f}_u)_{u\subseteq S})(\bsx) :=
\sum_{u\subseteq S} \Gamma(\alpha)^{-|u|}
\int_{[0,1]^{u}} \tilde{f}_u(\bst_u) \prod_{j\in u}
(x_j-t_j)^{\alpha-1}_+ \, \rd \bst_u.
\end{equation}
Here we used the convention that $L_p([0,1]^{\emptyset}) :=\RR$,
$\|\cdot\|_{L_p([0,1]^{\emptyset})}$ denotes the absolute value on $\RR$,
and
\begin{equation*}
\int_{[0,1]^{\emptyset}} \tilde{f}_\emptyset(\bst_\emptyset)
\prod_{j\in\emptyset} (x_j-t_j)^{\alpha-1}_+\,\rd \bst_\emptyset
:= \tilde{f}_\emptyset \in \RR.
\end{equation*}
Note that for $\alpha^{-1} < p_1 \le p_2 \le \infty$ we have $\HH_{\alpha,s,p_2} \subseteq \HH_{\alpha,s,p_1}$.

For $\bsx, \bsy\in [0,1]^s$ we write 
$(\bsx_u,\bsy)$ for $(z_1,\ldots,z_s)$
if $z_j=x_j$ for all $j\in u$ and $z_j = y_j$ otherwise.
For a function $f:[0,1]^s\to \RR$, $u\subseteq S$, and 
$\bst_u\in [0,1]^{u}$ we put
\begin{equation*}
f_{\forksnot, u}(\bst_u) := \sum_{v\subseteq u} (-1)^{|u\setminus v|}\,
f(\bst_v, \bszero).
\end{equation*}
The term $f_{\forksnot, u}$ is the $u$-th term of the
\emph{anchored decomposition} (also known as \emph{cut-HDMR}, cf. \cite{RA99}) 
$$
f=\sum_{u\subseteq S} f_{\forksnot, u}
$$ 
of
$f$ with anchor in $\bszero$, i.e., $f_{\forksnot, u}$ depends only on the variables $x_j$, $j\in u$, and  $f_{\forksnot, u}(\bsx) = 0$ if there exists  some $j\in u$ with $x_j=0$, see \cite{KSWW10}.
We define now for $u\subseteq S$ (formally) the fractional
derivative
\begin{equation}
\label{partfracder}
D^{(\bsalpha_u, \bszero_{S\setminus u})}f(\bst_u, \bszero) :=
\prod_{j\in u} \frac{\partial^{\alpha}}{\partial t^\alpha_j}
f_{\forksnot, u}(\bst_u),
\end{equation}
where
$\frac{\partial^{\alpha}}{\partial t_j^\alpha}$ is the $j$-th partial
derivative corresponding to the fractional derivative of order
$\alpha>0$ in (\ref{eq_fracderivative}).
Here we use the convention
\begin{equation*}
D^{(\bsalpha_\emptyset, \bszero_{S})}f(\bst_\emptyset,\bszero) 
:= f(\bszero).
\end{equation*}

In the proof of Lemma \ref{s-Derivative} (ii) we will see that any 
order of applying
the partial fractional derivatives in (\ref{partfracder}) to a
function $f\in \HH_{\alpha,s,p}$ leads to the same result. Hence an
analogon of the classical theorem of Schwarz for mixed classical
partial derivatives holds also for mixed fractional partial derivatives.

\begin{lemma}\label{s-Derivative}
Let $s\in \NN$, $\alpha\in (0,1]$, and $p>\alpha^{-1}$.
\begin{itemize}
\item[{\rm (}i{\rm )}]
Each $f\in \HH_{\alpha,s,p}$ is H\"older continuous with H\"older
exponent $\alpha-1/p$.
\item[{\rm (}ii{\rm )}]
If $f\in \HH_{\alpha,s,p}$ and
$f = \Phi((\tilde{f}_u)_{u\subseteq S})$, where
$\tilde{f}_u \in L_p([0,1]^{u})$ for all $u\subseteq S$,
then for each $u \subseteq S$ and almost all
$t_u\in [0,1]^{u}$ the fractional derivative
$D^{(\bsalpha_u, \bszero_{S\setminus u})}f(\bst_u, \bszero)$
exists and is equal to $\tilde{f}_u(t_u)$.
In particular, we have for all $f\in \HH_{\alpha,s,p}$ that
\begin{equation*}
f = \Phi((D^{(\bsalpha_u, \bszero_{S\setminus u})}
f(\cdot, \bszero))_{u\subseteq S}).
\end{equation*}
Moreover, for each $u\subseteq S$ the $u$-th term of the anchored decomposition
of $f$ with anchor in $\bszero$ is given by
\begin{equation}\label{anchored_decomp}
f_{\forksnot, u} = \frac{1}{\Gamma(\alpha)^{|u|}} \int_{[0,1]^u}  
D^{(\bsalpha_u, \bszero_{S\setminus u})} f(\bst_u,  \bszero) \prod_{j\in u} (x_j-t_j)_+^{\alpha -1}
\, {\rm d} \bst_u.
\end{equation}
\end{itemize}
\end{lemma}

\begin{proof}
Statement (i) of Lemma \ref{s-Derivative} 
was verified in \cite[Prop.~1]{D08}.

Hence it only remains to prove Statement (ii). Let
$f\in \HH_{\alpha,s,p}$ with $f=\Phi((\tilde{f}_u)_{u\subseteq S})$
for $\tilde{f}_u \in L_p([0,1]^{u})$, $u\subseteq S$. Then
\begin{equation*}
\begin{split}
\sum_{v\subseteq u} (-1)^{|u\setminus v|} f(\bsx_v, \bszero) =&
\sum_{v\subseteq u} (-1)^{|u\setminus v|} \left( \sum_{w\subseteq v}
\Gamma(\alpha)^{-|w|} \int_{[0,1]^{w}} \tilde{f}_w(\bst_w)
\prod_{j\in w}(x_j-t_j)^{\alpha-1}_+ \,{\rm d} \bst_w \right)\\
=& \sum_{w\subseteq u} \left( \sum_{w\subseteq v \subseteq u}
(-1)^{|u\setminus v|} \right)
\Gamma(\alpha)^{-|w|} \int_{[0,1]^{w}} \tilde{f}_w(\bst_w)
\prod_{j\in w}(x_j-t_j)^{\alpha-1}_+ \,{\rm d} \bst_w.
\end{split}
\end{equation*}
Since
\begin{equation*}
\sum_{w\subseteq v \subseteq u} (-1)^{|u \setminus v|} =
\sum_{v'\subseteq u\setminus w} (-1)^{|v'|} = (1-1)^{|u\setminus w|}
= \delta_{u,w},
\end{equation*}
we get
\begin{equation}\label{anchored_decomp_2}
f_{\forksnot, u}(\bsx_u) = \sum_{v\subseteq u} (-1)^{|u\setminus v|} f(\bsx_v, \bszero) =
\Gamma(\alpha)^{-|u|} \int_{[0,1]^{u}} \tilde{f}_u(\bst_u)
\prod_{j\in u}(x_j-t_j)^{\alpha-1}_+ \,{\rm d} \bst_u.
\end{equation}
Without loss of generality let us assume that $u=\{1,2,\ldots,\ell\}$
for some $\ell \in S$.
Then
\begin{equation*}
\begin{split}
&\prod_{j\in u} \frac{\partial^{\alpha}}{\partial x^\alpha_j}
\left( \sum_{v\subseteq u} (-1)^{|u\setminus v|} f(\bsx_v,\bszero)
\right)
= \frac{1}{\Gamma(1-\alpha)^{|u|}\Gamma(\alpha)^{|u|}}
\frac{\partial}{\partial x_1} \int^{x_1}_0 \,\times\\
&\left(
\ldots \frac{\partial}{\partial x_\ell} \int^{x_\ell}_0
\left( \int_{[0,1]^\ell} \tilde{f}_u(\bst_u) \prod_{j\in u}
(y_j-t_j)^{\alpha-1}_+ \,{\rm d} \bst_u \right) 
(x_\ell - y_\ell )^{-\alpha}
{\rm d} y_\ell  \ldots \right) (x_1-y_1)^{-\alpha}
\,{\rm d}y_1.
\end{split}
\end{equation*}
Using Fubini's theorem and the identity for the beta function \eqref{beta-integral} we get
\begin{equation*}
\begin{split}
&\frac{\partial}{\partial x_\ell} \int^{x_\ell}_0
\int_{[0,1]^\ell} \tilde{f}_u(\bst_u) \prod_{j=1}^\ell
(y_j-t_j)^{\alpha-1}_+ \,\bigotimes^\ell_{j=1} {\rm d}t_j
(x_\ell -y_\ell)^{-\alpha} \,{\rm d}y_\ell\\
= &\frac{\partial}{\partial x_\ell} \int_{[0,1]^{\ell-1}}
\prod^{\ell-1}_{j=1} (y_j-t_j)^{\alpha-1}_+
\left[ \int^{x_\ell}_0 \tilde{f}_u(\bst_u) \int^{x_\ell}_{t_\ell}
(y_\ell-t_\ell)^{\alpha-1}(x_\ell-y_\ell)^{-\alpha}\,{\rm d}y_\ell
\,{\rm d}t_\ell \right] \bigotimes^{\ell-1}_{j=1} {\rm d}t_j\\
= &\Gamma(\alpha)\Gamma(1-\alpha) \frac{\partial}{\partial x_\ell}
\int^{x_\ell}_0 \int_{[0,1]^{\ell-1}} \tilde{f}_u(\bst_u)
\prod^{\ell-1}_{j=1}(y_j-t_j)^{\alpha-1}_+
\,\bigotimes^\ell_{j=1} {\rm d}t_j\\
= &\Gamma(\alpha)\Gamma(1-\alpha) \int_{[0,1]^{\ell-1}} \tilde{f}_u
(t_1,\ldots,t_{\ell-1},x_\ell) \prod^{\ell-1}_{j=1}
(y_j-t_j)^{\alpha-1}_+ \,\bigotimes^{\ell-1}_{j=1}{\rm d}t_j
\end{split}
\end{equation*}
for almost all $x_\ell \in [0,1]$.
From this it is easy to see that
\begin{equation*}
\prod_{j\in u} \frac{\partial^{\alpha}}{\partial x^\alpha_j}
\left( \sum_{v\subseteq u} (-1)^{|u\setminus v|}
f(\bsx_v,\bszero) \right) = \tilde{f}_u(x_1,\ldots,x_\ell)
=\tilde{f}_u(\bsx_u)
\end{equation*}
for almost all $\bsx_u \in [0,1]^u$. Now \eqref{anchored_decomp} follows from \eqref{anchored_decomp_2}.
\end{proof}

\begin{definition}\label{Def:Space_Frac_Smooth_Mult}
For $p,q\in [1,\infty]$ and $\alpha\in (1/p,1]$  let the \emph{space of fractional smoothness}
$\HH_{\alpha,s,p,q}$ be the vector space $\HH_{\alpha,s,p}$ defined as in \eqref{def_frac_H_alpha_s_p},
endowed with the
norm
\begin{equation*}
\|f\|_{\alpha, s,p,q} := \left( \sum_{u\subseteq S}
\Gamma(\alpha)^{-q|u|} \|\tilde{f}_u\|^{q}_{L_p}
\right)^{1/q}.
\end{equation*}
\end{definition}

Note that the norm is well-defined due to Lemma~\ref{s-Derivative} (ii). 
It is easily seen that $\HH_{\alpha,s,p,q}$
is a complete normed space. 

\begin{remark}\label{s:q_egal}
Note that different choices of $q$ lead to equivalent norms
on the vector space $\HH_{\alpha,s,p}$.  Therefore we often suppress the reference to $q$ and write simply $\HH_{\alpha,s,p}$  instead of $\HH_{\alpha,s,p,q}$.
\end{remark}

Let us close the subsection by considering the case $p=2=q$.
In this case $\HH_{\alpha,s,2,2}$ is a Hilbert space with scalar
product
\begin{equation*}
\langle f,g \rangle_{\alpha,s,2,2} =
\sum_{u\subseteq S} \Gamma(\alpha)^{-2|u|}
\langle \tilde{f}_u,\tilde{g}_u \rangle,
\end{equation*}
where $f= \Phi((\tilde{f}_u)_{u\subseteq S})$,
$g = \Phi((\tilde{g}_u)_{u\subseteq S})$,
and $\Phi$ as defined in (\ref{Phi}).
It is easy to check that $\HH_{\alpha,s,2,2}$ is the complete $s$-fold
tensor product Hilbert space of $\HH_{\alpha,2,2}$.
Moreover, $\HH_{\alpha,s,2,2}$ is a reproducing kernel Hilbert space
and its reproducing kernel $K_{\alpha,s}$ is the $s$-fold product
of the reproducing kernel $K_\alpha$ of $\HH_{\alpha,2,2}$ from \eqref{eq:rk_1}, i.e.,
\begin{equation*}
K_{\alpha,s}(\bsx,\bsy) = \prod^s_{j=1} K_{\alpha}(x_j,y_j)
\hspace{2ex}\text{for all $\bsx, \bsy\in [0,1]^s$,}
\end{equation*}
cf. \cite{Aro50}.
In case $p \neq 2$ the situation is similar but more complicated.
Tensor products of general Banach spaces is a well-studied but not simple subject,
we refer to \cite{LiCh} and \cite{DF} as well as to our Appendix~\ref{Sec:Appendix}. 
We now present an extension of the tensor product result for $p=2=q$ to the case where $p=q$ and $1 < p < \infty$.
As before, we denote by $\delta_p$ the $p$-nuclear tensor norm, see Definition~\ref{tensdefi},
and use the convention $\HH_{\alpha,1,p} := \HH_{\alpha,p}$.

\begin{theorem}\label{tensor3}
 Let $1 <p<\infty$, $\alpha \in (1/p,1]$ and $s\in \N$, $s\ge 2$.
 Then
 \[
 \HH_{\alpha,s-1,p}  \otimes_{\delta_p}\HH_{\alpha,p} =  \HH_{\alpha,s,p} 
 \]
 in the sense of equivalent norms.
\end{theorem}

The proof can be found in Subsection~\ref{APP:TPSO} of the appendix on tensor products of Banach spaces.


\subsection{A Sharp Fractional Koksma-Hlawka Inequality}
\label{Subsec:SKHI}


For a real number $0<\alpha\le 1$, a (multi-)set $P=\{\bsx_0,\ldots,
\bsx_{N-1}\}$ in $[0,1)^s$ with $\bsx_n = (x_{n,1},\ldots,x_{n,s})$,
a set $\emptyset\neq u \subseteq S$, and a point
$\bst_u = (t_j)_{j\in u}\in [0,1]^{u}$ we define the
\emph{fractional discrepancy function} $\Delta_\alpha$ by
\begin{equation*}
\Delta_\alpha(\bst_u, u, P) := \alpha^{-|u|}
\prod_{j\in u}(1-t_j)^\alpha - \frac{1}{N} \sum^{N-1}_{n=0}
\prod_{j\in u}(x_{n,j}-t_j)_+^{\alpha-1}.
\end{equation*}
Obviously, $\Delta_\alpha(\cdot, u, P)$ is a measurable function.
If $\alpha=1$, then $\Delta_1(\cdot, u, P)
\in L_\infty ([0,1]^{u})$ and
$\|\Delta_1(\cdot,u,P)\|_{L_q} >0$ for all $q\in [1,\infty]$.
If $\alpha\in (0,1)$ and $q<(1-\alpha)^{-1}$, then we have $\Delta_\alpha(\cdot, u, P)
\in L_q([0,1]^{u})$ and
$\|\Delta_\alpha(\cdot,u,P)\|_{L_q} >0$.
Let us use the short hand
\begin{equation*}
I_s(f):= \int_{[0,1]^s} f(\bst)\,{\rm d}\bst.
\end{equation*}

The following fractional Koksma-Hlawka inequality was shown in \cite{D08} for $\alpha \in (0,1)$;
in the classical case $\alpha =1$ the result was known before, see, e.g., \cite{Hic98}.

\begin{theorem}\label{th_fracKH}
Let $s \ge 1$, $0 < \alpha \le 1$, $p > \alpha^{-1}$, and $q \ge 1$.
Then for any $f \in
\HH_{\alpha,s,p}$ we have
\begin{equation*}
\left|I_s(f) - \frac{1}{N}
\sum_{n=0}^{N-1} f(\bsx_n) \right|  \le  D^\ast_{\alpha,s,p',q'}(P)
V_{\alpha,s,p,q}(f),
\end{equation*}
where
\begin{equation*}
D^\ast_{\alpha,s,p',q'}(P) = \left(\sum_{\emptyset \neq u \subseteq S}
 \left(\int_{[0,1]^{u}}
|\Delta_\alpha(\bst_u,u,P)|^{p'} \,\rd\bst_u \right)^{q'/p'}
\right)^{1/q'}
\end{equation*}
and
\begin{equation*}
V_{\alpha,s,p,q}(f) =\left(\sum_{\emptyset \neq u \subseteq S} 
\Gamma(\alpha)^{-q|u|}
\left(\int_{[0,1]^{u}} \left| D^{(\bsalpha_u, \bszero_{S\setminus u})}f(\bst_u, \bszero)\right|^{p} \,\rd\bst_u
\right)^{q/p} \right)^{1/q}
\end{equation*}
with the obvious modifications in the case where $p,q$ or $q'=\infty$.
\end{theorem}

\begin{remark}
Note that
$V_{\alpha,s,p,q}(\cdot)$ is just the seminorm in the space
$\HH_{\alpha,s,p}$ and could be replaced by the norm
$\|\cdot\|_{\alpha,s,p,q}$.
Furthermore, for $\alpha\in (0,1)$ the condition $p > \alpha^{-1}$ implies $p' < (1-\alpha)^{-1}$, which
is also the condition required to guarantee that
$D^\ast_{\alpha,s,p',q'}(P) < \infty$. 
We refer to $D^\ast_{\alpha,s,p',q'}(P)$ as \emph{fractional $(p', q')$-discrepancy} 
of the point set $P$.
\end{remark}

Here we want to show that the fractional Koksma-Hlawka inequality is indeed
sharp. Let us therefore denote by $e^{\wor}(Q_P,\HH_{\alpha,s,p,q})$ the
worst-case error
of the quasi-Monte Carlo cubature $Q_P$ based on the set of
integration points $P$ on $\HH_{\alpha,s,p,q}$, which is analogously defined as in (\ref{wce}), i.e.,
\begin{equation*}
e^{\wor}(Q_P,\HH_{\alpha,s,p,q}) = \sup_{f\in \HH_{\alpha,s,p}; \|f\|_{\alpha,s,p,q}\le 1} \left| \int_{[0,1]^s} f(\bsx) \, \rd \bsx
- \frac{1}{N} \sum_{n=0}^{N-1} f(\bsx_n) \right|.
\end{equation*}
Then we have in fact the following theorem:

\begin{theorem}
\label{KoksmaHlawka}
Let $s\in\NN$, $0<\alpha\le 1$, $p>\alpha^{-1}$, and $q\in [1,\infty]$.
For all point sets $P\subset [0,1]^s$ we have
\begin{equation}
\label{koksmahlawka}
e^{\wor}(Q_P,\HH_{\alpha,s,p,q}) =  D^\ast_{\alpha,s,p',q'}(P).
\end{equation}
\end{theorem}

Theorem~\ref{KoksmaHlawka} leads to the following insight: 
If we rely on the fractional discrepancy as quality criterion for sample points, then the spaces of fractional smoothness are the appropriate function spaces for the analysis of QMC cubature rules.

\begin{proof}
Due to Theorem \ref{th_fracKH} we have 
$e^{\wor}(Q_P,\HH_{\alpha,s,p,q}) \le D^\ast_{\alpha,s,p',q'}(P)$, 
so it only remains to show that
$e^{\wor}(Q_P,\HH_{\alpha,s,p,q}) \ge  D^\ast_{\alpha,s,p',q'}(P)$.
To this purpose we use the well-known fact that in any Banach space $(X, \|\cdot\|_X)$
with dual space $(X', \|\cdot\|_{X'})$ we have
$\|x\|_X = \sup\{|x'(x)|: x'\in X', \|x'\|_{X'}=1\}$ and that due to
the Fischer-Riesz theorem $L_p([0,1]^s)$ is the dual space of
$L_{p'}([0,1]^s)$
with respect to the dual bracket $\langle \cdot, \cdot \rangle$. (Note
that $p'<\infty$, as $p > \alpha^{-1} \ge 1$.) Furthermore, the space $\ell_{q}^{2^s-1}$ of finite sequences of length $2^s-1$ is the dual space of $\ell_{q'}^{2^s-1}$ (even if $q'=\infty$).

For $\emptyset \neq u \subseteq S$ the function 
$\Delta_\alpha(\cdot,u,P)$
is in $L_{p'}([0,1]^{u})$. Then for given $\varepsilon >0$ we can
choose an
$\tilde{f}^{(\varepsilon)}_u\in L_p([0,1]^{u})$ such that
\begin{equation*}
\begin{split}
\int_{[0,1]^{u}} \tilde{f}^{(\varepsilon)}_u(\bst_u)
\Delta_\alpha(\bst_u,u,P) \,\rd \bst_u
 =& \left|\int_{[0,1]^{u}} \tilde{f}^{(\varepsilon)}_u(\bst_u)
\Delta_\alpha(\bst_u,u,P) \,\rd \bst_u
\right|\\
\ge& \|\tilde{f}^{(\varepsilon)}_u\|_{L_p}
\big( \|\Delta_\alpha(\bst_u,u,P)\|_{L_{p'}} - \varepsilon \big).
\end{split}
\end{equation*}
Let us put
\begin{equation*}
 f^{(\varepsilon)}(\bsx) := \sum_{\emptyset \neq u\subseteq S}
\Gamma(\alpha)^{-|u|}
\int_{[0,1]^{u}} \tilde{f}^{(\varepsilon)}_u(\bst_u) \prod_{j\in u} (x_j-t_j)^{\alpha-1}_+
\, \rd \bst_u.
\end{equation*}
Thus we have
$$\|f^{(\varepsilon)}\|_{\alpha,s,p,q} =
V_{\alpha,s,p,q}(f^{(\varepsilon)})
= \left( \sum_{\emptyset \neq u \subseteq S} \Gamma(\alpha)^{-q|u|} \|\tilde{f}^{(\varepsilon)}_u\|^{q}_{L_p} \right)^{1/q}.$$
Without loss of generality we can assume that $\|f^{(\varepsilon)}\|_{\alpha,s,p,q} = 1$ and
\begin{equation*}
\begin{split}
\sum_{\emptyset \neq u \subseteq S} \Gamma(\alpha)^{-|u|}
\|\tilde{f}^{(\varepsilon)}_u\|_{L_p}
\big( \|\Delta_\alpha(\bst_u,u,P)\|_{L_{p'}} -\varepsilon \big)& \ge \\
\left( \sum_{\emptyset \neq u \subseteq S}
\big( \|\Delta_\alpha(\bst_u,u,P)\|_{L_{p'}}-\varepsilon \big)^{q'}
\right)^{1/q'} -\varepsilon&
\end{split}
\end{equation*}
hold (by multipling each $\tilde{f}^{(\varepsilon)}_u$ with an
individual positive
factor if necessary).
Altogether we have
\begin{equation*}
 \begin{split}
\bigg| \int_{[0,1]^s} f^{(\varepsilon)}(\bst)\,\rd \bst -
\frac{1}{N} \sum_{n=0}^{N-1} f^{(\varepsilon)}&(\bsx_n) \bigg|
=
\sum_{\emptyset \neq u \subseteq S} \Gamma(\alpha)^{-|u|}
\left|  \int_{[0,1]^{u}} \tilde{f}^{(\varepsilon)}_u(\bst_u)
\Delta_\alpha(\bst_u,u,P) \,\rd \bst_u \right|\\
&\ge
\sum_{\emptyset \neq u \subseteq S} \Gamma(\alpha)^{-|u|}
\|\tilde{f}^{(\varepsilon)}_u\|_{L_p}
\big( \|\Delta_\alpha(\bst_u,u,P)\|_{L_{p'}} - \varepsilon \big)\\
&\ge
\left( \sum_{\emptyset \neq u \subseteq S}
\big( \|\Delta_\alpha(\bst_u,u,P)\|_{L_{p'}}-\varepsilon \big)^{q'}
\right)^{1/q'} -\varepsilon.
 \end{split}
\end{equation*}
Since $\varepsilon$ can be chosen arbitrarily small, we get
$e^{\wor}(Q_P,\HH_{\alpha,s,p,q}) \ge  D^\ast_{\alpha,s,p',q'}(P)$.
\end{proof}


\section{Besov Spaces and Bessel Potential spaces}
\label{Sec:Besov}



\subsection{The Univariate Case}
\label{Subsec:Besov_Univ}



\subsubsection{Preliminaries}


\paragraph{Besov Spaces.}

For a function $f \in L_p([0,1])$, $1\le p <
\infty$, and $0 < \delta < 1$ we define the \emph{integral modulus of
continuity} by
\begin{equation*}
\omega_p(f,\delta) :=  \sup_{h \le \delta} \left(\int_0^{1-h}
|f(x+h)-f(x)|^p \rd x \right)^{1/p},
\end{equation*}
and for $p = \infty$ we define the modulus of continuity for $f \in
L_\infty([0,1])$ by
\begin{equation*}
\omega_\infty(f,\delta) := \sup_{h \le \delta}  \esssup_{0 \le
x \le 1-h} |f(x+h)-f(x)|.
\end{equation*}

Let $\alpha\in (0,1)$ and $1 \le p, q \le \infty$.
Let us define the semi-norm
\begin{equation*}
|f|_{B^\alpha_{p,q}} :=  \left(\int_0^1 (t^{-\alpha} \omega_p(f,t))^q
\,\frac{\rd t}{t} \right)^{1/q}
\end{equation*}
for $q\in [1,\infty)$ and
\begin{equation*}
|f|_{B^\alpha_{p,\infty}} :=  \esssup_{0<t\le 1} t^{-\alpha}
\omega_p(f,t)
\end{equation*}
for $q=\infty$.

We say that a locally integrable 
 function $f:[0,1]\rightarrow \real$ is in the Besov space 
 $B^\alpha_{p,q} = B^\alpha_{p,q} ([0,1])$
whenever its norm
\begin{equation*}
\|f\|_{B^\alpha_{p,q}} :=   \|f\|_{L_p}  +
|f|_{B^\alpha_{p,q}}
\end{equation*}
is finite.

\begin{remark}
 \rm
There is a rich literature about Besov spaces. We refer, e.g.,  to Nikol'skij \cite{Ni}, Peetre \cite{Pee76}
and Triebel \cite{Tri83, Tri92, Tri06}.
\end{remark}

Sometimes it is convenient to consider a different norm
on $B^\alpha_{p,q}$. Let  $\rho>1$ and put
\begin{equation*}
|f|^{(\rho)}_{B^\alpha_{p,q}} := \left( \sum^\infty_{j=0}
\rho^{j\alpha q}
\omega_p(f, \rho^{-j})^q \right)^{1/q}
\end{equation*}
for $q\in [1,\infty)$ and
\begin{equation*}
|f|^{(\rho)}_{B^\alpha_{p,\infty}} :=
\sup_{j \in \mathbb{N}_0} \rho^{j\alpha} \omega_p(f,\rho^{-j})
\end{equation*}
for $q=\infty$. Then
$\|\cdot\|^{(\rho)}_{B^\alpha_{p,q}}$, given by
\begin{equation*}
\|f\|^{(\rho)}_{B^\alpha_{p,q}} := \|f\|_{L_p}
+ |f|^{(\rho)}_{B^\alpha_{p,q}},
\end{equation*}
is a norm on $B^\alpha_{p,q}$. 

\begin{lemma}
\label{Normeq}
Let $\rho>1$, $p,q\in [1,\infty]$, and $\alpha \in (0,1)$. Then the norms $\|\cdot\|_{B^\alpha_{p,q}}$ and
$\|\cdot\|^{(\rho)}_{B^\alpha_{p,q}}$ on $B^\alpha_{p,q} ([0,1])$
are equivalent.
\end{lemma}

In case $\rho = 2$ the statement of Lemma~9 is standard, we refer, e.g., to  \cite[Section~2]{DP88}, \cite[(4.6)]{Dah96} or 
 \cite[Thm.~2.5.12]{Tri83}. 
For general $\rho >1$ it becomes a consequence of  more general Fourier-analytic descriptions, see \cite[2.2.1]{Tri78}.
Alternatively one can also generalize the method used in \cite[Section~2]{DP88}.
For the convenience of the reader we provide a detailed proof in the spirit of this last reference.

\begin{proof}
Let $f\in L_p([0,1])$. Then $\omega_p(f,\cdot)$ is a monotonic
increasing function and bounded by $2\|f\|_{L_p}$. Now let
$q\in [1,\infty)$. Then
\begin{equation*}
\begin{split}
\int^1_0 [t^{-\alpha}\omega_p(f,t)]^q \,\frac{\rd t}{t}
&= \sum^\infty_{j=0} \int^{\rho^{-j}}_{\rho^{-j-1}}
[t^{-\alpha}\omega_p(f,t)]^q \,\frac{\rd t}{t}\\
&\le \sum^\infty_{j=0}
(\rho^{-j} - \rho^{-j-1})[\rho^{(j+1)\alpha}\omega_p(f,\rho^{-j})]^q
\rho^{j+1}\\
&\le \rho^{\alpha q+1} \sum^\infty_{j=0} [\rho^{j\alpha}
\omega_p(f, \rho^{-j})]^q.
\end{split}
\end{equation*}
Furthermore,
\begin{equation*}
\begin{split}
\int^1_0 [t^{-\alpha}\omega_p(f,t)]^q \,\frac{\rd t}{t}
&\ge \sum^\infty_{j=0}
(\rho^{-j} - \rho^{-j-1})[\rho^{j\alpha}\omega_p(f,\rho^{-j-1})]^q
\rho^{j}\\
&= \left( 1 - \frac{1}{\rho} \right) \rho^{-\alpha q}
\sum^\infty_{j=1} [\rho^{j\alpha}\omega_p(f,\rho^{-j})]^q \\
\end{split}
\end{equation*}
and
\begin{equation*}
[\rho^{j\alpha} \omega_p(f,\rho^{-j})]^q|_{j=0} = [\omega_p(f,1)]^q
\le 2^q\|f\|^q_{L_p}.
\end{equation*}
In the case where $q=\infty$, we see easily that
\begin{equation*}
|f|^{(\rho)}_{B^\alpha_{p,\infty}} \le
|f|_{B^\alpha_{p,\infty}} \le
\rho^{\alpha}|f|^{(\rho)}_{B^\alpha_{p,\infty}}.
\end{equation*}
From these inequalities follows that the norms
$\|\cdot\|_{B^\alpha_{p,q}}$ and
$\|\cdot\|^{(\rho)}_{B^\alpha_{p,q}}$ are equivalent.
\end{proof}

\paragraph{Bessel Potential Spaces.}

For the next theorem we need to introduce \emph{Bessel potential spaces}, 
which are also known as \emph{(fractional) Sobolev spaces}.
On $\real$ they are introduced by means of Bessel potentials.
Let $\alpha \in \RR$ and $1<p<\infty$. 
The tempered distribution $f\in \mathcal{S}' (\real)$
belongs to $H^\alpha _p (\real)$ if $\cfi [(1+|\xi|^2)^{\alpha/2} \cf f(\xi)]$ is in $ L_p (\real)$;
here $\cf$, $\cfi$ denote the Fourier transform and its inverse on $\mathcal{S}'(\RR)$.
It is well-known that for $\alpha \ge 0$ we have $H^\alpha _p (\real) \subseteq  L_p (\real)$. 
The norm on $H^\alpha _p (\real)$ is given by
\[
 \|f\|_{H^\alpha _p (\real)}:= \|\cfi [(1+|\xi|^2)^{\alpha/2} \cf f(\xi)](\, \cdot\, )\|_{L_p (\real)}.
\]
On the interval $[0,1]$ we could introduce the space $H^\alpha_p ([0,1])$ in at least two ways
which are fortunately equivalent.
The first variant
is based on 
restrictions.
Let $\alpha \ge 0$ and $1< p< \infty$.
For a function $g \in L_p (\real)$ we denote by $g_{|_{[0,1]}}$ the restriction of $g$ to the interval $[0,1]$. 
A function $f \in L_p ([0,1])$ belongs to $H^\alpha _p ([0,1])$
if there exists a function $g\in H^\alpha _p (\RR)$ such that $f=g|_{[0,1]}$, and the norm on
$H^\alpha _p ([0,1])$ is given by
\[
 \|f\|_{H^\alpha _p ([0,1])}:= \inf \Big\{ \|g\|_{H_p^\alpha (\real)}: ~ g \in H^\alpha_p (\real) \quad \mbox{and}\quad f=g_{|_{[0,1]}}\Big\}\, .
\]
There is another characterization by using differences. Let $\alpha \in (0,1)$ and $1<p<\infty$.
A function $f \in L_p ([0,1])$ belongs to $H^\alpha _p ([0,1])$
if 
\beq\label{neueq2}
&& \hspace*{-0.7cm}
\| f \|^*_{H^\alpha_p ([0,1])} :=  
\| f \|_{L_p ([0,1])}
\\
\nonumber
& + & \bigg(\int_{0}^1 \bigg[ \int_0^{\min(x,1-x)}   t^{-2\alpha} 
\bigg( \frac 1t \int_{\{h: ~|h|<t,~x+h \in (0,1)\}} |f(x+h)-f(x)|\, dh\bigg)^2 \frac{dt}{t}\bigg]^{p/2} dx\bigg)^{1/p}  
\eeq
is finite.
Furthermore, $\| \, \cdot \,  \|^*_{H^\alpha_p ([0,1])}$ and $\| \, \cdot \,  \|_{H^\alpha_p ([0,1])}$ are equivalent on 
$H^\alpha _p ([0,1])$. We refer to \cite[Thm.~4.10]{Tri06}. 
From \eqref{neueq2} it becomes evident that
\begin{equation}\label{monotonicity_bessel_potential_spaces}
H^\alpha _{p_2} ([0,1]) \hookrightarrow H^\alpha _{p_1} ([0,1])
\hspace{3ex}\text{for all $1<p_1 < p_2 < \infty$.}
\end{equation}

\subsubsection{Embedding Results} 

We start by studying embeddings of spaces of fractional smoothness into 
Bessel potential spaces.

\begin{theorem}\label{Winfried}
Let $\alpha\in (0,1)$ and $\alpha^{-1} < p < \infty$.
Then we have the continuous embedding
\begin{equation*}
\HH_{\alpha,p} \hookrightarrow H^\alpha_{p} ([0,1]).
\end{equation*}
\end{theorem}


\begin{proof}
{\em Step 1.} Tools from Fourier  analysis.
For $g \in C_0^\infty (\re)$ we consider  
\[
\ce_\alpha g(x):= \frac{1}{\Gamma (\alpha)} \int_{-\infty}^\infty \frac{g(t)}{(x-t)_+^{1-\alpha}}\, dt\, , \qquad x \in \re \, .
\]
Of course,  $\ce_\alpha  g$ is a convolution of $g$ with the kernel 
\[
k_\alpha (t):= \frac{1}{\Gamma (\alpha)} \, t_+^{-1+\alpha}\, , \qquad t \in \re \, .
\]
The  Fourier transform of $k_\alpha $ is given by 
\begin{equation}\label{1}
\cf k_\alpha (\xi ) := \frac{1}{\sqrt{2\pi}} \int_{-\infty}^\infty e^{-it\xi}\, k_\alpha(t)\, dt
= \frac{1}{\sqrt{2\pi}} \, |\xi|^{-\alpha} \, \left\{
\begin{array}{lll}
e^{-i \pi \alpha/2} & \qquad & \mbox{if}\quad \xi>0\, , 
\\
e^{i \pi \alpha/2} & \qquad & \mbox{if}\quad \xi<0\, , 
\end{array}
\right.
\end{equation}
see, e.g., \cite[page 400]{BN}. 
Further we need 
\begin{equation}\label{2}
\cf |t|^{-s}(\xi) = \sqrt{\frac 2\pi} \, \frac{\Gamma (1-s)}{|\xi|^{1-s}} \, \sin \frac{\pi s}{2}\, , \qquad \xi \in \re \, , s\in (0,1)\,,
\end{equation}
see \cite[page 516]{BN}. We want to apply the Hilbert transform.
Let
\[
Hg (x) := PV \Big[\frac{1}{\pi} \int_{-\infty}^\infty \frac{g(x-t)}{t}\, dt
\Big]\, , \qquad x \in \re\, , 
\]
see \cite[page 306]{BN}.
The  Fourier transform of the  Hilbert transform of a function $g$ is given by 
\begin{equation}\label{3}
\cf H g(\xi)= -i \, (\sign \xi) \, \cf g (\xi)\, .
\end{equation}
Moreover, if $1 <p< \infty$, then  
\begin{equation}\label{4}
\| \, Hg \,\|_{L_p} \le c_H \, \| \, g \, \|_{L_p } 
\end{equation}
holds for all  $g \in L_p (\re)$, where $c_{H}$ does not depend on $g$, we refer to \cite[Thm.~8.1.12]{BN}.
\\
{\em Step 2.} Homogeneous function spaces. For a while we have to calculate modulo polynomials. 
With $\cp$ we denote the set of all polynomials. 
Let $f$ be a tempered distribution, i.e., $f \in \cs' (\re)$. Then we put
\[
[f] := \{g: \quad \exists p \in \cp \quad \mbox{mit}\quad f+p =g\}\, .
\]
For $s \in \re $ and $1 <p<\infty$, let  $\dot{H}^s_p (\re)$ denote the collection of all equivalence classes  $[f] \in \cs'(\re)/\cp$ such that 
\[
\|\, [f]\, \|_{\dot{H}^s_p(\RR)} := \|\, \cfi [|\xi|^s \, \cf f (\xi)](\,\cdot \, )\|_{L_p}
<\infty \, .
\]
$\|\, [f]\,\|_{ \dot{H}^s_p (\re)}$ is a norm and $\dot{H}^s_p (\re)$  becomes a Banach space.
It is well-known that
\[
\dot{H}^s_p (\re) \cap L_p (\re) = H^s_p (\re)
\]
by using an appropriate interpretation (if the equivalence class  $[f]$ contains 
an element $g$ which  belongs to $L_p(\re)$, then this particular function is in ${H}^s_p (\re)$). 
For all this we refer to \cite[Chapter 6]{BL}.
\\
{\em Step 3.} The regularity of $\ce_\alpha g$.
Let $1<p<\infty$ and $0 < \alpha < 1$.
We claim that
$g \in L_p (\re)$ implies $[\ce_\alpha g] \in \dot{H}^{\alpha}_p(\re)$.
\\
Temporarily we assume $g \in C_0^\infty (\re)$.
Concerning the  Fourier transform of $\ce_\alpha g$ we know
\[
\cf \ce_\alpha g (\xi) = \sqrt{2 \pi} \, \cf k_\alpha (\xi ) \, \cf g(\xi)\, .
\]
Hence
\[
\cfi [|\xi|^{\alpha} \, \cf \ce_\alpha  g (\xi)](x) =   
\cfi [\sign^* (\xi)\, \cf g (\xi)](x)\, ,
\]
where
\[
\sign^*(\xi):= \left\{
\begin{array}{lll}
e^{-i \pi \alpha/2} & \qquad & \mbox{if}\quad \xi>0\, , 
\\
e^{i \pi \alpha/2} & \qquad & \mbox{if}\quad \xi<0\, .
\end{array}
\right.
\]
Because of  
\[
\sign^*(\xi) = \gamma \, \sign \xi + \beta
\]
where
\[
\beta := \cos \frac{\pi \alpha}{2} \qquad \mbox{and}
\qquad \gamma := -i \, \sin \frac{\pi \alpha}{2} 
\]
we get the identity
\[
\cfi [|\xi|^{\alpha} \, \cf \ce_\alpha  g (\xi)](x) =    \sin \frac{\pi \alpha}{2} Hg(x) + \cos \frac{\pi \alpha}{2}\, g(x)\, .
\]
Applying (\ref{4}), we conclude
\[
\| \, \cfi [|\xi|^{\alpha} \, \cf \ce_\alpha  g (\xi)](x)\, \|_{L_p} \le 
(c_H +1) \, \| \, g\, \|_{L_p}\, .
\]
A density argument completes the proof of our claim.
\\
{\em Step 4.} 
Let $f \in L_p (0,1)$. By $\tilde{f}$ we denote the extension  of $f$ by zero.
Hence $\tilde{f} \in L_p (\re)$.
Applying Step 3 we conclude
$[\ce_\alpha \tilde{f}] \in \dot{H}^{\alpha}_p(\re)$.
Observe
\[
\ce_\alpha \tilde{f} (x)= \frac{1}{\Gamma (\alpha)} \int_{-\infty}^\infty 
\tilde{f}(t) \, (x-t)_+^{\alpha-1}\, dt
= \frac{1}{\Gamma (\alpha)} \int_{0}^x f(t) \, (x-t)^{\alpha-1}\, dt
\, , \qquad x \in [0,1] \, .
\] 
By denoting $\ce_\alpha \tilde{f}|_{[0,1]}$ again by $\ce_\alpha \tilde{f}$, we obtain
from Lemma~\ref{Lemma0}(i)
in particular that $\ce_\alpha \tilde{f} \in L_p (0,1)$.
\\
{\em Step 5.} More about function spaces.
Following \cite[Thm.~2.5.11]{Tri83} an equivalence  class $[g]$ belongs to $\dot{H}^{\alpha}_p(\re)$
if and only if there exists some $G \in [g]$ such that
\begin{equation}\label{5}
\int_{-\infty}^\infty \Big[ \int_0^1   t^{-2\alpha} \Big( \frac 1t \int_{|h|<t} |G(x+h)-G(x)|\, dh\Big)^2 \frac{dt}{t}\Big]^{p/2} dx  <\infty\, .
\end{equation}
By Step 4 this implies that there is a representative $G \in [\ce_\alpha \tilde{f}]$ satisfying \eqref{5}.
Therefore we find a polynomial $P$ such that $G + P = \ce_\alpha \tilde{f}$.
For a function $F$ we put
\[
T(F) := \int_{0}^1 \Big[ \int_0^1   t^{-2\alpha} \Big( \frac 1t \int_{\{h : \: |h|<t, \: x+h \in (0,1)\} } 
|F(x+h)-F(x)|\, dh\Big)^2 \frac{dt}{t}\Big]^{p/2} dx \, .
\]
It is easy to see that 
$$T(\ce_\alpha \tilde{f}) \le 2^p \left( T(G) + T(P) \right).$$
Due to \eqref{5} we have $T(G) < \infty$, and the estimate
$$|P(x+h)-P(x)| \le |h| \max \big\{ |P'(y)| \,:\, y\in [-1,2] \big\} \hspace{3ex}\text{for}\hspace{2ex} x\in [0,1], h\in [-1,1],$$ 
shows immediately that $T(P)< \infty$. Hence we get
\[
\bigg(\int_{0}^1 \Big[ \int_0^1   t^{-2\alpha} \Big( \frac 1t \int_{\{h : \: |h|<t, \: x+h \in (0,1)\} } 
|\ce_\alpha \tilde{f}(x+h)-\ce_\alpha \tilde{f}(x)|\, dh\Big)^2 \frac{dt}{t}\Big]^{p/2} dx \bigg)^{1/2} <\infty\, .
\]
Because of $\ce_\alpha \tilde{f} \in L_p (0,1)$ this means  that 
 $\ce_\alpha \tilde{f}  \in H^\alpha_p (0,1)$, see \eqref{neueq2}.
This proves the claim.
\end{proof}

We recall a well-known embedding between Bessel potential and Besov spaces, see, e.g.,
\cite[2.3.2]{Tri83}.

\begin{lemma}\label{einbettung2}
Let $1<p<\infty$ and $\alpha \in (0,1)$.
Then we have  
$H^\alpha_p ([0,1]) \hookrightarrow B^\alpha_{p,\max \{p,2\}}([0,1])$.
\end{lemma}



We now show that for $1 \le p,q \le \infty$ the Besov space $B^\alpha_{p,q}$ is
continuously embedded in the Haar wavelet space $\mathcal{H}_{\wav,\alpha,p,q}$.
Recall that we constructed our wavelet spaces with respect to a fixed integer base $b\ge 2$.

\begin{theorem}\label{BesovHaar}
Let $\alpha\in (0,1)$, and let $1\le p,q \le \infty$. 
We have
\begin{equation*}
\|f\|_{\wav,\alpha,p,q} \le C_p \,
\|f\|^{(b)}_{B^\alpha_{p,q}}
\hspace{2ex}\text{for all $f\in B^\alpha_{p,q} ([0,1])$,}
\end{equation*}
where $C_p= \max\{1, 2^{1/p} (1-1/b) b^\alpha\}$.
In particular, we have the continuous embedding
\begin{equation*}
B_{p,q}^\alpha ([0,1]) \hookrightarrow \HH_{\wav,\alpha,p,q}.
\end{equation*}
\end{theorem}

\begin{proof}
Let $j\in\NN$. Then we obtain for $1 \le p < \infty$ with the help of H\"older's inequality
\begin{equation*}
 \begin{split}
  |\langle f, \psi^j_{i,k}\rangle | 
 &=
 \left| b^{\frac{j}{2}-1} \left( b
\int_{b^{-j}(bk+i)}^{b^{-j}(bk+i+1)} f(x) \,{\rm d}x 
- \int_{b^{1-j}k}^{b^{1-j}(k+1)} f(x) \,{\rm d}x 
 \right) \right|\\ 
&=
b^{\frac{j}{2}-1} \left| \sum_{\nu\in \nabla_j\setminus\{i\}}
\int_{b^{-j}(bk+\nu)}^{b^{-j}(bk+\nu+1)} \big(f(x+(i-\nu)b^{-j}) -f(x)\big)
\,{\rm d}x \right|\\
&\le b^{j(1/2 - 1/p')} b^{-1} \sum_{\nu\in \nabla_j\setminus\{i\}}
\left( \int_{b^{-j}(bk+\nu)}^{b^{-j}(bk+\nu+1)} |f(x+(i-\nu)b^{-j}) -f(x)|^p
\,{\rm d}x \right)^{1/p}.
 \end{split}
\end{equation*}
Note that $1/2-1/p' = 1/p-1/2$. Thus, using H\"older's inequality for sums, we have for $1\le p <\infty$ that
\begin{equation*}
 \begin{split}
  &\sum_{k\in\Delta_{j-1}} \sum_{i\in\nabla_j}
|\langle f, \psi^j_{i.k}\rangle |^p \\
&\le b^{p j (1/p - 1/2)} b^{-p} (b-1)^{p-1}  \sum_{k\in\Delta_{j-1}} \sum_{i\in\nabla_j}
\sum_{\nu\in \nabla_j\setminus\{i\}} \int_{b^{-j}(bk+\nu)}^{b^{-j}(bk+\nu+1)} |f(x+(i-\nu)b^{-j}) -f(x)|^p \,{\rm d}x \\ 
& \le 2 b^{p j(1/p-1/2)} b^{-p} (b-1)^{p-1}
\sum_{k\in\Delta_{j-1}} \sum_{\tau=1}^{b-1} \sum_{\nu=0}^{\tau-1}
\int_{b^{-j}(bk)}^{b^{-j}(bk+1)} |f(x+ \tau b^{-j}) -f(x + \nu b^{-j})|^p
\,{\rm d}x \\ 
&\le 2 b^{p j (1/p - 1/2)} b^{-p} (b-1)^{p-1} \sum_{\kappa=1}^{b-1}
\int_{0}^{1-\kappa b^{-j}} |f(x+ \kappa b^{-j}) -f(x)|^p \,{\rm d}x \\
&\le 2 b^{p j (1/p - 1/2)} b^{-p} (b-1)^{p} \sup_{0\le h \le (b-1)b^{-j}}
\int_0^{1-h} |f(x+h) -f(x)|^p \, {\rm d}x \\
&\le \left( 2^{1/p} b^{j (1/p - 1/2)} (1-1/b)\, \omega_p(f,(b-1)b^{-j}) \right)^p.
 \end{split}
\end{equation*}
Thus we have for $1\le q <\infty$
\begin{equation*}
 \begin{split}
\|f\|^q_{\wav, \alpha,p,q} &= \sum_{j=0}^{\infty} b^{q (\alpha-1/p+1/2) j} \left(
\sum_{k\in\Delta_{j-1}} \sum_{i\in\nabla_j} |\langle f, \psi_{i,k}^{j} \rangle
|^p \right)^{q/p} \\ &\le \left| \int_0^1 f(x)\,{\rm d}x \right|^{q} +
2^{q/p} (1-1/b)^q \sum_{j=1}^\infty b^{q (\alpha-1/p+1/2) j} b^{qj(1/p-1/2)} \omega_p(f,(b-1)b^{-j})^q\\ &\le \|f\|_{L_p}^q + 2^{q/p} (1-1/b)^q b^{q \alpha}
\sum_{j=0}^\infty b^{q\alpha j} \omega_p(f,b^{-j})^q \\
&\le \max \left\{1, 2^{1/p}(1-1/b) b^{\alpha} \right\}^q
\left( \|f\|_{L_p}^q + \sum_{j=0}^{\infty} b^{q \alpha j} \omega_p(f,b^{-j})^q \right).
\end{split}
\end{equation*}
Applying Jensen's inequality, we obtain
\begin{equation*}
\|f\|_{\wav, \alpha,p,q}  \le \max \left\{1, 2^{1/p}(1-1/b) b^{\alpha} \right\} \|f\|_{B^\alpha_{p,q}}^{(b)}.
\end{equation*}

For $q=\infty$ we have
\begin{equation*}
 \begin{split}
\|f\|_{\wav, \alpha,p,q} &= \sup_{j \in \mathbb{N}_0} b^{(\alpha-1/p+1/2) j} \left(
\sum_{k\in\Delta_{j-1}} \sum_{i\in\nabla_j} |\langle f, \psi_{i,k}^{j} \rangle
|^p \right)^{1/p} \\ &\le \max \left\{ \left| \int_0^1 f(x)\,{\rm d}x \right|,
2^{1/p} (1-1/b) \sup_{j \in \mathbb{N}} b^{(\alpha-1/p+1/2) j} b^{j(1/p-1/2)} \omega_p(f,(b-1)b^{-j}) \right\}\\ &\le \|f\|_{L_p} + 2^{1/p} (1-1/b) b^{\alpha}
\sup_{j \in \mathbb{N}_0} b^{\alpha j} \omega_p(f,b^{-j}) \\
& \le \max \left\{1, 2^{1/p}(1-1/b) b^{\alpha} \right\}
\left( \|f\|_{L_p} + \sup_{j \in \mathbb{N}_0} b^{\alpha j} \omega_p(f,b^{-j}) \right) \\ & \le \max \left\{1, 2^{1/p}(1-1/b) b^{\alpha} \right\} \|f\|_{B^\alpha_{p,\infty}}^{(b)}.
\end{split}
\end{equation*}

Furthermore, we have for $p=\infty$ and $j\in\NN$
\begin{equation*}
 \begin{split}
  |\langle f, \psi^j_{i.k}\rangle | &=
b^{\frac{j}{2}-1} \left| \sum_{\nu\in \nabla_j\setminus\{i\}}
\int_{b^{-j}(bk+\nu)}^{b^{-j}(bk+\nu+1)} \big(f(x+(i-\nu)b^{-j}) -f(x)\big)
\,{\rm d}x \right|\\
&\le b^{-j/2-1} \sum_{\nu\in \nabla_j\setminus\{i\}}
\esssup_{x \in [b^{-j} (bk+\nu), b^{-j} (bk+\nu+1)]} |f(x+(i-\nu)b^{-j}) -f(x)|
\\ & \le b^{-j/2} (1-1/b)\max_{\nu \in \nabla_j} \esssup_{x \in [b^{-j} (bk+\nu), b^{-j} (bk+\nu+1)]} |f(x+(i-\nu)b^{-j}) -f(x)|.
 \end{split}
\end{equation*}
Thus for $j \in \mathbb{N}$ we have
\begin{equation*}
 \begin{split}
  & \max_{k \in \Delta_{j-1}} \max_{i\in\nabla_j}
|\langle f, \psi^j_{i.k}\rangle | \\
&\le b^{- j/2} (1-1/b) \max_{k\in\Delta_{j-1}} \max_{i\in\nabla_j}
\max_{\nu\in \nabla_j} \esssup_{x \in [b^{-j}(bk+\nu), b^{-j}(bk+\nu+1)]} |f(x+(i-\nu)b^{-j}) -f(x)| \\
& \le b^{-j/2}(1-1/b) \sup_{0 < h \le (b-1) b^{-j} } \esssup_{x \in [0,1-h]} |f(x+h)-f(x)| \\
& \le b^{-j/2}(1-1/b) \omega_\infty(f,(b-1)b^{-j}).
 \end{split}
\end{equation*}

Hence, for $p = \infty$ and $1\le q <\infty$ we have
\begin{equation*}
\begin{split}
\|f\|_{\wav, \alpha,\infty,q}^q &= \sum_{j=0}^\infty b^{q(\alpha+1/2)j} \left(\max_{k \in \Delta_{j-1}} \max_{i \in \nabla_j} |\langle f, \psi^j_{i,k}\rangle | \right)^{q} \\
& \le \left|\int_0^1 f(x) \,\mathrm{d} x \right|^q +  (1-1/b)^q
\sum_{j=1}^\infty b^{q \alpha j} \omega_\infty(f,(b-1)b^{-j})^q \\
& \le \|f\|_{L_\infty}^q + b^{q \alpha} (1-1/b)^q
\sum_{j= 0}^\infty b^{q \alpha j} \omega_\infty(f, b^{-j})^q \\
& \le \max\{1, b^\alpha(1-1/b)\}^q \left(\|f\|_{L_\infty}^q + \sum_{j=0}^\infty b^{q \alpha j}
\omega_\infty(f, b^{-j})^q \right).
\end{split}
\end{equation*}
Using Jensen's inequality we obtain
\begin{equation*}
\|f\|_{\wav, \alpha, \infty, q} \le \max\{1, b^\alpha(1-1/b)\}
\|f\|^{(b)}_{B^\alpha_{\infty,q}}.
\end{equation*}
The same result can also be shown for $p,q=\infty$ with some minor modifications.
\end{proof}

\begin{remark}
 \rm
 Theorem \ref{BesovHaar} has a counterpart for the classical Haar system, we refer to 
 \cite[Prop.~2.37]{Tri10}. There it is formulated for spaces defined on $\real$.
Further let us mention that for the  classical Haar system there is a complete understanding of the relationships between
 sequence spaces (of Haar coefficients) and Besov spaces including all limiting situations, see \cite{GSU21}.
 \end{remark}

We summarize the for us most important embedding results of this section in the following corollary.

\begin{corollary}\label{DirEmb}
Let $\alpha\in (0, 1)$ and $\alpha^{-1} < p  < \infty$. We have the continuous embeddings
\begin{equation}\label{univariate_embedding_series}
\mathcal{H}_{\alpha, p} \hookrightarrow H^\alpha_{p}([0,1])
\hookrightarrow 
\begin{cases}
B^\alpha_{p, 2}([0,1]) \hookrightarrow \HH_{\wav,\alpha,p, 2}, & \text{if $p < 2$,}\\
 B^\alpha_{2, 2}([0,1]) \hookrightarrow \HH_{\wav,\alpha, 2, 2}, & \text{if $p \ge 2$.}
\end{cases}
\end{equation}
\end{corollary}

\begin{proof}
We simply combine Theorem~\ref{Winfried}, Lemma~\ref{einbettung2}, and Theorem~\ref{BesovHaar}, 
and use in the case $p > 2$ additionally \eqref{monotonicity_bessel_potential_spaces} with $p_1:=2$ and $p_2:= p$.
\end{proof}


As the next proposition shows, the embeddings in \eqref{univariate_embedding_series} 
cannot simply be extended to the case $\alpha =1$.

\begin{proposition}\label{prop1}
Let $1\le p \le \infty$ and $1\le q < \infty$. The function $f:[0,1]\to \RR, x\mapsto x$ is in $\HH_{1,p}$, but not in
$\HH_{\wav,\alpha,p,q}$ for any
$\alpha \ge 1$. Thus the space $\HH_{1,p}$
is not included in $\HH_{\wav,\alpha,p,q}$ for any $\alpha\ge 1$.
\end{proposition}


\begin{proof}
For all $1\le p \le \infty$ the function $f(x)=x$ is in $\HH_{1,p}$,
since we have for $\tilde{f} :=\psi^0_{0,0} =1_{[0,1)} \in L_p([0,1]).$ that
\begin{equation*}
f(x) = \int^x_0 \tilde{f}(t) \,{\rm d}t
= \frac{1}{\Gamma(1)} \int^1_0 \tilde{f}(t)(x-t)^{0}_+ \,{\rm d}t.
\end{equation*}
Furthermore, we have $\langle f, \psi^0_{0,0} \rangle = 1/2$
and for $j\ge 1$, $i\in\nabla_j$, and $k\in\Delta_{j-1}$
$$\langle f, \psi^j_{i,k} \rangle =
\frac{1}{2} (2i+1-b)b^{-\frac{3}{2}j}.$$
Hence we obtain, with the obvious modifications for $p=\infty$,
\begin{equation*}
\begin{split}
\|f\|^{q}_{\wav,\alpha, p,q} =& | \langle f, \psi^0_{0,0} \rangle|^{q}
+ \sum_{j=1}^\infty b^{q(\alpha-1/p+1/2) j} \left( \sum^{b-1}_{i=0}
\sum_{k\in\Delta_{j-1}} |\langle f, \psi^j_{i,k} \rangle|^{p} \right)^{q/p}\\
=& 2^{-q}  \left( 1 + \left( \sum^{b-1}_{i=0} |2i+1-b|^{p} \right)^{q/p}
b^{-q/p} \sum_{j=1}^\infty b^{q(\alpha - 1)j} \right).
\end{split}
\end{equation*}
Thus $\|f\|_{\wav,\alpha,p,q}$ is finite if and only if $\alpha <1$.
\end{proof}


\subsection{The Multivariate Case} 
\label{Subsec:Besov_Multi}


\subsubsection{Preliminaries}

We recall some results from \cite{SU1,SU2}
on tensor products of Besov and Bessel potential spaces.
Therefore we need to introduce function spaces of dominating mixed smoothness.
Most transparent are tensor products of Bessel potential spaces.

\paragraph{Bessel Potential Spaces.}

Let $1<p<\infty$ and $\alpha \in \real$.
The \emph{Bessel potential space} or \emph{(fractional) Sobolev space with dominating mixed smoothness}
$S^\alpha_pH(\real^s)$ is  the collection of all tempered distributions $f
\in \mathcal{S}'(\real^s)$ such that
\[
\| \, f \, \|_{S^{\alpha}_p H(\real^s)}: = \Big\|\, \cfi
\Big[\prod\limits_{i=1}^s (1+|\xi_i|^2)^{\alpha/2}\, \cf f (\xi)
\Big](\, \cdot \, )\Big\|_{L_p(\real^s)}
\]
is finite.  Similarly as in the univariate case, we have $S^\alpha_pH(\real^s) \subseteq L_p(\real^s)$ for $\alpha \ge 0$.

If we take $\alpha =1$, then we obtain $S^{1}_p H(\real^s) = S^{1}_p W(\real^s)$ in the sense of equivalent norms,
where $S^{1}_p W(\real^s)$ denotes the first order Sobolev space of dominating mixed smoothness.
This space is the collection of all $f\in L_p (\real^s)$ such that
\[
 \|f\|_{S^{1}_p W(\real^s)} := \|f\|_{L_p(\real^s)} + \sum_{{{\bf m}\in \N_0^s} \atop \|{\bf m}\|_\infty \le 1}
 \|D^{\bf m} f\|_{L_p (\real^s)}
<\infty\, .
 \]
Here $D^{\bf m} f$ is understood as the corresponding distributional derivate of $f$.
In many problems the mixed derivative $D^{(1, 1, \ldots \, 1)} f$ is the dominating one, which is  the reason for the name.


Later on we shall need also a wavelet characterization of theses spaces.
It will be convenient to work with compactly supported tensor product wavelets.

Let $u >0$.
Let $\phi \in C^u (\real)$ be a univariate compactly supported scaling function
and denote by   $\vartheta \in C^u (\real)$ an associated univariate compactly supported 
wavelet satisfying a moment condition 
\[
 \int_{\real} \vartheta  (t)\, t^m \, dt =0 \quad \mbox{for all}\quad m \in \NN_0 
\quad \mbox{with} \quad m< u\, .
\]
It follows that
$\{ 2^{j/2}\, \vartheta_{j,m} \,|\, \quad j \in \NN_0, \: m \in \ZZ\}$, where
\[
\vartheta_{j,m} (t):= \left\{ \begin{array}{lll}
	\phi (t-m) & \qquad & \mbox{if}\quad j=0, \: m \in \zz\, , 
	\\
	\sqrt{1/2} \, \vartheta (2^{j-1}t-m)& \qquad & \mbox{if}\quad j\in \N\, , \: m \in \zz\, , 
\end{array} \right.
\]
is an orthonormal basis in $L_2 (\re)$, see \cite{Woj}. Let $s\ge 2$.
Furthermore, let ${\nu} =(\nu_1, \ldots ,\nu_s) \in \mathbb{N}_0^s$ and ${m}= (m_1, \ldots , m_s)\in \mathbb{Z}^s$.
We put
\[
\Theta_{{\nu}, m} (x) := \prod_{\ell=1}^{s} 2^{\nu_\ell/2} \vartheta_{\nu_\ell, m_\ell} (x_\ell)\, . 
\]
Then $\{ \Theta_{{\nu}, {m}} \,|\, {\nu} \in \NN_0^s, \, {m} \in \ZZ^s\, 
\}$
is a tensor product wavelet basis of $L_2 (\real^s)$.
We 
put
\[
Q_{{\nu},{m}} := \Big\{x\in \real^s \, \Big|\,  \quad 2^{-\nu_\ell} \, m_\ell \le  x_\ell < 2^{-\nu_\ell}\, (m_\ell+1)\, , \: \ell = 1, \, \ldots \, , s\Big\} \, .
\]
As before, we denote by $1_{Q_{{\nu},{m}}}(x)$ the characteristic function of $Q_{{\nu},{m}}$ and use the notation $|j| := j_1 + \ldots \, + j_s $ for $j \in \NN_0^s$.

If $1<p < \infty$, $ \alpha \in \mathbb{R}$ and
	$\lambda:=\lbrace \lambda_{{\nu},{m}}\in\mathbb{C} \,|\, {\nu}\in \mathbb{N}_0^s,\ {m}\in \mathbb{Z}^s \rbrace$,
	then we define
\begin{equation}\label{def:s_alpha_p_h}
	s_{p}^\alpha h := \Big\lbrace\lambda \,\Big|\, \| \lambda \|_{s_{p}^\alpha h} =
	\Big\|\Big(\sum_{{\nu}\in \mathbb{N}_0^s}\sum_{{m}\in \mathbb{Z}^s}|2^{|{\nu}| (\alpha + 1/2)} 	\lambda_{{\nu},{m}} 1_{Q_{\nu},{m}}(.) 
	|^2\Big)^{\frac{1}{2}}\Big\|_{ L_p(\mathbb{R}^s)} <\infty \Big\rbrace. 
\end{equation}
For a tempered distribution $f\in \mathcal{S}'(\real^s)$ and a rapidly decreasing function $\varphi \in \mathcal{S}(\real^s)$ we denote their dual bracket by $\langle f, \varphi \rangle := f(\varphi)$. Clearly, if $f\in L_p(\real^s)$, then the dual bracket is nothing but the integral of the product $f \cdot \varphi$ over $\real^s$.

\begin{proposition}\label{waveBesovS}
Let $1 <p< \infty$ and $\alpha \in \real$. Suppose $ u> |\alpha|$.
Then the mapping $\mathcal{W}:~ f \mapsto \, (\langle f, \Theta_{\nu,m} \rangle)_{\nu,m}$
is a continuous linear isomorphism of $S^\alpha_p H(\real^s)$ onto $s^\alpha_p h(\real^s)$
\end{proposition}

We refer to Vybiral \cite{Vy} and Triebel \cite[Thm.~1.54]{Tri10}.


Now we turn to the tensor product interpretation.
We shall use the convention $S^{\alpha}_{p}H(\real):= H^{\alpha}_{p}(\real)$.
For the next result see, e.g., \cite{SU1, SU2}. 

\begin{proposition}\label{tensorb} 
Let $s\ge 1$, $\alpha \in \real$ and $1<p<\infty$. Then we have
 \begin{equation*}\label{eq-5}
    \nonumber H^\alpha_p(\real) \otimes_{\delta_p} S^\alpha_pH(\real^s) = S^\alpha_pH(\real^s) \otimes_{\delta_p} H^\alpha_p(\real) = \nonumber S^\alpha_pH(\real^{s+1})
 \end{equation*}
in the sense of equivalent norms.
\end{proposition}

\begin{remark}\label{Lebesgue}
The above identities contain for $\alpha =0$  the well-know relations
\begin{equation}\label{ws-02}
 L_p (\real)\otimes_{\delta_p} L_p (\real^s) 
= L_p (\real^s)\otimes_{\delta_p} L_p (\real) 
= L_p (\real^{s+1})\, ,
\end{equation}
see, e.g., \cite{LiCh}. But \eqref{ws-02} is true also for $p=1$.
\end{remark}

Since we want to deal with spaces on the unit cube $[0,1]^{s}$, 
we shall apply restrictions  to switch from $\real^{s}$  to the cube.  
Let $\alpha \ge 0$. 
The space
$S^{\alpha}_{p}H([0,1]^s)$ is the space of all $f\in
   L_p([0,1]^s)$ such that there exists a $g\in
   S^{\alpha}_{p}H(\real^s)$ satisfying $f = g|_{[0,1]^s}$. It is endowed  with the quotient norm
   $$
 \|\,f\,\|_{S^{\alpha}_{p}H([0,1]^s)} := \, \inf\{\|\, g \,\|_{S^{\alpha}_{p}H(\real^s)} \,|\, g\in S^{\alpha}_{p}H(\real^s)\,,\, g_{|_{[0,1]^s}} =f\}\,.
   $$

\begin{lemma}\label{Win_Lemma_13}
 Let $s\ge 1$,  $\alpha \ge 0$, and $1 < p<\infty$. Then the following statements hold. 
 \\
{\rm (i)} 
$ S^{\alpha}_{p}H([0,1]^s)$ is a Banach space.
\\
{\rm (ii)} The scale $ S^{\alpha}_{p}H([0,1]^s)$
is monotone in $p$ and $\alpha$, i.e., for any $\varepsilon >0$
we have 
\[
 S^{\alpha+\varepsilon}_{p}H([0,1]^s) \hookrightarrow  S^{\alpha}_{p}H([0,1]^s)
\]
 and 
 \[
 S^{\alpha}_{p+\varepsilon} H([0,1]^s) \hookrightarrow  S^{\alpha}_{p}H([0,1]^s)\, .
\]
 \end{lemma}

\begin{proof}
 Part (i) can be derived from the fact that 
 $ S^{\alpha}_{p}H(\real^s) $ is a Banach space.
 
 
Now we turn to 
 part (ii).
 We claim that there exists a constant $C>0$ such that each
 $f \in S^{\alpha}_{p}H([0,1]^s)$ has an extension $g$ to $\real^s$
 satisfying $\supp g \subset [-1,2]^s$ and
 \[
  \| g  \|_{S^{\alpha}_{p}H(\real^s)} \le C\, \| f  \|_{S^{\alpha}_{p}H([0,1]^s)}\,,
  \]
  and therefore
 \begin{equation}\label{eq:asymp_equiv}
  \| g \|_{S^{\alpha}_{p}H(\real^s)} \asymp \| f \|_{S^{\alpha}_{p}H([0,1]^s)} \, .
  \end{equation}
Indeed, let $\psi \in C_0^\infty (\real^s)$ such that $\psi (x)=1$ on 
$[0,1]^s$ and $\supp \psi \subset [-1,2]^s$. 
Then $\psi$ is a pointwise multiplier on  $S^{\alpha}_{p}H(\real^s)$, i.e.,
the mapping $h \mapsto  h \cdot \psi$ is a linear endomorphism on $S^{\alpha}_{p}H(\real^s)$ with finite operator norm $C_\psi$, see, e.g.,  \cite{NS}.
Let $g$ denote an arbitrary extension of $f$.
Then
\[
 \| f \|_{S^{\alpha}_{p}H([0,1]^s)}  \le  
 \| g\, \cdot \, \psi  \|_{S^{\alpha}_{p}H(\real^s)} 
 \le C_{\psi}\,
 \| g  \|_{S^{\alpha}_{p}H(\real^s)}\,.
 \]
 By definition of the norm on $S^{\alpha}_{p}H(\real^s)$ we may choose $g$ in such a way that  we have
 \[
  \| g \|_{S^{\alpha}_{p}H(\real^s)}  \le  2\, 
 \| f  \|_{S^{\alpha}_{p}H([0,1]^s)} \, , 
 \]
which concludes the proof of the claim. 

Now we  apply the wavelet characterization of the classes  $S^{\alpha}_{p}H(\real^s)$, see \eqref{def:s_alpha_p_h} and Proposition~\ref{waveBesovS}.  
Let $g$ be the extension of $f$ with support in $[-1,2]^s$ from above.
Observe that the function 
\[
\Big(\sum_{{\nu}\in \mathbb{N}_0^s}\sum_{{m}\in \mathbb{Z}^s}|2^{|{\nu}| (\alpha + 1/2)} \,	\langle g, \Theta_{\nu,m} \rangle 1_{Q_{\nu},{m}}(.) 
|^2\Big)^{\frac{1}{2}}
\]
increases monotonically with respect to the parameter $\alpha$, which establishes the first claimed embedding. Furthermore, it 
has its support in a ball with centre in the origin and radius $L$, where $L$ only depends on the generators $\phi$, $\vartheta$ of the wavelet system. 
Hence, we can  trace the second claimed embedding back to the well-known relation 
$L_{p+\varepsilon}(\Omega) \hookrightarrow L_{p}(\Omega)$
for bounded measurable domains $\Omega \subset \real^s$.
\end{proof}

\paragraph{Besov Spaces.}

Next we recall a corresponding result for Besov spaces.
Here we would like to avoid to give  the definition of the spaces $S^\alpha_{p,q}B(\real^s)$.

For us it will be sufficient to define these spaces by their wavelet characterization.
Therefore we will use the same notation as in case of the spaces $S^\alpha_{p}H(\real^s)$.

If $1\le p,q\leq \infty$, $\alpha \in \mathbb{R}$ and
	$\lambda:=\lbrace \lambda_{{\nu},{m}}\in\mathbb{C}:~{\nu}\in \mathbb{N}_0^s,  \, {m}\in \mathbb{Z}^s \rbrace$,
	then we define
$$
s_{p,q}^\alpha b := \Big\lbrace\lambda: \| \lambda \|_{s_{p,q}^\alpha b} =
\Big(\sum_{{\nu}\in \mathbb{N}_0^s}2^{|{\nu}| (\alpha -\frac{1}{p}+\frac 12 )q}\big(\sum_{{m}\in \mathbb{Z}^s}|\lambda_{{\nu},{m}} 
|^p\big)^{\frac{q}{p}}\Big)^{\frac{1}{q}}<\infty \Big\rbrace$$
with the usual modification for $p$ or/and q equal to $\infty$.

By using this notation Besov spaces of dominating mixed smoothness can now be introduced as follows.

\begin{proposition}\label{wavelet}
	Let $1\le  p,q \le\infty$ and $\alpha \in \real$. Suppose $u > |\alpha|$. 	
Then the  mapping 
\begin{equation}\label{wave}
		{\mathcal W}: \quad f \mapsto (\langle f, \Theta_{{\nu},{m}} \rangle)_{{\nu} \in \N_0^d\, , \, {m} \in \ZZ^s} 
\end{equation}
is a continuous linear isomorphism of $S^\alpha_{p,q}B(\real^s)$ onto $s^\alpha_{p,q}b$.
\end{proposition}

A proof of this characterization can be found in 
Vybiral \cite{Vy}, but  see also \cite[Thm.~1.54]{Tri10}.
We may take it as a definition here.
The independence of this definition from the generators $(\phi,\vartheta)$
of the wavelet system 
follows immediately from the equivalence to Fourier analytic characterizations, see again  \cite{Vy} and \cite{Tri10}.


There exist  further characterizations of these classes, e.g., by differences, we refer to 
\cite{ST}, \cite{Tri18} and \cite{Ul10}.

Also in case of Besov spaces we shall use the convention $S^{\alpha}_{p, q}B(\real):= B^{\alpha}_{p, q}(\real)$.

We proceed  as in  case of the Bessel potential spaces. 
The space
$S^{r}_{p,q}B([0,1]^s)$ is the space of all $f\in
   L_p([0,1]^s)$ such that there exists a $g\in
   S^{r}_{p,q}B(\real^s)$ satisfying $f = g|_{[0,1]^s}$. It is endowed  with the quotient norm
   $$
 \|\,f\,\|_{S^{r}_{p,q}B([0,1]^s)} := \inf\, \{\|\, g \,\|_{S^{r}_{p,q}B(\real^s)} \,|\, g_{|_{[0,1]^s}} =f\}\,.
   $$
 
Now we can formulate the characterization we will use below.

\begin{proposition}\label{tensor4} 
Let $s\ge 1$. 
{\rm (i)} Let   $\alpha \in [0,1]$ and $1 < p<\infty$.
Then
 \begin{equation}\label{eq-5b}
H^\alpha_p([0,1]) \otimes_{\delta_p} S^\alpha_pH([0,1]^s) = S^\alpha_pH([0,1]^s) \otimes_{\delta_p} H^\alpha_p([0,1])=  S^\alpha_pH([0,1]^{s+1})
 \end{equation}
 holds in the sense of equivalent norms.
 \\
{\rm (ii)} Let   $\alpha \in (0,1)$ and $1 \le  p<\infty$.
Then
 \begin{equation}\label{eq-5c}
B^\alpha_{p,p}([0,1]) \otimes_{\delta_p} S^\alpha_{p,p} B([0,1]^s) = S^\alpha_{p,p}B([0,1]^s) \otimes_{\delta_p} B^\alpha_{p,p}([0,1])=  S^\alpha_{p,p}B([0,1]^{s+1})
 \end{equation}
 holds in the sense of equivalent norms.
\end{proposition}

A proof of Proposition \ref{tensor4} can be found in the appendix in \cite{SU2}.

For both types of spaces there exist intrinsic characterizations as well. We omit the details.

\subsubsection{Embedding Results} 


In a first step, we want to embed our function spaces of fractional smoothness into corresponding
Bessel potential spaces.
Due to the fact that the $p$-nuclear tensor norm $\delta_p$ is a uniform cross norm and due to
Theorem~\ref{tensor3} and~\ref{Winfried}, and 
to Proposition~\ref{tensor4}(i),
we obtain the following corollary.

\begin{corollary}\label{Winfried2}
Let $s\in \N$,  
$1 <  p < \infty$, and $\alpha\in (1/p,1)$.
Then we have the continuous embedding
\begin{equation*}
\HH_{\alpha,s,p} \hookrightarrow S^\alpha_{p}H([0,1]^s).
\end{equation*}
\end{corollary}

The proof Corollary~\ref{Winfried2} can be found in Section~\ref{APP:TPSO} of the appendix on tensor products of 
Banach spaces.

In a second step we want to embed Bessel potential spaces into suitable Besov spaces over
the $s$-dimensional unit cube and vice versa.
To this purpose we make use of the known embeddings for $1<p <\infty$ and $\alpha\in (0,1)$
\be\label{elementaryembedding}
S^\alpha_{p,\min(p,2)}B(\real^s) \hookrightarrow S^\alpha_{p}H(\real^s)\hookrightarrow 
S^\alpha_{p,\max(p,2)}B(\real^s)\, , 
\ee
see \cite[Prop.~2.2.3/2(iii)]{ST}. 
This result and the definition of the norms on $S^\alpha_{p}H([0,1]^s)$ and $S^\alpha_{p,q}B([0,1]^s)$
via restriction immediately yield the next corollary.

\begin{corollary}\label{Cor:Neu_Equation_60+}
Let $s\in \N$,  
$1 < p < \infty$, and $\alpha\in (0,1)$. 
Then we have the continuous embedding
\begin{equation*}
S^\alpha_{p,\min(p,2)}B([0,1]^s) \hookrightarrow S^\alpha_{p}H([0,1]^s)\hookrightarrow 
S^\alpha_{p,\max(p,2)}B([0,1]^s)\,.
\end{equation*}
\end{corollary}






In a third step, we want to establish embeddings between Besov spaces of dominated mixed
smoothness and corresponding Haar wavelet spaces.

\begin{remark}\label{Rem:Emb_Bes_Wav_pp}
We may argue parallel to Corollary \ref{Winfried2} to obtain a first embedding result between Besov and Haar wavelet spaces: Since $\delta_p$ is a uniform cross norm, 
Proposition~\ref{tensor4}, 
Theorem~\ref{BesovHaar}, 
and 
Theorem~\ref{tensor17}
imply the following result.\\
Let $s\in \N$, $\alpha\in (0,1)$, and $1 \le p< \infty$.
Then we have the continuous embedding
\begin{equation*}
S^{\alpha}_{p,p} B([0,1]^{s}) \hookrightarrow  \HH_{\wav,\alpha,s,p,p} \, .
\end{equation*}
It follows easily from Proposition~\ref{prop1} that
this 
embedding result 
cannot be extended to
the case where $\alpha =1$. 
\end{remark}

The next theorem generalizes the result in Remark~\ref{Rem:Emb_Bes_Wav_pp} substantially.



\begin{theorem}\label{tensor18a}
 Let $1 \le p, q \le \infty$, $\alpha \in (0,1)$ 
 and $s\in \N$.
 Then we have the continuous embedding
 \be\label{tensor19a}
 S^\alpha_{p,q} B ([0,1]^s) \hookrightarrow 
 \HH_{\wav,\alpha,s,p,q}\,.
 \ee
\end{theorem}

\begin{proof}
	{\em Step 1.} Let us fix the base $b=2$.
We will apply a general estimate of local means in the framework of Besov spaces
of dominating mixed smoothness, due to Vybiral \cite{Vy} and Triebel \cite[Thm.~1.52]{Tri10}.
The formulation in \cite[Thm.~1.52]{Tri10} is slightly more convenient for us, but restricted there to $s=2$. We will use it here  for general $s$, $A=0$ and $B=1$.
We extend the system 
$\psi^{j}_{i,k}$, $j \in \N_0$,  $k \in \Delta_{j - 1}$, $i \in \nabla_{j}$, defined in Section~\ref{HAAR-1}, 
from values $k \in \Delta_{j - 1}$ to arbitrary values $k \in \zz$.
Let $\bsj \in \N_0^s$, $\bsk \in \zz^s$, $\bsi \in \{0,1\}^s$, and $x\in \real^{s}$. 
Then we put again
\beqq
\Psi^{\bsj}_{\bsi,\bsk} (x)& := &  \prod_{\ell=1}^s
\psi^{j_\ell}_{i_\ell,k_\ell} (x_\ell)
\\
& = &  \prod_{\ell=1}^s 
 2^{j_\ell/2-1} \left( 2\, \cdot 1_{[2^{-j_\ell}(2k_\ell+i_{\ell}), 2^{-j_\ell}(2k_\ell+i_\ell+1))}(x_{\ell})
- 1_{[2^{1-j_\ell}k_\ell, 2^{1-j_\ell}(k_\ell+1))}(x_{\ell}) \right)\, .
\eeqq
Let $0 < \alpha <1$ and $1\le p,q \le \infty$.
Then it follows from \cite[Thm.~1.52]{Tri10} that  there exists a positive constant $C$ such that 
\be\label{localm}
\Big(
\sum_{\bsj \in \N^{s}_{0}} 2^{|\bsj|(\alpha -1/p)q} \Big(\sum_{k \in \zz^s} \Big|
\int_{\real^s} 2^{|\bsj|/2}\, \Psi^{\bsj}_{{\bf 0},\bsk} (x) \, f(x) \, dx
\Big|^p\Big)^{q/p} \Big)^{1/q} \le C\, \|f\|_{S^\alpha_{p,q}B(\real^s)}
\ee
holds for all functions $f\in S^\alpha_{p,q}B(\real^s)$,  
where ${\bf 0} = (0,\ldots,0)$.
\\
Recall that  $S^\alpha_{p,q}B([0,1]^s)$ is defined by restrictions.
Let $g\in S^\alpha_{p,q}B([0,1]^s)$ and denote by  $Eg$ any admissible extension of $g$ to $\re^s$
such that 
\[
\| Eg\|_{S^\alpha_{p,q}B(\real^s)}\le 2\, \| g\|_{S^\alpha_{p,q}B([0,1]^s)} \, .
\]
Then \eqref{localm} 
yields
\beq\label{n00}
\| g \|_{\HH_{\wav,\alpha,s,p,q}} & = &
\Big(
\sum_{\bsj \in \N^{s}_0} 2^{|\bsj|(\alpha -1/p+1/2)q} \Big(\sum_{\bsk \in \Delta_{j-1}} 
 \sum_{\bsi \in \nabla_{\bsj}}
\Big|\int_{[0,1]^s} \Psi^{\bsj}_{\bsi,\bsk} (x) \, g(x) \, dx
\Big|^p\Big)^{q/p} \Big)^{1/q}
\nonumber
\\
& = & 
2^{s/p}\, 
\Big(
\sum_{\bsj \in \N^{s}_0} 2^{|\bsj|(\alpha -1/p+1/2)q} \Big(\sum_{\bsk \in \Delta_{\bsj-1}} \Big|
\int_{\real^s}  \Psi^{\bsj}_{{\bf 0},\bsk} (x) \, Eg(x) \, dx
\Big|^p\Big)^{q/p} \Big)^{1/q}
\nonumber
\\
&\le &
2^{s/p}\,
\Big(
\sum_{\bsj \in \N^{s}_0} 2^{|\bsj|(\alpha -1/p+1/2)q} \Big(\sum_{\bsk \in \zz^s} \Big|
\int_{\real^s}  \Psi^{\bsj}_{{\bf 0},\bsk} (x) \, Eg(x) \, dx
\Big|^p\Big)^{q/p} \Big)^{1/q}
\nonumber
\\
&\le & 2^{s/p}\, C\, \|Eg\|_{S^\alpha_{p,q}B(\real^s)} \le 2^{1+ s/p} \, C \, \|g\|_{S^\alpha_{p,q}B([0,1]^s)}\, . 
\eeq
{\em Step 2.} Let $b>2$. Here we give a very rough sketch of the proof only.
Let $\rho \in C_0 ^\infty (\real)$ be a function such that $\rho (t)=1$ if $|t|\le 1$ and 
$\rho (t)=0$ if $|t|\ge 3/2$. Then we define
\[
\varphi_0 := \rho \, , \quad \varphi_1 (t) := \rho (t/b) - \rho (t) \, \quad 
\mbox{and}\quad \varphi_j (t) := \varphi_{1} (b^{-j+1} t)\, , \quad t \in \real\, , 
\]
and $j \in \NN$.
As a consequence we obtain a smooth decomposition of unity, i.e.,
\[
\sum_{j=0}^\infty \varphi_j (t) = 1 \qquad \mbox{for all}\quad t \in \real\, .
\]
Now we take tensor products.
Hence
\[
\sum_{j\in \NN_0^s} \varphi_{ j_1} (x_1) \, \cdot \, \ldots \, \cdot \, \varphi_{j_s} (x_s)= 1 \qquad \mbox{for all}\quad x \in \real^s\, .
\]
For brevity we put
\[
\Phi_j (x) :=  \varphi_{ j_1} (x_1) \, \cdot \, \ldots \, \cdot \, \varphi_{j_s} (x_s)\, ,\qquad 
 x \in \real^s\, .
\]
For those systems $(\Phi_j)_{j \in \NN_0}$ we can apply the theory of regular coverings, described in \cite[2.2.1]{Tri78}.
As a consequence we obtain that $S^\alpha_{p,q}B(\real^s)$ can be described as follows.
This space is the collection of all tempered distributions $f$ such that 
\[
\|\, f \, \|_{S^\alpha_{p,q}B(\real^s)}^*:= \Big(\sum_{j \in \NN_0^s} b^{|j| \alpha q}\, \|\mathcal{F}^{-1}[\Phi_j \mathcal{F}f ] \|_{L_p (\real^s)}^q\Big)^{1/q} <\infty\, .
\] 
Moreover, the norm $\|\, \cdot \, \|_{S^\alpha_{p,q}B(\real^s)}^*$ is equivalent to the norm $\|\, \cdot \, \|_{S^\alpha_{p,q}B(\real^s)}$. There is no reference for the dominating mixed case, but for the isotropic case we refer to Proposition 2.2.1 in \cite{Tri78}.
There is not much difference between these two cases with this respect.
Now one has to follow the proof of the local mean characterization in Vybiral \cite{Vy}
step by step, taking into account the modifications indicated by Triebel in \cite{Tri10} (proof of Theorem 1.52 and Remark 1.53). By working with the norm $\|\, \cdot \, \|_{S^\alpha_{p,q}B(\real^s)}^*$
all quantities in the needed generalization of \eqref{n00} depend on a common scaling factor $b>1$.
Having established this generalization, we obtain
$S^\alpha_{p,q}B([0,1]^s) \hookrightarrow \HH_{\wav,\alpha,s,p,q}$ for all $b >1$.
 \end{proof}

We summarize the main embedding results of this section in the following corollary.

\begin{corollary}\label{DirEmb_Multi}
Let $s\in\N$, $\alpha\in (0, 1)$, and $\alpha^{-1} < p  < \infty$. We have the continuous embeddings
\begin{equation}\label{multivariate_embedding_series}
\mathcal{H}_{\alpha, s, p} \hookrightarrow S^\alpha_{p}H([0,1]^s)
\hookrightarrow 
\begin{cases}
S^\alpha_{p, 2}B([0,1]^s) \hookrightarrow \HH_{\wav,\alpha, s, p, 2}, & \text{if $p < 2$,}\\
S^\alpha_{2, 2}(B[0,1]^s) \hookrightarrow \HH_{\wav,\alpha, s, 2, 2}, & \text{if $p \ge 2$.}
\end{cases}
\end{equation}
\end{corollary}

\begin{proof}
We combine
Corollary~\ref{Winfried2} and~\ref{Cor:Neu_Equation_60+}, and Theorem~\ref{tensor18a}, and in the case $p>2$ we additionally use Lemma~\ref{Win_Lemma_13}.
\end{proof}


\subsection{Quasi-Monte Carlo Integration and Fractional Discrepancy}
\label{Subsec:QMC_Frac_Disc}


We want to transfer our optimal results on QMC-integration for Haar wavelet spaces from 
 Theorem~\ref{Nets} with the help of the embedding results summarized in Corollary~\ref{DirEmb_Multi} 
 to Besov and Sobolev spaces of dominating mixed smoothness and to our spaces of fractional smoothness and to fractional discrepancies.
We start with Besov and Sobolev spaces of dominating mixed smoothness.

First we state 
upper error bounds for 
QMC rules based on $(t,m,s)$-nets for Besov and Sobolev spaces of dominated mixed smoothness
in Corollary~\ref{Cor:Nets_Besov_Spaces} and~\ref{Cor:Pre_Hauptresultat}, and discuss afterwards
previously known results in Remark~\ref{Rem:Known_Int_Besov_Sobolev} and~\ref{Rem:Known_QMC_Int_Disc}, which indeed show that our error bounds are optimal.

Then we present the desired results on 
numerical integration on our spaces of fractional smoothness and on the behavior of fractional discrepancies.

\begin{corollary}\label{Cor:Nets_Besov_Spaces}
Let $s\in \N$, $p,q\in [1,\infty]$, and $\alpha \in (0,1)$. Let $\alpha \ge 1/p$ if $q=1$ and $\alpha > 1/p$ if $q>1$.
Let $m\ge 1$, and let $P$ be a $(t,m,s)$-net in base $b$. 
Furthermore, put $N:=|P|=b^{m}$.
Then there exists a constant $C$, independent of $N$, such that
the QMC cubature rule $Q_P$, cf. \eqref{qmc-alg},  
satisfies
\begin{equation}
 \label{nets_besov}
e^{\wor}(Q_P, S^{\alpha}_{p,q}B([0,1]^s))
\le C N^{-\alpha}
\ln(N)^{\frac{s-1}{q'}}.
\end{equation}
In particular, we obtain for $q=1$ that the convergence rate of the worst-case error does not depend on the dimension $s$:
\begin{equation}
 \label{nets_besov_q=1}
e^{\wor}(Q_P, S^{\alpha}_{p,1}B([0,1]^s))
\le C N^{-\alpha}.
\end{equation}
\end{corollary}

\begin{proof}
The result follows immediately from Theorem~\ref{Nets}  and~\ref{tensor18a}.
\end{proof}

From  
Corollary~\ref{DirEmb_Multi} and~\ref{Cor:Nets_Besov_Spaces} we quickly obtain the 
following result. 

\begin{corollary}\label{Cor:Pre_Hauptresultat}
Let $s\in\N$, $1 < p < \infty$. and  $\alpha \in (1/p, 1)$. 
Let $m\ge 1$. Let $P$ be a $(t,m,s)$-net in base $b$ and denote the corresponding QMC-cubature  by $Q_P$ (cf.  \eqref{qmc-alg}). 
Furthermore, put $N:= |P| = b^m$.
\begin{itemize}
\item[{\rm (}i{\rm )}] Let additionally $\alpha > 1/2$. Then there exists a constant $C$, depending only on $\alpha$,
$s$, $p$, $b$, and $t$, such that
\begin{equation}\label{nets_sobolev}
e^{\wor}(Q_P, S^{\alpha}_{p}H([0,1]^s)) 
\le C N^{-\alpha} \ln(N)^{\frac{s-1}{2}}.
\end{equation}
\item[{\rm (}ii{\rm )}] Let additionally $\alpha \le 1/2$. Then for each $\widetilde{p} \in (\alpha^{-1}, p]$ there exists a constant $C$, depending only on $\alpha$,
$s$, $p$, $\widetilde{p}$, $b$, and $t$, such that
\begin{equation}\label{nets_sobolev_alpha_small}
e^{\wor}(Q_P, S^{\alpha}_{p}H([0,1]^s)) 
\le C N^{-\alpha} \ln(N)^{\frac{s-1}{\widetilde{p}'}}.
\end{equation}
\end{itemize}
\end{corollary}

\begin{proof}
Statement~(i) follows for $p\le 2$ from the second embedding in Corollary~\ref{Cor:Neu_Equation_60+} 
and Corollary~\ref{Cor:Nets_Besov_Spaces}. If $p>2$, we use first the embedding
$S^{\alpha}_{p}H([0,1]^s) \hookrightarrow S^{\alpha}_{2}H([0,1]^s)$, see Lemma~\ref{Win_Lemma_13}, and then the second embedding in Corollary~\ref{Cor:Neu_Equation_60+} 
as well as Corollary~\ref{Cor:Nets_Besov_Spaces} for the ``new parameter'' $p=2 > \alpha^{-1}$.\\
Statement~(ii) follows from first employing the embeddings $S^{\alpha}_{p}H([0,1]^s) \hookrightarrow S^{\alpha}_{\widetilde{p}}H([0,1]^s) \hookrightarrow S^{\alpha}_{\widetilde{p}, \widetilde{p}}B([0,1]^s)$, cf. Lemma~\ref{Win_Lemma_13} and Corollary~\ref{Cor:Neu_Equation_60+} and note that $\widetilde{p} > 2$, and Corollary~\ref{Cor:Nets_Besov_Spaces}.
\end{proof}


\begin{remark}\label{Rem:Known_Int_Besov_Sobolev} 
We denote by $X_N := \{\bsx^1, \ldots \, , \bsx^N\} \subset [0,1]^s$ a generic $N$-point set, and by
\[
\Lambda_N (f,X_N):= \sum_{i=1}^N \lambda_i \, f(\bsx^i)
\]
a generic cubature rule based on $X_N$ 
and on integration weights $\lambda_1, \ldots, \lambda_N \in \real$. 
For $F \subset C([0,1]^s)$ we put
\[
 \Lambda_N (F,X_N):= \sup_{f\in F}\, \bigg| 
I_s(f)  - \Lambda_N (f,X_N)\bigg| 
\]
and 
\[
 \kappa_N (F):= \inf_{X_N \subset [0,1]^s}\, \inf_{\lambda_1, \ldots \, \lambda_N \in \real}
 \, \Lambda_{N} (F,X_N).
\]
For better reference we need to distinguish between function classes defined on $[0,1]^s$ and the $s$-dimensional torus denoted by $\mathbb{T}^s$ and identified with $[0,1]^s$.
Whenever we consider a space defined on  $\mathbb{T}^s$ then this means that we only work with  in each variable $1$-periodic  functions. Recall that for the function spaces we have in mind, i.e., Bessel potential spaces and Besov spaces
with positive smoothness, the periodic spaces are always continuously embedded into their non-periodic counterpart.  
 First let $1\le p,q \le \infty$ and $\alpha >1/p$ if $q>1$ and $\alpha \ge 1/p$ if $q=1$. Furthermore,  let $F$ denote the norm unit ball in $S^\alpha_{p,q}B(\mathbb{T}^s)$. 
 Then it is known that 
 \begin{equation}\label{eq:kappa_besov}
 \kappa_N (F) \asymp N^{-\alpha}\, (\log N)^{(d-1)/q'}\, , 
 \end{equation}
 cf. \cite[Theorem~8.5.1]{DTU}. The lower error bound for arbitrary $q$  and $\alpha > 1/p$  was pointed out in \cite{Tem15}, and for $q=1$ and $\alpha=1/p$  in \cite{UU15}. For the upper error bound see \cite{NUU17}.\\ 
 Now let $1<  p < \infty$ and $\alpha > \max(1/p, 1/2)$, and let $F$ denote the norm unit ball in $S^\alpha_{p}H(\mathbb{T}^s)$. 
 Then it is known that 
 \begin{equation}\label{eq:kappa_sobolev}
 \kappa_N (F) \asymp N^{-\alpha}\, (\log N)^{(d-1)/2}\, , 
 \end{equation}
 cf. again \cite[Theorem~8.5.1]{DTU}. The lower error bound was proved in \cite{Tem15}, and the proof of the upper error bound can again be found  \cite{NUU17}. By switching back to the non-periodic spaces we are dealing with here, we observe that the lower bounds in \eqref{eq:kappa_besov}, \eqref{eq:kappa_sobolev} carry over thanks to the continuous embedding we mentioned before.

Note that the lower error bounds in \eqref{eq:kappa_besov} and \eqref{eq:kappa_sobolev}
imply that our upper error bounds in \eqref{nets_besov}, 
\eqref{nets_besov_q=1}, and \eqref{nets_sobolev} are best possible and that QMC rules based on arbitrary $(t,m,s)$-nets achieve in the corresponding settings the optimal convergence rates on periodic as well as on non-periodic Besov and Sobolev spaces of dominated mixed smoothness. \\
The upper bounds in  \eqref{eq:kappa_besov} and \eqref{eq:kappa_sobolev} were achieved by using 
\emph{Frolov cubature formulas}. Frolov cubature rules are, like QMC rules, equal-weight cubatures, but have the drawback that their weights, in general, do not sum up to $1$. That means that they do not integrate constant functions exactly. In contrast, cubature rules based on $(t,m,s)$-nets do integrate constant functions exactly. 
\end{remark}

\begin{remark}\label{Rem:Known_QMC_Int_Disc} 
As mentioned above, QMC cubature rules have the advantage that they are equal-weight rules and that their weights sum up to one. In \cite{Mar13} it was shown that the upper error bound~\eqref{nets_besov} 
holds for the Chen-Skriganov construction of $(t,m,s)$-nets $P$ for large base $b$, see~\cite{CS02, Skr06}, in the case where $1\le p,q \le \infty$ and $1/p < \alpha <1$.  \\
Note that we proved that the bound~\eqref{nets_besov} holds actually for arbitrary $(t,m,s)$-nets $P$, regardless of the size of the base $b$ or of the specific construction principle.
\end{remark}


\begin{corollary}\label{Hauptresultat}
Let $s\in\N$, $1 < p < \infty$, $1\le q \le \infty$, and  $\alpha \in (1/p, 1)$. 
Let $m\ge 1$. Let $P$ be a $(t,m,s)$-net in base $b$ and denote the corresponding QMC-cubature  by $Q_P$ (cf.  \eqref{qmc-alg}). 
Furthermore, put $N:= |P| = b^m$.
\begin{itemize}
\item[{\rm (}i{\rm )}] Let additionally $\alpha > 1/2$. Then there exists a constant $C$, depending only on $\alpha$,
$s$, $p$, $q$, $b$, and $t$, such that
\begin{equation}\label{frac_disc_estimate}
e^{\wor}(Q_P, \HH_{\alpha,s,p, q}) = D^*_{\alpha,s,p',q'}(P) 
\le C N^{-\alpha} \ln(N)^{\frac{s-1}{2}}.
\end{equation}
\item[{\rm (}ii{\rm )}] Let additionally $\alpha \le 1/2$. Then for each $\widetilde{p} \in (\alpha^{-1}, p]$ there exists a constant $C$, depending only on $\alpha$,
$s$, $p$, $\widetilde{p}$, $q$, $b$, and $t$, such that
\begin{equation}\label{frac_disc_estimate_2}
e^{\wor}(Q_P, \HH_{\alpha,s,p, q}) = D^*_{\alpha,s, p',q'}(P) 
\le C N^{-\alpha} \ln(N)^{\frac{s-1}{{\widetilde{p}}'}}.
\end{equation}
\end{itemize}
\end{corollary}

\begin{proof}
The statements follow immediately from Theorem~\ref{KoksmaHlawka},
Corollary~\ref{Winfried2},
and Corollary~\ref{Cor:Pre_Hauptresultat}.
\end{proof}

\begin{remark}\label{Rem:No_Ext_to_alpha_1}
In contrast to the situation where $\alpha < 1$, the upper bound \eqref{frac_disc_estimate} in Corollary~\ref{Hauptresultat} 
does not hold for \emph{arbitrary} $(t,m,s)$-nets $P$  in the case where $\alpha =1$.
Indeed, 
for $\alpha = 1$ and $1\le p < \infty$ the quantity $D^*_{1,s,p,p}(P)$ is essentially the $L_{p}$-star discrepancy (meaning if $T(P)$ denotes the point set resulting from applying the mapping $x\mapsto 1-x$ to each coordinate of all points of $P$, then $D^*_{1,s,p,p}(T(P))$ is the $L_{p}$-star discrepancy of $P$), and 
it is well-known that for $s\ge 2$ arbitrary $(t,m,s)$-nets $P$ in base $b$ do not necessarily satisfy \eqref{frac_disc_estimate} with a constant only depending on $s$, $p$, $b$, and $t$.
For instance, for $s=2$ the van der Corput set $P$ with $2^m$ points is a $(0,m,s)$-net in base $2$, but $D^*_{1,2,p,p}(T(P))$ is proportional to $N^{-1} \ln(N)$, see, e.g. \cite[Section~2.2]{Mat10}.

But note that for \emph{specific} $(t,m,s)$-nets the upper bound 
\eqref{frac_disc_estimate} 
extends to the case where $\alpha =1$.
Indeed, 
it was shown in \cite{Dav56} for $s=2$ and, with the help of probabilistic arguments, in \cite{Roth79, Roth80, Fro80, Che81} for arbitrary dimension $s$ that 
the $L_{p}$-star discrepancy of \emph{suitable} point sets $P$ satisfies \eqref{frac_disc_estimate}.
Note that the first explicit deterministic constructions for dimension $s\ge 3$, the Chen-Skriganov constructions already mentioned in Remark~\ref{Rem:Known_QMC_Int_Disc}, 
were published much later in \cite{CS02} for $p=2$ and in \cite{Skr06} for arbitrary $1<p< \infty$, and these constructions are highly non-trivial. Infinite sequences with optimal $L_p$-star discrepancy were studied in \cite{DHMP17}. The explicit constructions therein are based on order 2 digital nets, see \cite{DHMP17} for more details.

Moreover, it is well known that for $\alpha = 1$ the $L_{p}$-star-discrepancy and integration error bound  \eqref{frac_disc_estimate} is best possible for any point set $P$ and the QMC cubature rule based on $P$, respectively; see the matching lower discrepancy bounds provided in \cite{Roth54, Sch77}. 
\end{remark}


\section{Appendix - Tensor Products of Banach Spaces}
\label{Sec:Appendix}


In this appendix we provide the necessary background on tensor products of Banach spaces and  prove Theorem~\ref{tensor17} and \ref{tensor3}.


\subsection{Preliminaries}
\label{tpbs}


In this section we follow \cite{LiCh}, but see also \cite{DF}.
Let $X$ and $Y$ be Banach spaces. $X'$ denotes the dual space of $X$
and $\cl (X,Y)$ the space of continuous linear operators $T: X \to Y$.
Consider the set of all formal expressions
$\sum_{i=1}^n f_i \otimes g_i$, $n \in \N$, $f_i \in X$ and $g_i \in Y$.
We introduce an equivalence relation by means of
\[
\sum_{i=1}^n f_i \otimes g_i \sim \sum_{j=1}^m u_j \otimes v_j
\]
if 
\be\label{eq-111}
\sum_{i=1}^n \varphi(f_i) \, g_i = \sum_{j=1}^m \varphi(u_j) \, v_j \qquad
\mbox{for all}\quad \varphi \in X'\, .
\ee
Then the algebraic tensor product $X \otimes Y$ of $X$ and $Y$
is defined to be the set of all such equivalence classes.
One can equip 
the algebraic tensor product $X \otimes Y$
with several different cross norms, i.e., norms 
$\alpha(\cdot, X,Y)$ on $X\otimes Y$ that satisfy $\alpha(f\otimes g, X,Y) = \|f\|_X \|g\|_Y$ for all $f\in X$, $g\in Y$.
We are interested in so-called uniform cross norms only.
A uniform cross norm $\alpha$ assigns to each pair of Banach spaces $(X,Y)$ a cross norm $\alpha(\cdot, X,Y)$ such that if
$X_1,X_2,Y_1,Y_2$ are Banach spaces 
then
\[
\alpha\big((T_1\otimes T_2)h,Y_1,Y_2 \big) \leq \|T_1\|_{ \mathcal{L}(X_1,Y_1)}
\cdot
    \|T_2\|_{ \mathcal{L}(X_2,Y_2)}
    \cdot \alpha\big(h,X_1,X_2\big)
\]
holds for all $ h = \sum\limits_{j=1}^n f_j\otimes g_j \in X_1 \otimes X_2$
and all $T_1 \in \cl (X_1,Y_1)$, $T_2 \in \cl (X_2,Y_2)$,
where the tensor product $T_1\otimes T_2$ of the operators $T_1$ and $T_2$ is  
defined by
\be\label{eq-114}
(T_1 \otimes T_2) h := \sum_{i=1}^n (T_1 f_i)\otimes (T_2 g_i) \, , \qquad
h = \sum\limits_{i=1}^n f_i\otimes g_i \in X_1 \otimes X_2\, .
\ee
The completion of $X\otimes Y$ with respect to the tensor norm $\alpha$ will be denoted by
$X\otimes_{\alpha} Y$. If $\alpha$ is uniform, then  $T_1\otimes T_2$ has a unique
continuous linear 
extension to $X_1\otimes_{\alpha} X_2$ which we again denote by $T_1 \otimes T_2$.
Simple, but important, is the next property we need.

\begin{lemma}\label{tops1}
Let $X_1,X_2, Y_1,Y_2$ be Banach spaces, let
 $\alpha(\cdot,X,Y)$ be a uniform tensor norm, and let $T_1 \in \cl (X_1,Y_1)$ and $T_2 \in \cl (X_2,Y_2)$.
If $T_1$ and $T_2$ are linear isomorphisms, then the operator $T_1 \otimes T_2$ is a continuous linear isomorphism from
$X_1 \otimes_{\alpha} X_2$ onto $Y_1 \otimes_{\alpha} Y_2$.
If $T_1$ and $T_2$ are linear isometries, then $T_1\otimes T_2$ is a linear isometry from
$X_1 \otimes_{\alpha} X_2$ into $Y_1 \otimes_{\alpha} Y_2$.
\end{lemma}

\noindent
Next we recall one well-known construction of a tensor norm, namely the $p$-nuclear norm.

\begin{definition}\label{tensdefi}   Let $X$ and $Y$ be Banach spaces.
\begin{itemize}
\item[{\rm (}i{\rm )}]
Let $1 < p <  \infty$ and let $1/p+1/p' = 1$. Then the \emph{$p$-nuclear tensor norm
$\delta_p(\cdot,X,Y)$} or  \emph{left $p$-Chevet-Saphar norm} is given by
\begin{equation*}
\begin{split}
    &\delta_p(h,X,Y) := \\
    &~~~\inf\left\{\Big(\sum\limits_{i=1}^n \|f_i \|^p_{X}\Big)^{1/p}\cdot
    \sup\Big\{\Big(\sum\limits_{i=1}^n |\psi(g_i)\Big|^{p'}\Big)^{1/p'}:
    \psi \in Y', \|\psi \|_{Y'} \leq 1 \Big\}\right\}\,,
  \end{split}
\end{equation*}
where the infimum is taken over all representations of $h$ of the form
\[
h = \sum_{i=1}^n f_i \otimes g_i\, , \qquad f_i \in X\, , \quad g_i \in Y\, , \quad i=1, \ldots \, ,n\, .
\]
\item[{\rm (}ii{\rm )}] 
The \emph{projective tensor norm} $\delta_1(\cdot,X,Y)$ is defined by
\[
\delta_1(h,X,Y) := \inf \Big\{ \sum\limits_{j = 1}^n
\, \|f_j\|_{X} \, \|g_j\|_{Y}: \quad f_j \in X, \, g_j \in Y, \,
h = \sum\limits_{j=1}^n f_j \otimes g_j\Big\}\, .
\]
\end{itemize}
\end{definition}

\begin{remark}\label{tensrem1}
\rm
The expression $\delta_p(\cdot,X,Y)$ defines a norm; we refer to
\cite[Chapt.~1]{LiCh}.
In the case where $1<p<\infty$ one can replace
\be\label{f6}
    \sup\Big\{\Big(\sum\limits_{i=1}^n |\psi(g_i)|^{p'}\Big)^{1/p'}:
    \psi\in Y', \|\psi \|_{Y'} \leq 1 \Big\}
\ee
by
\be\label{f5}
    \sup\Big\{\Big\|\sum\limits_{i=1}^n \lambda_i g_i \Big\|_{Y}:
    \Big(\sum\limits_{i=1}^n |\lambda_i|^p\Big)^{1/p}\leq 1\Big\}\,,
\ee
see \cite[Lem.~1.44]{LiCh}. 
\end{remark}

\noindent
Tensor products of functions and sequences have been the original source
for the introduction of the abstract tensor product.
Here we only deal with those ``classical'' tensor products.

\begin{remark}\label{Rem:Tensor_L_p}
Let $1\le p< \infty$ and $s\ge 2$. Then we have
\be\label{tensorlp}
 L_{p} ([0,1]^{s-1}) \otimes_{\delta_p} L_p ([0,1]) =  L_p ([0,1]^s) 
\ee
with equality of norms, see \cite[Corollary~1.52]{LiCh} or \cite{DF}. 
Furthermore, it is trivial to verify that we have 
$$L_p([0,1]^s) \otimes_{\delta_p} \RR = \RR \otimes_{\delta_p} L_p([0,1]^s)  
= L_p( [0,1]^s)$$ and
$  \RR \otimes_{\delta_p} \RR = \RR$ with equalities of norms.
\end{remark}

\subsection{Tensor Products of Sequence Spaces and Haar Wavelet Spaces} 
\label{APP:TPSS}

In this subsection we provide a result on tensor products of sequence spaces and the proof
of Theorem~\ref{tensor17}, which deals with tensor products of Haar wavelet spaces.


Let $a^\alpha_{s,p,q}$ denote the collection of all sequences $\lambda =(\lambda_{\bsj,\bsk,\bsi})_{\bsj,\bsk,\bsi}$
such that
\begin{equation*}
\|\lambda\|^q_{a^\alpha_{s,p,q}} :=
\sum^\infty_{L=0} b^{q(\alpha-1/p+1/2)L} \sum_{|\bsj|=L}
\left( \sum_{\bsk\in \Delta_{\bsj -{\bf 1}}} 
\sum_{\bsi \in \nabla_{\bsj}}
|\lambda_{\bsj,\bsk,\bsi}|^p \right)^{q/p} <\infty\, ,
\end{equation*}
where $\Delta_{\bsj -{\bf 1}} $ and $\nabla_{\bsj}$ are depending on $s$, see \eqref{neu-ws-1}.
Obviously these sequence spaces are closely related to the definition of our wavelet spaces 
$\HH_{\wav,\alpha,s,p,q}$.
As a consequence of results on tensor products of weighted $\ell_p$ spaces, see \cite{SU1},
one can derive the following.

\begin{lemma}\label{tensor5}
 Let $1 \le p < \infty$, $\alpha \in \re$ and $s\in \N$, $s\ge 2$.
 Then
 \[
 a^\alpha_{s-1,p,p}  \otimes_{\delta_p} a^\alpha_{1,p,p} =  a^\alpha_{s,p,p} 
 \]
 with equality of norms. 
\end{lemma}

\begin{remark}\label{Rem:Wav_Folgenraum}
Let $\alpha >0$.
The mapping 
$$\iota_1: \HH_{\wav,\alpha,p,q} \to a^\alpha_{1,p,q}\,,\, f\mapsto (\langle f, \psi^j_{i,k}\rangle )_{j,k,i}$$
and its $s$-variate generalization
$$\iota_s: \HH_{\wav,\alpha,s,p,q} \to a^\alpha_{s,p,q}\,,\, f\mapsto (\langle f, \Psi^{\bsj}_{\bsi,\bsk}\rangle )_{\bsj,\bsk,\bsi}$$
are obviously linear isometries, but, due to \eqref{fourier_sum_zero}, not surjective.
Consequently, the corresponding tensor operator
\[
\iota_{s-1} \otimes \iota_1: \HH_{\wav, \alpha,s-1,p,p}  \otimes_{\alpha_p}   \HH_{\wav, \alpha,p,p}
\to a^\alpha_{s-1,p,p}  \otimes_{\alpha_p} a^\alpha_{1,p,p} 
\]
is, due to Lemma~\ref{tops1}, also a linear isometry, but in general not surjective.
\end{remark}



\begin{proof}[Proof of Theorem~\ref{tensor17}] 
Step 1. Let $s=2$.

We first prove $
  \HH_{\wav,\alpha,1,p,p} \otimes_{\delta_p}  \HH_{\wav,\alpha,1,p,p} \hookrightarrow
   \HH_{\wav,\alpha,2,p,p} \, .$
From Lemma \ref{lemmalp} we know that $\HH_{\wav,\alpha,p,p} \hookrightarrow L_p ([0,1])$.
Together with Lemma~\ref{tops1}
and \eqref{tensorlp} 
this implies  
\[
 \HH_{\wav,\alpha,p,p} \otimes_{\delta_p} \HH_{\wav,\alpha,p,p} \hookrightarrow 
 L_p ([0,1])\otimes_{\delta_p}  L_p ([0,1]) = L_p ([0,1]^2)\, .
\]
In particular, we conclude that 
$ \HH_{\wav,\alpha,p,p} \otimes_{\delta_p} \HH_{\wav,\alpha,p,p} 
\subseteq L_1 ([0,1]^2)$.
 
Let $h\in \HH_{\wav,\alpha,p,p} \otimes \HH_{\wav,\alpha,p,p}$.
By definition 
\[
 h= \sum_{m=1}^N  f^{(m)} \otimes g^{(m)}   
\]
for some suitable $f^{(m)}, g^{(m)}\in \HH_{\wav,\alpha,p,p}$.
Then we have for all $\bsj \in \N^2_0$, $\bsk \in \Delta_{\bsj -{\bf 1}} $,  and $\bsi \in \nabla_{\bsj}$
\[
 \lambda_{\bsj,\bsk,\bsi} := 
 \Big\langle \sum_{m=1}^N f^{(m)} \otimes g^{(m)}, \Psi^{\bsj}_{\bsi,\bsk} \Big\rangle
=  \sum_{m=1}^N 
\langle f^{(m)}, \psi^{j_1}_{i_1,k_1} \rangle \cdot 
\langle g^{(m)}, \psi^{j_2}_{i_2,k_2} \rangle.
\]
Due to Remark~\ref{Rem:Wav_Folgenraum} and Lemma~\ref{tensor5}
\begin{equation*}
 \|h \|_{\HH_{\wav,\alpha,2,p,p}} = \|\lambda\|_{a^\alpha_{2,p,p}} 
  = \left\| \lambda  \right\|_{ a^\alpha_{1,p,p}  \otimes_{\delta_p} a^\alpha_{1,p,p} } 
=  \left\| \sum_{m=1}^N  f^{(m)} \otimes g^{(m)} \right\|_{ \HH_{\wav,\alpha,p,p} \otimes_{\delta_p}  \HH_{\wav,\alpha,p,p}}\,.
\end{equation*}
Since the algebraic tensor product $\HH_{\wav,\alpha,p,p} \otimes \HH_{\wav,\alpha,p,p}$
is dense in $\HH_{\wav,\alpha,p,p} \otimes_{\delta_p} \HH_{\wav,\alpha,p,p}$, we obtain the isometric embedding
\begin{equation}\label{isom_embed}
 \HH_{\wav,\alpha,p,p} \otimes_{\delta_p} \HH_{\wav,\alpha,p,p} \hookrightarrow 
\HH_{\wav,\alpha,2,p,p}.
\end{equation}

Now we prove $
\HH_{\wav,\alpha,2,p,p} \hookrightarrow  \HH_{\wav,\alpha,1,p,p} \otimes_{\delta_p}  \HH_{\wav,\alpha,1,p,p}$.
Let $h\in \HH_{\wav,\alpha,2,p,p}$, i.e., $h\in L_1 ([0,1]^2)$ and 
\[
 h =
\sum^\infty_{\bsj \in \N_0^2} 
\sum_{\bsk\in \Delta_{\bsj-1}}
\sum_{\bsi \in \nabla_{\bsj}} 
\lambda_{\bsj,\bsk,\bsi} 
\, \Psi^{\bsj}_{\bsi,\bsk} 
\]
with 
 $\lambda_{\bsj,\bsk,\bsi} = \langle h, \Psi^{\bsj}_{\bsi,\bsk} \rangle$, 
 cf. Lemma~\ref{Rem_meaningful}.
 Recall that the sum on the right hand side converges unconditionally in $\HH_{\wav,\alpha,2,p,p}$ to $h$, cf. again Lemma~\ref{Rem_meaningful}.
For $L\in \N$ we set $\Lambda: = \{ \bsj \in \N_0^2\,:\, j_1,j_2 \le L \}$.  
Additionally, we put 
$w(j):= b^{(\alpha-1/p + 1/2)j}$.
Then the finite sum $h_L:= \mathcal{S}_{\Lambda}(h)$ (cf. \eqref{partial_sum_S} ) satisfies
\beqq
 h_L & = &
\sum_{j_1=0}^L \sum_{k_1\in \Delta_{j_1-1}} \sum_{i_1 \in \nabla_{j_1}} 
\psi^{j_1}_{i_1,k_1} \, \otimes \Big(
\sum_{j_2=0}^L   \sum_{k_2\in \Delta_{j_2-1}}
\sum_{i_2 \in \nabla_{j_2}}  \lambda_{\bsj,\bsk,\bsi} 
\,  \psi^{j_2}_{i_2,k_2} \Big)
\\
&=&
\sum_{j_1=0}^L \sum_{k_1\in \Delta_{j_1-1}} \sum_{i_1 \in \nabla_{j_1}}
\bigg( \frac{1}{w(j_1)}\, 
\psi^{j_1}_{i_1,k_1}  \otimes \Big(  w(j_1)
\sum_{j_2=0}^L   \sum_{k_2\in \Delta_{j_2-1}}
\sum_{i_2 \in \nabla_{j_2}}  \lambda_{\bsj,\bsk,\bsi} 
\,  \psi^{j_2}_{i_2,k_2}\Big)\bigg)\, .
\eeqq
We put
\[
g^{j_1,k_1,i_1} := \frac{1}{w(j_1)}\,  \psi^{j_1}_{i_1,k_1}
\]
and 
\[
 f^{j_1,k_1,i_1}:=  w(j_1)\, \sum_{j_2=0}^L    \sum_{k_2\in \Delta_{j_2-1}}
\sum_{i_2 \in \nabla_{j_2}}  \lambda_{\bsj,\bsk,\bsi} 
\,  \psi^{j_2}_{i_2,k_2}\, .
\]
Clearly, $g^{j_1,k_1,i_1},  f^{j_1,k_1,i_1} \in \HH_{\wav,\alpha,p,p} \subseteq L_1([0,1])$.
It follows from \eqref{HHY7} that 
\be\label{tensor13}
\Big\|  g^{j_1,k_1,i_1} \Big\|_{\HH_{\wav,\alpha,p,p}}^{p} 
= \left\{
\begin{array}{lll}
(b-1) b^{-p} + (1-1/b)^p &\qquad & \mbox{if}\quad j_1\in \N\, , 
\\
1 && \mbox{if}\quad j_1=0\,.
\end{array}
\right.
\ee
Furthermore,
\beq\label{tensor14}
\Big\|  f^{j_1,k_1,i_1}
 \Big\|_{\HH_{\wav,\alpha,p,p}}^{p}
& = & 
\sum_{j_2=0}^L b^{j_2(\alpha-1/p+1/2)p}  \sum_{k_2\in \Delta_{j_2-1}}
\sum_{i_2 \in \nabla_{j_2}}   | w(j_1) \lambda_{\bsj,\bsk,\bsi}|^p \,.
\eeq
Recall, if $1<p<\infty$ the  tensor norm
$\delta_p(\cdot,X,Y)$ is given by
\[
\delta_p(h,X,Y) := \inf\left\{\Big(\sum\limits_{i=1}^n \|f_i \|_X^p\Big)^{1/p}\cdot
 \sup\Big\{\Big\|\sum\limits_{i=1}^n \mu_i g_i \Big\|_Y:
    \Big(\sum\limits_{i=1}^n |\mu_i|^p\Big)^{1/p}\leq 1\Big\}   
\right\}\,,
\]
where the infimum is taken over all representations of $h$ of the form
\[
h = \sum_{i=1}^n f_i \otimes g_i\, , \qquad f_i \in X\, , \quad g_i \in Y\, , \quad i=1, \ldots \, ,n\, .
\]
The two estimates \eqref{tensor13} and \eqref{tensor14}
yield
\beqq
&&\hspace*{-0.7cm}
\delta_p(h_{L},\HH_{\wav,\alpha,p,p},\HH_{\wav,\alpha,p,p})^p  
\\
& \le &  
\sum_{j_1=0}^L \sum_{k_1\in \Delta_{j_1-1}} \sum_{i_1 \in \nabla_{j_1}} 
\sum_{j_2=0}^L  b^{j_2(\alpha-1/p + 1/2)p} \sum_{k_2\in \Delta_{j_2-1}}
\sum_{i_2 \in \nabla_{j_2}}  |w(j_1)\, \lambda_{\bsj,\bsk,\bsi}|^p 
\\
&& \times \sup  \left\{\Big\|
\sum_{j_1=0}^L \sum_{k_1\in \Delta_{j_1-1}} \sum_{i_1 \in \nabla_{j_1}} 
\mu_{j_1,k_1,i_1}  \, g^{j_1,k_1,i_1}\Big\|_{\HH_{\wav,\alpha,p,p}}^{p}:~ \|\mu\|_{\ell_p}\le 1
\right\}\, .
\eeqq
Observe that 
\beqq
&&
\hspace*{-0.7cm}
\Big\|
\sum_{j_1=0}^L \sum_{k_1\in \Delta_{j_1-1}} \sum_{i_1 \in \nabla_{j_1}} 
\mu_{j_1,k_1,i_1}  \, g^{j_1,k_1,i_1}\Big\|^p_{\HH_{\wav,\alpha,p,p}}
=
\sum_{j_1=0}^L   \sum_{k_1\in \Delta_{j_1-1}}
\sum_{i_1 \in \nabla_{j_1}}  \left| \mu_{j_1,k_1,i_1}  - \frac{1}{b} \sum_{i \in \nabla_{j_1}} \mu_{j_1,k_1,i}  \right|^p \, .
\eeqq
Due to Minkowski's inequality and Jentzen's inequality we obtain
\begin{equation*}
\begin{split}
&\left( \sum_{j_1=0}^L   \sum_{k_1\in \Delta_{j_1-1}}
\sum_{i_1 \in \nabla_{j_1}}  \left| \mu_{j_1,k_1,i_1}  - \frac{1}{b} \sum_{i \in \nabla_{j_1}} \mu_{j_1,k_1,i}  \right|^p \right)^{1/p} \\
\le &\left( \sum_{j_1=0}^L   \sum_{k_1\in \Delta_{j_1-1}}
\sum_{i_1 \in \nabla_{j_1}}  \left| \mu_{j_1,k_1,i_1}  \right|^p \right)^{1/p}
+ \left( \sum_{j_1=0}^L   \sum_{k_1\in \Delta_{j_1-1}} b 
 \left|  \frac{1}{b} \sum_{i \in \nabla_{j_1}} \mu_{j_1,k_1,i}  \right|^p \right)^{1/p}\\
 \le &2 \left( \sum_{j_1=0}^L   \sum_{k_1\in \Delta_{j_1-1}}
\sum_{i_1 \in \nabla_{j_1}}  \left| \mu_{j_1,k_1,i_1}  \right|^p \right)^{1/p}
\le 2 \,.
\end{split}
\end{equation*} 
Hence
\beqq
&&\hspace*{-0.7cm} 
\delta_p(h,\HH_{\wav,\alpha,p,p},\HH_{\wav,\alpha,p,p})^p
\\
&\le & 
2^{p} \sum_{j_1=0}^L \sum_{k_1\in \Delta_{j_1-1}} \sum_{i_1 \in \nabla_{j_1}} 
\sum_{j_2=0}^L  b^{p(\alpha-1/p + 1/2)(j_1+j_2)} \sum_{k_2\in \Delta_{j_2-1}}
\sum_{i_2 \in \nabla_{j_2}}  |\lambda_{\bsj,\bsk,\bsi}|^p
\\
&\le & 2^{p} \|h \|^p_{\HH_{\wav,\alpha,2,p,p}}\, .
\eeqq
This proves 
\begin{equation}\label{invers_embed}
\HH_{\wav,\alpha,2,p,p} \hookrightarrow \HH_{\wav,\alpha,p,p} \otimes_{\delta_p} \HH_{\wav,\alpha,p,p}
\end{equation}
if $1<p<\infty$.
\\
It remains to deal with $p=1$.
Recall that the projective tensor norm $\delta_1 (h, X, Y)$ is defined by
\[
\delta_1(h, X, Y) = \inf \left\{
\sum_{m=1}^n \|f^m \|_X \cdot \|g^m\|_Y :~ f^m \in X, g^m \in Y, ~h = \sum_{m=1}^n
f^m \otimes  g^m \right\}\,  .
\]
In our concrete situation we obtain from \eqref{tensor13} and \eqref{tensor14}
\beqq
&&\hspace*{-0.7cm}
\delta_1(h,\HH_{\wav,\alpha,1,1},\HH_{\wav,\alpha,1,1})  
\\
& \le &  \frac{2(b-1)}{b} \, 
\sum_{j_1=0}^L \sum_{k_1\in \Delta_{j_1-1}} \sum_{i_1 \in \nabla_{j_1}} 
\sum_{j_2=0}^L  b^{j_2(\alpha-1/2)}\sum_{k_2\in \Delta_{j_2-1}}
\sum_{i_2 \in \nabla_{j_2}}  |w(j_1)\, \lambda_{\bsj,\bsk,\bsi}| 
\\
& \le & \frac{2(b-1)}{b} \, \|h\|_{\HH_{\wav,\alpha,2,1,1}}\, .
\eeqq
Hence we obtain \eqref{invers_embed} also in the case $p=1$. 
Notice that \eqref{invers_embed} implies, in particular, that the embedding in 
 \eqref{isom_embed} is surjective; since it is an isometry, we see the norms in $\HH_{\wav,\alpha,2,p,p}$ and $\HH_{\wav,\alpha,p,p} \otimes_{\delta_p} \HH_{\wav,\alpha,p,p}$ are actually equal. 
This finishes the proof in case $s=2$.\\
Step 2. 
The remaining cases $s\ge 3$ can be proved by induction on $s$.
\end{proof}


\subsection{Tensor Products of 
Sobolev-Type Spaces} 
\label{APP:TPSO}


In this subsection we provide 
 the proof
of Theorem~\ref{tensor3}, which deals with tensor products of spaces of fractional smoothness.



\begin{proof}[Proof of Theorem~\ref{tensor3}]
For notational convenience, we provide the proof only for the case $s=2$.
From this it will become obvious that the proof approach works also for general $s\ge 2$.\\
{\em Step 1.}
To begin with we shall verify the continuous embedding
\be\label{stp1H}
 \HH_{\alpha,p} \otimes_{\delta_p}\HH_{\alpha,p}
\hookrightarrow \HH_{\alpha,2,p} \, .
\ee
The definition of the norm $\|\cdot\|_{\alpha,p,p}$ implies that the linear mappings
$$\delta_{0}:  \HH_{\alpha,p,p} \to \RR\,,\, f \mapsto f(0)$$
and 
$$D^{\alpha}:  \HH_{\alpha,p,p} \to L_p([0,1])\,,\, f \mapsto D^{\alpha}f$$
are continuous with operator norms $\|\delta_0\| = 1$ and $\| D^{\alpha} \| = \Gamma(\alpha)$. 
Now let $h = \sum_{i=1}^n f_i\otimes g_i $, $f_i,  g_i \in  \HH_{\alpha,p}$, $i=1, \ldots \, , n$.
Due to Lemma~\ref{Lemma0} we have for $\tilde{f_i} : = D^{\alpha} f_i$, $\tilde{g_i} := D^{\alpha}g_i$, $i=1,\ldots,n$,
\beqq
h(x,y) & = & \sum_{i=1}^n \Big(f_i (0) + \frac{1}{\Gamma (\alpha)} \, \int_0^1 \tilde{f}_i(t)\,  (x-t)_+^{\alpha-1}\, dt\Big) \, \cdot \,  
 \Big(g_i (0) + \frac{1}{\Gamma (\alpha)} \, \int_0^1 \tilde{g}_i(t)\,  (y-t)_+^{\alpha-1}\, dt\Big)
\\
& = & \sum_{i=1}^n f_i (0)\, g_i (0)  + 
\sum_{i=1}^n \Big[\frac{f_i(0)}{\Gamma (\alpha)} \, \int_0^1 \tilde{g}_i(t)\,  (y-t)_+^{\alpha-1}\, dt
 + \frac{g_i (0)}{\Gamma (\alpha)} \, \int_0^1 \tilde{f}_i(t)\,  (x-t)_+^{\alpha-1}\, dt\Big]
 \\
& + & \quad \sum_{i=1}^n \frac{1}{\Gamma (\alpha)^2} \int_{[0,1]^2} \,  \tilde{f}_i (t_1) \tilde{g}_i(t_2)\,  (x-t_1)_+^{\alpha-1} \, (y-t_2)_+^{\alpha-1}\, dt_1 \, dt_2\, , 
  \eeqq
$x,y \in (0,1)$.
Hence we obtain
\[
h (x_1,x_2) = 
\sum_{u \subseteq \{1,2\}} \Gamma (\alpha)^{-|u|} \, \int_{[0,1]^u} \tilde{h}_{u} (t_u) \prod_{j \in u}   (x_j-t_j)_+^{\alpha-1} \,  dt_u \, ,
\]
where
\beqq
\tilde{h}_{\emptyset} &:= & \sum_{i=1}^n f_i (0)\, g_i (0) = (\delta_0 \otimes \delta_0)h\,
\\
\tilde{h}_{\{1\}} (t) &:= & \sum_{i=1}^n g_i (0) \,  \tilde{f}_i(t) = (D^{\alpha} \otimes \delta_0)h
\, ,
\\
\tilde{h}_{\{2\}} (t) &:= & \sum_{i=1}^n f_i (0) \,  \tilde{g}_i(t) = (\delta_0 \otimes D^{\alpha})h\,   , 
\\
\tilde{h}_{\{1,2\}} (t_1,t_2) &:= & \sum_{i=1}^n \,  \tilde{f}_i (t_1) \tilde{g}_i(t_2) = (D^{\alpha} \otimes D^{\alpha})h\,  .
\eeqq
Consequently,
\beqq
\| \, h \, \|_{\alpha,2,p,p}^p  & = &
\sum_{u \subseteq \{1,2\}} \Gamma (\alpha)^{-|u|p} \| \, \tilde{h}_{u} \, \|_{L_p ([0,1]^u)}^p \\
& = & 
| (\delta_0 \otimes \delta_0)h |^p + \Gamma(\alpha)^{-p} 
\| (D^{\alpha} \otimes \delta_0)h \|^p_{L_p ([0,1])} \\
& & + \Gamma(\alpha)^{-p} \| (\delta_0 \otimes D^{\alpha})h \|^p_{L_p ([0,1])} 
+ \Gamma(\alpha)^{-2p} \| (D^{\alpha} \otimes D^{\alpha})h \|^p_{L_p ([0,1]^2)}\\
& \leq & 4 \delta_p (h,  \HH_{\alpha,p,p},  \HH_{\alpha,p,p})^p\,,
\eeqq
where the final inequality follows from Remark~\ref{Rem:Tensor_L_p}
and the fact that $\delta_p$ is a uniform tensor norm.

Thus we have the norm bound
\begin{equation}\label{neu-ws-2}
\| \, h \, \|_{\alpha,2,p,p} \le 4^{1/p}  \delta_p(h,  \HH_{\alpha,p,p},  \HH_{\alpha,p,p})\,.
\end{equation}

{\em Step 2.} We still have to prove
\be\label{stp2H}
 \HH_{\alpha,2,p} 
\hookrightarrow  \HH_{\alpha,p} \otimes_{\delta_p}\HH_{\alpha,p}\, .
\ee
Let 
\beqq
h (x_1,x_2) &= &
\sum_{u \subseteq \{1,2\}} \Gamma (\alpha)^{-|u|} \, \int_{[0,1]^u} 
\tilde{h}_{u} (t_u) \prod_{j \in u}   (x_j-t_j)_+^{\alpha-1} \,  dt_u \, .
\\
&=& h(0,0) + \Gamma (\alpha)^{-1} \, \int_{[0,1]} 
\tilde{h}_{1} (t) (x_1-t)_+^{\alpha-1} \,  dt
+
\Gamma (\alpha)^{-1} \, \int_{[0,1]} 
\tilde{h}_{2} (t)    (x_2-t)_+^{\alpha-1} \,  dt
\\
&+&
\Gamma (\alpha)^{-2} \, \int_{[0,1]^2} 
\tilde{h}_{1,2} (t_1,t_2) \, (x_1-t_1)_+^{\alpha-1} \, (x_2-t_2)_+^{\alpha-1} \,  dt_1\, dt_2
\eeqq
be a  given function in $\HH_{\alpha,2,p} $.
We need to find an appropriate decomposition of this function.
Temporarily we switch to a dense subset in  $\HH_{\alpha,2,p} $ by assuming that
\be\label{w5}
\tilde{h}_{1,2} (t_1,t_2) = \sum_{i=1}^n \tilde{f}_i (t_1)\, \tilde{g}_i(t_2) 
\ee
where $\tilde{f}_i , \, \tilde{g_i} \in L_p ([0,1])$, $i=1, \ldots \, n$, 
\[
\| \tilde{g}_i\|_{L_p} = \Gamma (\alpha)\qquad \mbox{and}\qquad
\supp \tilde{g}_i \cap \supp \tilde{g}_j = \emptyset\, , \quad i \neq j\, . 
\]
Furthermore, we define
\beqq
f_j (x_1)& := & 
\Gamma (\alpha)^{-1} \, \int_{[0,1]} 
\tilde{f}_{j} (t)    (x_1-t)_+^{\alpha-1} \,  dt\, , \qquad j=1, \ldots \, n\, ,
\\
g_j (x_2)& := & 
\Gamma (\alpha)^{-1} \, \int_{[0,1]} 
\tilde{g}_{j} (t)    (x_2-t)_+^{\alpha-1} \,  dt\, , \qquad j=1, \ldots \, n\, ,
\\
f_{n+1} (x_1) &:= & \Gamma (\alpha)^{-1} \, \int_{[0,1]} 
\tilde{h}_{1} (t) (x_1-t)_+^{\alpha-1} \,  dt\, , \qquad \tilde{f}_{n+1} := \tilde{h}_1\, , 
\\
g_{n+1} (x_2) &:= & 1\, , \qquad \tilde{g}_{n+1} := 0\, , 
\\
f_{n+2} (x_1) &:=& \frac{1}{c}\, ,  \qquad \tilde{f}_{n+2} := 0\, , \\
g_{n+2} (x_2) &:=&  c\, \Gamma (\alpha)^{-1} \, 
\int_{[0,1]} \tilde{h}_{2} (t)    (x_2-t)_+^{\alpha-1} \,  dt \, , 
\qquad \tilde{g}_{n+2} := c\, \tilde{h}_2\, , 
\\
f_{n+3} (x_1) &:= & h(0,0)\, ,  \qquad \tilde{f}_{n+3} := 0\, , 
\\
g_{n+3} (x_2) &:= & 1\, , \qquad \tilde{g}_{n+3} := 0 \, ,
\eeqq
where we choose $c := \Gamma (\alpha)/\|\tilde{h_2}\|_{L_p (0,1)}$ if $\|\tilde{h_2}\|_{L_p (0,1)}>0$.
Otherwise the modifications will be obvious.
It follows
\[
h = \sum_{i=1}^{n+3} f_i \otimes g_i\, , \qquad f_i, \, g_i \in \HH_{\alpha,p, p}
\, , \qquad i =1, \ldots \, , n+3\, .
\]
Since $f_j(0)=0$ for $j=1,\ldots,n+1$, and $g_j(0)=0$ for $j=1,\ldots,n$ and $j=n+2$, we get
\begin{eqnarray}\label{f4}
&& \hspace{-0.7cm} \alpha_p(h,  \HH_{\alpha,p,p} , \HH_{\alpha,p,p} )^p
\\
&\leq& \Big(
\sum\limits_{i=1}^{n+3} |f_i (0)|^p + \| \tilde{f}_i/\Gamma (\alpha) \|_{L_p(0,1)}^p\Big) \cdot
\sup\limits_{\|\lambda\|_{\ell_p}\leq 1} \, \left\{
\Big|\sum\limits_{i=1}^{n+3} \lambda_i\, g_i(0)\Big|^p  + 
\Big\|\sum\limits_{i=1}^{n+3} \lambda_i\, \tilde{g}_i/\Gamma (\alpha) \Big\|^p_{L_p(0,1)}\right\}
\nonumber
\\
&=& \Big(\| \tilde{h}_2/\Gamma (\alpha) \|_{L_p(0,1)}^p +|h(0,0)|^p + \| \tilde{h}_1/\Gamma (\alpha) \|_{L_p(0,1)}^p +
\sum\limits_{i=1}^{n}  \| \tilde{f}_i/\Gamma (\alpha) \|_{L_p(0,1)}^p\Big) 
\nonumber
\\
&&
\times 
\sup\limits_{\|\lambda\|_{\ell_p}\leq 1} \, \left\{
|\lambda_{n+1} + \lambda_{n+3}|^p  + 
\Big\| c\, \lambda_{n+2} \, 
\tilde{h}_{2}/\Gamma (\alpha) + \, \sum\limits_{i=1}^n \lambda_i\, 
\tilde{g}_i/\Gamma (\alpha) \Big\|^p_{L_p(0,1)}\right\}\, .
\nonumber
\end{eqnarray}
Observe, by the restrictions on the supports of the $\tilde{g}_i$ and the choice of 
$\| \tilde{g}_i\|_{L_p (0,1)}$, 
\beqq
|\lambda_{n+1} + \lambda_{n+3}|^p   &+ &
\Big\| c\, \lambda_{n+2} \, \tilde{h}_{2}/\Gamma (\alpha) + \, \sum\limits_{i=1}^n \lambda_i\, 
\tilde{g}_i/\Gamma (\alpha) \Big\|^p_{L_p(0,1)}
\\
&\le &
2^p + \Big( 1 + \Big(\frac{1}{\Gamma (\alpha)^p}\, \sum_{i=1}^n |\lambda_i|^p \, 
\| \tilde{g}_i \|_{L_p (0,1)}^p \Big)^{1/p}\Big)^p \le 2^{p+1} \, .
\eeqq
Furthermore, using again the disjointnes of the supports of the $\tilde{g}_i$,
\beqq
\sum\limits_{i=1}^{n}  \| \tilde{f}_i/\Gamma (\alpha) \|_{L_p(0,1)}^p
& = & \sum\limits_{i=1}^{n}  \| \tilde{f}_i/\Gamma (\alpha) \|_{L_p(0,1)}^p\, 
\| \tilde{g}_i/\Gamma (\alpha) \|_{L_p(0,1)}^p
\\
& = & \Gamma (\alpha)^{-2p}\, 
\sum\limits_{i=1}^{n}  \int_{\supp \tilde{g}_i} 
\int_{[0,1]} |\tilde{f}_i (t_1)\, \tilde{g}_i (t_2) |^p\, dt_1\, dt_2
\\
& = & \Gamma (\alpha)^{-2p}\,  
\int_{[0,1]^2} \Big|\sum_{i=1}^n \tilde{f}_i (t_1)\, \tilde{g}_i (t_2) \Big|^p\, dt_1\, dt_2\, .
\eeqq
Both together inserted into \eqref{f4}, yields
\beqq
&& \delta_p(h,  \HH_{\alpha,p,p} , \HH_{\alpha,p,p} )^p
\\
& \leq & 2^{p+1} \, \Big(\| \tilde{h}_2/\Gamma (\alpha) \|_{L_p(0,1)}^p +|h(0,0)|^p + 
\| \tilde{h}_1/\Gamma (\alpha) \|_{L_p(0,1)}^p +
\| \tilde{h}_{1,2}/\Gamma^2 (\alpha) \|_{L_p(0,1)}^p\Big)
\\
&= & 2^{p+1}\, \| \, h \, \|_{\alpha,2,p,p}^p \, .
\eeqq
We conclude
\be\label{tensor-eq2}
\delta_p(h,  \HH_{\alpha,p ,p} , \HH_{\alpha,p,p} ) 
 \le 2^{1+1/p}\, \| \, h \, \|_{\alpha,2,p,p}\, .
\ee
for all $h$ where $\tilde{h}_{1,2}$ is as in \eqref{w5}.
A density argument completes the proof of the claim.
\end{proof}

\begin{remark}
 \rm
The inequalities \eqref{neu-ws-2} and \eqref{w5} yield
\be\label{w6}
\frac{1}{4^{1/p}}\, \| \, h \, \|_{\alpha,2,p,p}\le
 \delta_p(h,  \HH_{\alpha,p ,p} , \HH_{\alpha,p,p} ) 
 \le 2^{1+1/p}\, \| \, h \, \|_{\alpha,2,p,p}
\ee
for all $h$. It would be of certain interest to clarify whether
the norms $\| \, \cdot \, \|_{\alpha,2,p,p}$ and 
$\|\, \cdot \, \|_{\HH_{\alpha,p} \otimes_{\delta_p}\HH_{\alpha,p}}$
are not only equivalent but coincide.
\end{remark}

\begin{proof}[Proof of Corollary~\ref{Winfried2}] 
From the definition of the  tensor norm $\delta_p(\cdot,X,Y)$ for $1<p<\infty$, cf. Definition~\ref{tensdefi}, 
it is immediate that $X_0 \hookrightarrow X$ and 
$Y_0 \hookrightarrow Y$ with corresponding embedding constants $C_X, C_Y>0$ imply
\[
\delta_p(h,X,Y) \le C_X C_Y \delta_p(h,X_0,Y_0) 
\]
for all
\[
h = \sum_{i=1}^n f_i \otimes g_i\, , \qquad f_i \in X_0\, , \quad g_i \in Y_0\, , 
\quad i=1, \ldots \, ,n\, .
\]
Hence $X_0 \otimes_{\delta_p} Y_0 \hookrightarrow X \otimes_{\delta_p} Y$.
This will be applied with 
$X_0 = Y_0 = \HH_{\alpha,p}$ and  $X= Y= H^\alpha_p ([0,1])$,
see Theorem~\ref{tensor3}, 
Theorem~\ref{Winfried}, 
and Proposition~\ref{tensor4}.
As a result we have proved 
\[
\HH_{\alpha,2,p} \hookrightarrow S^{\alpha}_{p}H([0,1]^2). 
 \]
Now we proceed with induction on $s$. 
\end{proof}

\subsection*{Acknowledgment}
The authors thank Karsten Eggers, Goran Nakerst, Hans Triebel, and Marcin Wnuk for valuable comments. 

The first and the fourth author would like to thank the Isaac Newton Institute for Mathematical Sciences (INI), Cambridge, for support and hospitality during the programme Multivariate approximation, discretization, and sampling recovery, where work on this paper was undertaken. This work was supported by EPSRC grant EP/R014604/1, ARC grant DP220101811 and ELES. 

\footnotesize{

}

\end{document}

We point out that in the case $0<\alpha<1/p$ the Haar wavelet spaces coincide with the spaces used in \cite{Tri10}, \cite{Hin10} (dyadic versions) and \cite{Mar13} (see Theorem 2.11) where they are 
called Besov spaces with dominating mixed smoothness. The bases from the aforementioned works are not the same ones as used here. 
For instance the bases from \cite{Mar13} are orthogonal. To realise the equivalence of the norms of our Haar wavelet spaces and the ones from \cite{Mar13} 
we point out that the support of the functions are identical, namely elementary intervals $E_k^j$. One can characterize the wavelet functions with piecewise constant 
(on elementary intervals) functions and therefore the equivalence follows.

\begin{lemma}\label{tensor0}
Let $1 <p<\infty$.
Let
\[
K:= \Big\{ \sum_{j=1}^n \eta_j \, \cx_j : \quad n \in \N, \: \eta_j \in \C\, , \: |\supp \cx_j \cap \supp \cx_i|=0\, , \: i \neq j\Big\} \, .
\]
Then 
\[
\Big\{f : \quad \exists c\in \re, \: \exists \tilde{f} \in K \quad \mbox{s.t.}\quad  f= c + \int_0^1 \tilde{f} (t)\, (x-t)_+^{\alpha-1}\, dt\Big\} 
\]
is a dense subset in $\HH_{\alpha,p}$.
\end{lemma}